\documentclass{elsart}
\usepackage{amssymb, amsmath, amsfonts}
\usepackage[mathscr]{euscript}
\usepackage{mathrsfs}
\usepackage[english]{babel}
\usepackage{bm}
\usepackage[all]{xy}

\font\tenDDl=msbm10  
\font\sevenDDl=msbm7 
\font\fiveDDl=msbm5 
        \newfam\DDlfam  
                \textfont\DDlfam=\tenDDl
\scriptfont\DDlfam=\sevenDDl
        \scriptscriptfont\DDlfam=\fiveDDl
\def\CE{{\bf C}}

\long\def\nodo#1{{}}

\def\genfd{{\bm k}}
\def\genrg{{\bm k}} 

\def\functcat{{\underline{\rm Funct}}}
\def\setscat{{\rm Sets}}
\def\cgen{\bar{C}}
\def\MR#1{{ MR#1}} 
\def\Tr{\operatorname{Tr}}
\def\ovM{\overline{M}}
\def\ovb{\overline{b}}
\def\nxpoint{\refstepcounter{subsection}%
  \makepoint{\thesubsection}}
\def\nxsubpoint{\refstepcounter{subsubsection}%
  \makepoint{\thesubsubsection}}

\let\refpt=\refpoint
\let\myunderline=\textbf
\def\id{\operatorname{id}}
\def\ad{\operatorname{ad}}
\def\gr{\operatorname{gr}}
\def\pr{\operatorname{pr}}
\def\Ob{\operatorname{Ob}}
\def\Hom{\operatorname{Hom}}
\def\End{\operatorname{End}}
\def\Pro{\operatorname{Pro}}

\def\Der{\operatorname{Der}}
\def\Ker{\operatorname{Ker}}

\def\Coder{\operatorname{Coder}}
\def\Dercont{\operatorname{Der.cont}}
\def\iHom{\operatorname{\myunderline{Hom}}}
\def\iEnd{\operatorname{\myunderline{End}}}
\def\iIsom{\operatorname{\myunderline{Isom}}}
\def\iAut{\operatorname{\myunderline{Aut}}}
\def\Spf{\operatorname{Spf}}
\def\Spec{\operatorname{Spec}}
\def\uVect{\operatorname{\myunderline{Vect}}}
\def\calE{\mathcal{E}}
\def\calP{\mathcal{P}}
\def\calD{\mathscr{D}}
\def\calL{\mathscr{L}}
\def\Uenv{\mathcal{U}}
\def\gg{\mathfrak{g}}
\def\bbN{\mathbb{N}}
\def\bbQ{\mathbb{Q}}
\def\bbZ{\mathbb{Z}}
\def\bbC{\mathbb{C}}
\def\bbR{\mathbb{R}}
\def\bbG{\mathbb{G}}
\def\counit{{\bm\eta}}
\def\Seq{\mbox{\sl Seq}}
\def\itilde{i'}
\def\Ahat{{\hat A}}
\def\Shat{{\hat S}}
\def\uu#1{{\bm #1}}
\def\quotedprojlim{\operatorname{\mbox{``$\varprojlim$''}}}
\let\phi=\varphi
\let\epsilon=\varepsilon
\let\emptyset=\varnothing
\let\injlim=\varinjlim
\let\projlim=\varprojlim
\def\makepoint#1{\medbreak\noindent{\bf #1. }}

\def\relPk{{}_{\genrg\backslash}{\mathcal P}}
\def\relPkp{{}_{\genrg'\backslash}{\mathcal P}}
\def\functcat{\operatorname{\myunderline{Funct}}}
\def\setscat{{\rm Sets}}
\def\uuO{\mathscr O}  
\def\uW{{\myunderline{W}}}    
\def\uWom{\uW^\omega}
\def\Lie{\operatorname{Lie}}
\def\uLie{\operatorname{\myunderline{Lie}}}
\def\Expp{\operatorname{\myunderline{Exp}_+}}
\def\Expt{\operatorname{\myunderline{Exp}_\times}}
\def\Vect{\operatorname{Vect}}
\begin{document}\begin{frontmatter}

\title{A universal formula for representing Lie algebra generators as
formal power series with coefficients in the Weyl algebra}

\author{Nikolai Durov}
\address{Max Planck Institut f\"ur Mathematik, 
P.O.Box 7280, D-53072 Bonn, Germany\\
and Department of Mathematics and Mechanics, St. Petersburg 
State University, 198504 St. Petersburg, Russia}
\ead{durov@mpim-bonn.mpg.de}
\author{Stjepan Meljanac, Andjelo Samsarov, Zoran \v{S}koda}
\address{Theoretical Physics Division,
Institute Rudjer Bo\v{s}kovi\'{c}, Bijeni\v{c}ka cesta~54, P.O.Box 180,
HR-10002 Zagreb, Croatia
}
\ead{meljanac@irb.hr, asamsarov@irb.hr, zskoda@irb.hr}

\date{{}}

\maketitle
\begin{abstract}
\noindent Given a $n$-dimensional Lie algebra $\gg$
over a field $\genfd \supset \mathbb Q$,
together with its vector space basis $X^0_1,\ldots, X^0_n$, 
we give a formula, depending only on the structure constants,
representing the infinitesimal generators, 
$X_i =  X^0_i t$ in $\gg\otimes_\genfd \genfd [[t]]$, 
where $t$ is a formal variable, as a formal
power series in $t$ with coefficients in the Weyl algebra $A_{n}$. 
Actually, the theorem is proved for Lie algebras
over arbitrary rings $\genfd\supset \mathbb Q$.

We provide three different proofs, each of which is expected to be useful
for generalizations. The first proof is obtained by direct calculations 
with tensors. This involves a number of interesting combinatorial formulas 
in structure constants. The final step in calculation
is a new formula involving Bernoulli numbers and arbitrary derivatives
of $\coth(x/2)$. The dimensions of certain spaces 
of tensors are also calculated. The second method of proof
is geometric and reduces to a calculation of formal right-invariant
vector fields in specific coordinates, 
in a (new) variant of formal group scheme theory. 
The third proof uses coderivations and Hopf algebras.
\end{abstract}
\end{frontmatter}
\section{Introduction and the statement of the main theorem}

We consider here a remarkable special case of the following problem:
given a commutative ring $\genfd$, 
when may a given associative $\genfd$-algebra $U$ with $n$ generators,
say $X_1^0,\ldots, X_n^0$, 
be represented as a formal 1-parameter deformation
of a commutative polynomial subalgebra $\genfd[x_1,\ldots,x_n]$ of 
the Weyl algebra $A_{n,\genfd} := 
\genfd[x_1,\ldots,x_n,\partial^1,\ldots,\partial^n]/ 
\langle\partial^j x_i - x_i \partial^j - \delta^j_i \rangle$,
where the whole deformation is within the Weyl algebra itself.
More explicitly, we look for the deformations 
of the form
$$X_i^0 t = x_i + \sum_{N = 1}^\infty P_{N,i} t^N, \,\,\,\,\,\,\,\,\,\,\,\,
P_{N,i} \in  A_{n,\genfd},$$
where $t$ is a deformation parameter.
Because the deformation can introduce only infinitesimal noncommutativity,
we rescaled the generators $X_i^0$ by factor $t$ in the very formulation
of the problem, i.e. the actual algebra realized as
a deformation is $U_t$, what is the positive degree 
$U_t = \oplus_{i > 0} U t^i$ part of 
$U\otimes_\genfd \genfd[[t]]$. 
Such a deformation, if it exists, does not need to be unique.
In this paper we find a universal formula
which provides such a deformation when $U$ is the enveloping
algebra of any Lie algebra $\gg$ over any unital ring $\genfd$ 
containing the field ${\mathbb Q}$ of rational numbers. 
Still, the underlying $\genfd$-module of $\gg$ will be assumed free
(only projective when formula given in invariant form) 
and finitely generated. 

Given a basis $X^0_1,\ldots, X^0_n$ of a free finite-rank $\genfd$-module
underlying a Lie algebra $\gg$,
the structure constants $(C^0)_{ij}^k$ of $\gg$ are defined by 
$[X^0_i, X^0_j] =  \sum_{k = 1}^n (C^0)_{ij}^k X_k$ 
and are clearly antisymmetric in the {\it lower} two indices. 
As usual, the choice of the basis will be considered 
as an isomorphism $X^0 : \genfd^n\to \gg$
given by $X^0(e_i) = X^0_i$, where $e_1,\ldots, e_n$ 
is the standard basis of $\genfd^n$. 

Let $\gg\otimes_\genfd \genfd[[t]] = \oplus_{i = 0}^\infty \gg t^i$
be the Lie algebra $\gg$ but with scalars extended to 
include formal power series in one variable. 
Its positive degree part $\gg_t := \oplus_{i = 1}^\infty \gg t^i$ is
a Lie subalgebra of $\gg\otimes_\genfd \genfd[[t]]$ over $\genfd[[t]]$
with basis $X_1,\ldots, X_n$ where $X_i = X^0_i t$. Then
$[X_i, X_j] = C_{ij}^k X_k$ where the new structure
constants $C_{ij}^k := (C^0)^k_{ij} t$ are also of degree 1 in $t$ and
may be interpreted as infinitesimal. Denote by ${\mathcal U}(\gg_t)$
the universal enveloping $\genfd[[t]]$-algebra of $\gg_t$. 
It naturally embeds into $({\mathcal U}(\gg))_t :=
\oplus_{i > 0} {\mathcal U}(\gg) t^i$. 

Define a matrix $\CE$ over $A_{n,\genfd}$ by
$\CE^i_j = \sum_{k = 1}^n C^i_{jk}\partial^k$.

{\bf Main theorem. 
In above notation, if the structure constants are
totally antisymmetric, then for any number $\lambda \in \genfd$,
the formula
$$
X_i \mapsto \sum_\alpha x_\alpha \phi^\alpha_i
$$
where
\begin{equation}\label{eq:phi}
\phi^\alpha_\beta := \sum_{N = 0}^\infty 
\frac{(-1)^N B_{N}}{N!} (\CE^N)^\alpha_\beta
\in\, \genfd[\partial^1,\ldots,\partial^n][[t]]\hookrightarrow 
A_{n,\genfd}[[t]],
\end{equation}
and $B_n$ are Bernoulli numbers, 
extends to an embedding of associative $\genfd[[t]]$-algebras
$\Phi_\lambda : {\mathcal U}(\gg_t)\hookrightarrow A_{n,\genfd}[[t]]$.
If $\genfd = {\mathbb C}$ (or ${\mathbb Q}[\sqrt{-1}]$),
and if the basis is chosen such that $C^i_{jk} \in {\mathbb R}\sqrt{-1}$,
and $\lambda \in {\mathbb R}$, then
the same holds for the more general formula
}
$$
X_i \mapsto \sum_\alpha \lambda x_\alpha \phi^\alpha_i + 
(1-\lambda) \phi^\alpha_i x_\alpha.
$$
{\it Note that, for $\lambda = \frac{1}{2}$, 
the expressions for $X_i$ are hermitean (invariant with 
respect to the standard antilinear involution on $A_{n,{\mathbb C}}$).
}
In particular, for $\lambda = 1$
$$
X_i \mapsto x_i + \frac{1}{2} C^k_{ij} x_k\partial^j + \frac{1}{12}
C^{k'}_{ij} C^{k}_{k' j'} x_{k} \partial^j \partial^{j'} 
- \frac{1}{720} C^{k'}_{ij} C^{k''}_{k' j'} C^{k}_{k'' j''} x_{k} 
\partial^j \partial^{j'}\partial^{j''}
+\ldots
$$

{\it Clearly, the image 
$\Phi_\lambda({\mathcal U}(\gg_t))/ (t\Phi_\lambda({\mathcal U}(\gg_t)))$
modulo the subspace of all elements of degree $2$ and higher in $t$ is 
the polynomial algebra in ${\rm dim}_\genfd\,\gg$ {\it commuting} variables
$x_i = \Phi_\lambda(X_i)$. 
Thus this embedding may be considered 
as a realization of ${\mathcal U}(\gg_t)$ as
a {\bf deformation} of the commutative algebra $\genfd[x_1,\ldots,x_n]$,
where the whole deformation is taking place
within the Weyl algebra $A_{n,\genfd}$. 
}

In Section~\ref{sec:lambdared} we show that the generalization to
general $\lambda$ (when $\genfd = \bbC$) is easy.
In Sections 3--6 three of us (S.M., A.S., Z.\v{S}.)
motivate and prove the theorem
by direct computation with tensors.
In Sections \ref{s:Formal}--\ref{s:EndOfProof}
the first author (N.D.) gives an alternative proof
and interpretation using formal geometry.
In this second part, a completed Weyl algebra is used
instead of working with a deformation parameter~$t$
to make sense of the power series expressions in our formulas.
Over an arbitrary ring $\genfd\supset \bbQ$, the completed
Weyl algebra is identified with the algebra of
formal vector fields on a formal neighborhood of the origin
in our Lie algebra~$\gg$, considered here as a formal variety.
Similarly to the classical Lie theory over $\genfd = \bbC$ or $\bbR$,
where the elements of a Lie algebra can be interpreted as
(say, right) invariant vector fields on a Lie group,
we identify the elements of
a Lie algebra~$\gg$ over any ring $\genfd\supset\bbQ$
to the right-invariant formal vector fields
on a suitable formal group,
and compute them in terms of the coordinate chart given by
an appropriate version of the exponential map.
As a consequence, in this part (Sections \ref{s:Formal}--\ref{s:EndOfProof}),
we actually construct a deformation of the
abelian subalgebra generated by the~$\partial^i$ within the Weyl algebra,
rather than the subalgebra of coordinates,
but the difference is inessential: the automorphism of the Weyl algebra
mapping $x_i\mapsto -\partial^i$, $\partial^i \mapsto x_i$
interchanges the formulas between the first and the second parts of the work.
However, we kept the different conventions as the deformation
of ``space'' coordinates $x_i$ is our initial motivation,
while the representation via vector fields is
also a valuable geometric point of view.
In Section~\ref{sec:thirdproof},
N.D.\ adds a third proof using coalgebra structure
and coderivations. In some sense this proof is obtained by ``dualizing''
the previous geometric proof; this makes the proof shorter but more difficult
to understand.

As we learned from {\sc D. Svrtan} after completing
our first proof, one can find a superimposable
formula in {\sc E. Petracci}'s work~\cite{Petracci}
on representations by coderivations.
More precisely, Theorem 5.3 and formula (20) of her work,
once her formulas (13) and (15) are taken into account,
correspond to our~\refpt{p:mf.inv.coder}, i.e.\ the invariant form
of the Main Formula of present work,
expressed in the language of coderivations,
prior to any Weyl algebra identifications
(Weyl algebras are in fact never mentioned in \cite{Petracci})
and explicit coordinate computations.
Moreover, our formulas of~\refpt{p:adxpairing} essentially appear in
\cite{Petracci}, Remark~3.4.

{\it Notation.} Throughout the paper, for $\rho = 1,\ldots, n$, 
we use the $\genfd$-linear derivations 
$\delta_\rho := \frac{d}{d(\partial^\rho)}$ of $A_{n,\genfd}$. 
If $r$ is a real number, then $\lfloor r \rfloor$ 
denotes the largest integer smaller or equal to $r$
(integer part or floor of $r$). 

All considered modules over unital rings will be unital. 
For a $\genfd$-module $M$, $S(M), T(M)$ will denote its symmetric
and tensor $\genfd$-algebra, respectively.

\section{Reduction to $\lambda = 1$}
\label{sec:lambdared}

\nxpoint \label{p:gencommxfi}
Let $\psi$ and $\chi$ be matrices of expressions 
depending on $\partial$-s. Assume Einstein convention (summation over
each pair of repeated indices). Then
$$\begin{array}{lcl}
\lbrack x_\alpha \psi^\alpha_\mu, x_\beta \chi^\beta_\nu \rbrack & = &
x_\gamma ((\delta_\rho \psi^\gamma_\mu) \chi^\rho_\nu 
- (\delta_\rho \chi^\gamma_\nu) \psi^\rho_\mu ) \\
\lbrack x_\alpha \psi^\alpha_\mu,\chi_\nu^\beta x_\beta  \rbrack & = &
 x_\gamma ((\delta_\rho\psi^\gamma_\mu)\chi^\rho_\nu
- (\delta_\rho \chi^\gamma_\nu)\psi^\rho_\mu)
-\psi^\alpha_\mu(\delta_\alpha\delta_\beta \chi^\beta_\nu)\\
\lbrack \psi^\alpha_\mu x_\alpha , \chi_\nu^\beta x_\beta \rbrack & = &
x_\gamma((\delta_\rho \chi^\rho_\nu)\psi^\gamma_\mu - 
(\delta_\rho \psi^\rho_\nu)\chi^\gamma_\nu)
\\
\lbrack \psi^\alpha_\mu x_\alpha , x_\beta \chi^\beta_\nu \rbrack & = &
x_\gamma((\delta_\rho \chi^\rho_\nu)\psi^\gamma_\mu - 
(\delta_\rho \psi^\rho_\nu)\chi^\gamma_\nu)
+ \psi^\alpha_\mu(\delta_\alpha\delta_\beta \chi^\beta_\nu)
\end{array} $$
\nxpoint There is a $\genfd$-linear antiautomorphism $\dagger$ 
of the Weyl algebra $A_{n,\genfd}$
given on generators by $x \mapsto x$ and $\partial \mapsto -\partial$.
In particular, $(x_\alpha \phi^\alpha_\mu)^\dagger = 
(\phi^\alpha_\mu)^\dagger x_\alpha$.

In the case when $\genfd = {\mathbb C}$, we also have the 
conjugation -- the antilinear involution which will be also
denoted by $\dagger$. In that case, it is easy to check 
that for $\phi$ given by formula~\eqref{eq:phi} 
we have $\phi = \phi^\dagger$.

\nxpoint From the formulas in~\refpt{p:gencommxfi}, 
or, even easier, directly, we obtain
\begin{equation}\label{eq:refl2plus2}
[ x_\alpha \phi_{\alpha\mu}, x_\beta \phi_{\beta\nu}]
+ [ \phi_{\alpha\mu}x_\alpha, \phi_{\beta\nu} x_\beta ] 
= [ x_\alpha \phi_{\alpha\mu},  \phi_{\beta\nu} x_\beta]
+ [ \phi_{\alpha\mu}x_\alpha,  x_\beta \phi_{\beta\nu}]
\end{equation}
Therefore 
$$\begin{array}{l}
[\lambda x_\alpha \phi^{\alpha}_{\mu} +
(1-\lambda)(x_\alpha \phi^{\alpha}_\mu)^\dagger,
\lambda x_\beta \phi^\beta_\nu +
(1-\lambda)(x_\beta \phi_\nu^\beta)^\dagger ]\,\,\,\,= \\ 
\,\,\,\,\,\,\,\,\,\,\,\,\,\,\,\,\,\,\, = \,
[ \lambda x_\alpha \phi^\alpha_\mu + (1-\lambda)\phi^\alpha_\mu x_\alpha,
\lambda x_\beta \phi^\beta_\nu + (1 - \lambda) \phi^\beta_\nu x_\beta ]\\
\,\,\,\,\,\,\,\,\,\,\,\,\,\,\,\,\,\,\, \stackrel{\eqref{eq:refl2plus2}}= \, 
\lambda (\lambda + 1-\lambda) 
[ x_\alpha \phi^\alpha_\mu, x_\beta \phi^\beta_\nu]\, + \\
\,\,\,\,\,\,\,\,\,\,\,\,\,\,\,\,\,\,\,\,\,\,\,\,\,\,\,\,\,\,\,\,
 +\, (1-\lambda) (1-\lambda +\lambda) 
[ \phi^\alpha_\mu x_\alpha, \phi^\beta_\nu x_\beta]\\
\,\,\,\,\,\,\,\,\,\,\,\,\,\,\,\,\,\,\, = \, 
\lambda [ x_\alpha \phi^\alpha_\mu, x_\beta \phi^\beta_\nu]
+ (1-\lambda)([x_\beta \phi^\beta_\nu,x_\alpha \phi^\alpha_\mu])^\dagger
\end{array}$$
Thus it is sufficient to prove the $\lambda = 1$ identity
\begin{equation}\label{eq:C1}
[x_\alpha \phi^\alpha_\mu, x_\beta \phi^\beta_\nu ] = 
C_{\mu\nu}^\rho x_\gamma \phi^\gamma_\rho
\end{equation}
and the general identity is 
$\lambda\mbox{Eq.}(\ref{eq:C1})+(1-\lambda)\mbox{Eq.}(\ref{eq:C1})^\dagger$. 
\section{Covariance and the universal case} 
\label{sec:covariance}

\nxpoint In calculational approach, we first try a more general Ansatz, 
and then gradually inspect various identities a formula 
should satisfy in each order in the deformation parameter in order to
provide a Lie algebra representation. After long
calculations we conclude that specializing the coefficients to 
Bernoulli numbers ensures that the identities hold.

Our more general Ansatz, with nice covariance
properties under $GL(n,\genfd)$-action, is a special case of
a more general requirement of functoriality.
We want that our formula be universal under change of rings, and
universal for all Lie algebras over a fixed ring. 
Thus we build the Ansatz from tensors in structure constants 
of certain kind, and we want the shape to be controlled 
(covariant in some sense) under the morphisms of Lie algebras, 
where we allow the underlying ring, the Lie algebra, 
and its basis to change. To this end
we will define certain universal ring, which is not a field, 
and later a universal Lie algebra over 
it where our calculations in fact take place. 
By specialization, the formulas then imply the formulas
for ``concrete'' Lie algebras. 

\nxpoint {\bf Definition.} 
{\it Consider the affine $\genfd$-space $\genfd^{n^3}$. 
Define the affine function algebra $\genfd[{\mathcal C}_n]$ of 
the affine variety ${\mathcal C}_n$ 
} 
(``the variety of generic structure constants 
of generic rank $n$ Lie $\genfd$-algebra'') 
{\it 
to be the polynomial algebra in $n^3$-variables
$\cgen^i_{jk}$, $i,j,k = 1,\ldots, n$, 
modulo the homogeneous relations 
\newline (i) $\cgen^i_{jk} = \cgen^i_{kj}$ 
(antisymmetry in lower indices) 
\newline (ii) $\sum_\alpha \cgen_{ij}^\alpha \cgen_{\alpha k}^l + 
\cgen_{jk}^\alpha \cgen_{\alpha i}^{l} + 
\cgen_{ki}^{\alpha}\cgen_{\alpha j }^{l} = 0\,\, \forall i,j,k,l$
(Jacobi identity).

Let ${\mathcal L}_n$ be the Lie $\genfd[{\mathcal C}_n]$-algebra 
over $\genfd[{\mathcal C}_n]$, free as a $\genfd[{\mathcal C}_n]$-module
with basis 
$X : (\genfd[{\mathcal C}_n])^n\stackrel{\cong}\rightarrow{\mathcal L}_n$ 
and bracket $[X_k, X_l] = \cgen^i_{kl} X_i$. 
It will be sometimes also called universal.
}

\nxpoint The correspondence which to each $\genfd$-algebra
associates the set of all Lie brackets on 
the free $\genfd$-module $\genfd^n$ of rank $n$ extends to a covariant
functor from the category of (unital associative) 
$\genfd$-algebras to the category of sets.
It is clearly represented by $\genfd[{\mathcal C}_n]$.
If we take an arbitrary rank $n$ free Lie algebra $\gg$ over $\genfd$ as
a Lie algebra, and fix a basis $X^0 = (X^0_1,\ldots,X^0_n) : \genfd^n \to \gg$,
it can be therefore considered as a point of affine $\genfd$-variety 
${\mathcal C}_n$. 
The map of $\genfd$-algebras 
${\rm ev}_\gg :=  {\rm ev}_{\gg, X^0} : \genfd[{\mathcal C}_n]\to \genfd$ 
determined by $\bar{C}^i_{jk} \mapsto C^i_{jk}$ 
is called the {\bf evaluation map}.  

\nxpoint Let $A_{n,\genfd} [[t]]$ be the $\genfd[[t]]$-algebra of 
formal power series in one indeterminate
$t$ with coefficients in $A_{n,\genfd}$. 
Any choice of a basis in $\gg$ provides an isomorphism 
of $\genfd$-modules from 
$S(\gg)\otimes S(\gg^*)$ to $A_{n,\genfd}$, 
where $S(\gg)$ is the symmetric (polynomial) algebra in
$X_1^0,\ldots, X_n^0$. Algebra $S(\gg)\otimes S(\gg^*)$ acts on the left and 
right on $S(\gg)$, namely the elements of $S(\gg)$ act by multiplication,
and the elements of $\gg^*$ act by derivations. 

\nxpoint $GL_\genfd(\gg) \cong GL(n,\genfd)$ naturally acts on 
$\gg$, $\gg^*$, $T(\gg)\otimes T(\gg^*)$ and $S(\gg)\otimes S(\gg^*)$. 
We will take our Lie algebra $\gg$ to be free as $\genfd$-module
to be able to work with tensor components.
As our main interest is in the formulas for generators, 
given ${\mathcal O}\in GL_\genfd(\gg)$, we find more convenient to
consider the matrix elements ${\mathcal O}^\alpha_i$
for the expansion of a new basis ${\mathcal O}X$
in terms of old $X$, rather than
the more customary matrix elements for the
expansion of the contragradient vector components:
$({\mathcal O}X)_i =: \sum_{\alpha} {\mathcal O}^\alpha_i  X_\alpha$
The  structure constants $C^{{(\mathcal O})i}_{jk}$ in the new frame 
$({\mathcal O} X)_1,\ldots,({\mathcal O} X)_n$ can be easily described
\begin{equation}\label{eq:Cequiv}
\left\lbrack ({\mathcal O} X)_i, ({\mathcal O} X)_j\right\rbrack = 
  C^{({\mathcal O})\sigma}_{ij} ({\mathcal O}X)_\sigma,
\end{equation}
where $C^{({\mathcal O})\sigma}_{ij} := 
{\mathcal O}^\alpha_i {\mathcal O}^\beta_j 
C^\gamma_{\alpha\beta} ({\mathcal O}^{-1})^\sigma_\gamma$
(clearly the structure constants make a tensor which may 
be considered as living in $\gg^*\otimes \gg^*\otimes \gg$, but we here
present everything in coordinates).

\nxpoint The natural $GL(\gg,\genfd)$-action on $S(\gg)\otimes S(\gg^*)$
transports to the Weyl algebra $A_{n,\genfd}$ via the
identification of their underlying $\genfd$-modules.
Both actions may be considered as factored from $T(\gg)\otimes T(\gg^*)$.
It is crucial that the induced action is compatible 
both with the product in $A_{n,\genfd}$ and 
with the product in $S(\gg)\otimes S(\gg^*)$,
in the sense that 
${\mathcal O}(x\cdot y) = {\mathcal O}(x){\mathcal O}(y)$ 
for any $x,y \in A_{n,\genfd}$, 
${\mathcal O}\in GL(\gg,\genfd)\cong GL_n(\genfd)$.
This is because the $GL(\gg,\genfd)$-action
is factored from the action on $T(\gg)\otimes T(\gg^*)$.
Namely, the defining ideal, both for $S(\gg)\otimes S(\gg^*)$ and 
for $A_{n,\genfd}$, is $GL_n(\genfd)$-invariant 
(in the case of $A_{n,\genfd}$,
$\partial^j x_i - x_i \partial^j -\delta_i^j$ are components
of a tensor for which all components are included in the ideal). 
For this $GL_n(\genfd)$-action, $x = (x_1,\ldots,x_n)$ 
are $GL_n(\genfd)$-cogredient with respect to $X$,
and $\partial = (\partial^1,\ldots,\partial^n)$ are contregredient, 
i.e. $({\mathcal O}x)_i  = {\mathcal O}^\alpha_i x_\alpha$
and $({\mathcal O}\partial)^i =  ({\mathcal O}^{-1})^i_\alpha 
\partial^\alpha$.

Contractions of tensors with respect to the product in $A_{n,\genfd}$
have the same covariance, as if they would be contractions 
with respect to the product in $S(\gg)\otimes S(\gg^*)$. 
Given a contraction, $c : (\gg^*)^m\otimes \gg^n \to (\gg^*)^{m-1}\otimes 
\gg^{n-1}$ (e.g. the pairing $\gg^*\otimes \gg\to \genfd$), 
and tensors $A \in \gg^*$, $B\in \gg$ one usually considers the
behaviour or $c(A,B)$ under the action 
$c(A,B) \mapsto c({\mathcal O}A,{\mathcal O}B)$ and likewise for
multiple contractions. 
Writing down $\CE$ in terms of $C^i_{jk}$ and $\partial^l$
and using induction it is direct to show that  
$$\begin{array}{l}
({\mathcal O}\CE^N)^\sigma_\rho = 
\sum_{i,j} {\mathcal O}^\sigma_i  (\CE^N)^i_j ({\mathcal O}^{-1})^j_\rho \\
{\mathcal O}(\sum_\alpha A_\alpha x_\alpha (\CE^N)^\alpha_\beta)
= \sum_\sigma {\mathcal O}^\sigma_\beta  
\sum_\alpha A_\alpha x_\alpha (\CE^N)^\alpha_\sigma,
\end{array}$$
{\it Thus looking for the solution in terms of a series
$\sum_\alpha A_\alpha x_\alpha (\CE^N)^\alpha_\beta$ 
is seeking for a solution which behaves covariantly 
with respect to the change of coordinates under $GL_n(\genfd)$.
}

\nxpoint In the case of the universal rank $n$ Lie algebra ${\mathcal L}_n$
over $\genfd[{\mathcal C}_n]$, the components of the structure tensor 
are identical to the generators of the ground ring. Hence ${\mathcal O}$
induce an automorphism of $\genfd[{\mathcal C}_n]$ as a $\genfd$-module.
It is a remarkable fact, however, that $C^{({\mathcal O})i}_{jk}$
satisfy the same relations as $C^i_{jk}$, i.e. ${\mathcal O}$ induces
an $\genfd$-{\it algebra} automorphism of $\genfd[{\mathcal C}_n]$. 
For example, if we want to show the Jacobi identity, we consider
$$
{\mathcal O} (C^\alpha_{ij} C^l_{\alpha k})
= {\mathcal O}^{r}_i {\mathcal O}^s_j C^\sigma_{rs} 
({\mathcal O}^{-1})^\alpha_\sigma {\mathcal O}^t_\alpha 
{\mathcal O}^u_k C^\tau_{tu} ({\mathcal O}^{-1})^l_\tau
= {\mathcal O}^{r}_i {\mathcal O}^s_j
{\mathcal O}^u_k  C^t_{rs} C^\tau_{tu} ({\mathcal O}^{-1})^l_\tau
$$
and, two other summands, cyclically in $(i,j,k)$. 
Now rename $(r,s,t)$ apropriately in the two other summands to force the same
factor ${\mathcal O}^r_i {\mathcal O}^s_j {\mathcal O}^u_k$ in all
three summands. 
Then the $CC$-part falls in a form where the Jacobi identity can be readily
applied. The covariance follows from the functoriality under 
the canonical isomorphism of $\genfd[{\mathcal C}_n]$-Lie algebras 
from the pullback ${\mathcal O}^* {\mathcal L}_n$ to ${\mathcal L}_n$
for all ${\mathcal O}$ in $GL_n(\genfd)$.

\nxpoint 
A more general covariant Ansatz 
$X_\beta = \sum_{IJK} A_{IJK} x_\alpha (\Tr \CE^{2I})^J 
(\CE^{K})^\alpha_\beta$ is likely useful for finding new representations
for specific Lie algebras. 
In the universal case ($\gg = {\mathcal L}_n$),
however, the traces are contraction-disconnected 
from the rest of expression, and different trace factors 
can not be mixed. Namely, the defining ideal of $\genfd[{\mathcal C}_n]$ 
does not have ``mixed'' (different $J$ and $K$) elements, and commuting
$x_\alpha$ with such tensors also does not produce them either. 

\section{Differential equation and recursive relations}

Introduce the ``star'' notation for higher order 
connected tensors. Namely, in order to spot better just the relevant
indices, if we have the contraction of an upper index 
in one of the $C$-factors with one lower index of the next $C$-factor,
we may just write the $*$ symbol on the two places. 
In that notation, the derivative 
$\delta_\rho = \frac{\partial}{\partial (\partial^\rho)}$ 
applied to $\CE^I$ equals
$$\delta_\rho (\CE^I)^\gamma_\mu = 
C^{*}_{\mu\rho} \CE^*_* \CE^*_* \cdots \CE^\gamma_*
+ \CE^{*}_{\mu} C^*_{*\rho} \CE^*_* \cdots \CE^\gamma_*
+\ldots + \CE^{*}_{\mu} \CE^*_{ *} \CE^*_* \cdots C^\gamma_{*\rho}
$$

From 
$[x_\alpha \phi^\alpha_\mu, x_\beta \phi^\beta_\nu ] = C^\sigma_{\mu\nu}
x_\gamma \phi^\gamma_\rho$ (C1) (see Section~\ref{sec:lambdared}), 
equating the polynomials in $\partial$-s 
in front of $x_\gamma$, we obtain the system of 
differential equations manifestly antisymmetric with respect to the
interchange $\mu\leftrightarrow\nu$:
$$\fbox{$
(\delta_\rho \phi^\gamma_\mu)\phi^\rho_\nu -
(\delta_\rho \phi^\gamma_\nu)\phi^\rho_\mu =
C^\sigma_{\mu\nu} \phi^\gamma_\sigma
$}$$
We use Ansatz $\phi^i_j = \sum_{N = 0}^\infty A_I (\CE^I)^i_j$
where $A_0 = 1$ and $(\CE^0)^i_j = \delta^i_j$. 
In any order $N\geq 1$ in expansion in $t$ (hence in $C^i_{jk}$-s and also 
in $\partial^i$-s) we thus obtain
\begin{equation}\label{eq:diffphi}
\sum_{I=1}^N A_I A_{N-I} \left\lbrace [\delta_\rho (\CE^I)^\gamma_\mu]
(\CE^{N-I})^\rho_\nu - [\delta_\rho (\CE^I)^\gamma_\nu]
(\CE^{N-I})^\rho_\mu\right\rbrace 
= A_{N-1} C^\sigma_{\mu\nu}(\CE^{N-1})^\gamma_\sigma
\end{equation}
Notice that for $N = I$,
$$\begin{array}{lcl}
\delta_\rho (\CE^I)^\gamma_\mu 
(\CE^{N-I})^\rho_\nu & = &\delta_\nu (\CE^N)^\gamma_\mu \\
& = & C^{*}_{\mu\nu} \CE^*_* \cdots \CE^\gamma_*
+ \CE^{*}_{\mu} C^*_{*\nu} \CE^*_* \cdots \CE^\gamma_*
+\ldots + \CE^{*}_{\mu} \CE^*_* \cdots C^\gamma_{*\nu}\\
& = : & M'_0 + M'_1 + \ldots + M'_{N-1}
\end{array}$$
where we use the star notation, as explained above.
Of course, the tensors $M'_I = (M'_I)^\gamma_{\mu\nu}$ 
have three {\it supressed} indices $\mu,\nu,\gamma$.
In equation~(\ref{eq:diffphi}) these tensors 
come in the combination 
$(M_I)^\gamma_{\mu\nu} := (M'_I)^\gamma_{\mu\nu} -  (M'_I)^\gamma_{\nu\mu}$
antisymmetrized in the lower two indices. In detail,
$$\begin{array}{lcl}
M_0 &=& (M_0)^\gamma_{\mu\nu}:= 
2 C^{*}_{\mu\nu} \CE^*_* \CE^*_* \cdots \CE^\gamma_*, \\
M_I &=& (M_I)^\gamma_{\mu\nu}:= 
(\CE^{I-1})^{*}_{\mu} C^*_{*\nu} (\CE^{N-I-1})^\gamma_*
- (\mu\leftrightarrow\nu),
\,\,\,\,1\leq I \leq N-1,\\
M_{N-1} & = & (M_{N-1})^\gamma_{\mu\nu}:= 
\CE^*_\mu (\CE^{N-2})^*_* \CE^\gamma_\nu 
- \CE^*_\nu (\CE^{N-2})^*_* \CE^\gamma_\mu.
\end{array}$$
One can easily show that if the Lie algebra in question is 
$\mathfrak{su}(2)$ then $M_0 = M_1 = \ldots = M_{N-1}$ 
and the rest of the proof is much simpler. 
We checked using {\it Mathematica}$^\textsc{ TM}$ that 
already for  $\mathfrak{su}(3)$ the $M_I$-s are mutually different.  

In particular, the $N=1$ equation is obvious: 
$A_1 (C^\gamma_{\mu\nu} - C^\gamma_{\nu\mu}) = C^\gamma_{\mu\nu}$
what by antisymmetry in $(\mu\leftrightarrow \nu)$
forces $A_1 = \frac{1}{2}$.

The Jacobi identity $[[X_i, X_j], X_k] + cyclic = 0$ will be
used in the form
$$
\sum_\alpha C_{ij}^{\alpha} C_{\alpha k}^{ \beta} + 
C_{jk}^{\alpha} C_{\alpha i}^{\beta} + 
C_{ki}^{\alpha}C_{\alpha j}^\beta = 0.
$$
Together with the antisymmetry in the lower indices it implies
\begin{equation}\label{eq:1-2prim}
\CE^*_i C^\gamma_{*k} =
\CE^*_k C^\gamma_{*i} + C^*_{ik} \CE^\gamma_*
\end{equation}
In the 2nd order ($N=2$),
$$
A_1^2 C^\gamma_{*\mu} \CE^*_\nu
+ A_2 (C^*_{\nu\mu}\CE^\gamma_* + \CE^*_\mu C^\gamma_{\nu *})
- (\mu \leftrightarrow \nu) = A_1 \, C_{\mu\nu}^{*} \CE^\gamma_*.
$$
After applying~(\ref{eq:1-2prim}) with $\mu = i$, $\nu = k$,
we obtain $A_1^2 + 3 A_2 = A_1$ (up to a common factor, which is in
general nonzero), hence $A_2 = 1/12$.

In higher order, we will recursively show that the odd coefficients
are zero (except $A_1$). That means that, 
among all products $A_I A_{N-I}$, only the terms $I = 1$ and 
$I = N-1$, where the product equals $A_1 A_{N-1} = 1/2 A_{N-1}$ survive. 
Now, $A_{N-1}$ is at the both sides of the equation, 
and as we suppose those to be non-zero (what is justified afterwards), 
we divide the equation by $A_{N-1}$ to obtain
$$
C^\gamma_{\mu\rho} (\CE^{2K})^\rho_\nu 
+ \left\lbrack  \CE^\rho_\nu \delta_\rho 
(\CE^{2K})^\gamma_\mu  \right\rbrack -(\mu\leftrightarrow \nu) 
= 2 C^\sigma_{\mu\nu} (\CE^{2K})^\gamma_\sigma 
$$
where $2K+1 = N$. 

In even order $N = 2k\geq 4$, the RHS is zero. 
Because of the different shape
of the tensors involved, we split the LHS
into the part $I = N$, and the rest, 
which we then move to the RHS to obtain
\begin{equation}\label{eq:MnoK}
 -A_N (M_0 + M_1 +\ldots + M_{N-1}) = 
\sum_{I=1}^{N-1} A_I A_{N-I} \left\lbrace \lbrack
\delta_\rho (\CE^I)^\gamma_\mu\rbrack
(\CE^{N-I})^\rho_\nu - (\mu\leftrightarrow\nu)\right\rbrace
\end{equation}
for $N$ even.
The expression in the curly brackets will be denoted
by $(K_{I,N-I})^\gamma_{\mu\nu}$. 

As $\gamma, \mu, \nu$ will be fixed, and we prove the identities
for all triples $(\gamma, \mu, \nu)$ we will just 
write $K_{I,N-I}$ without indices.
In this notation~(\ref{eq:MnoK}) reads,
\begin{equation}\label{eq:withK}
-A_N (M_0 + M_1 +\ldots + M_{N-1}) = 
\sum_{I=1}^{N-1} A_I A_{N-I} K_{I,N-I}.
\end{equation}
As an extension of this notation, we may also denote
$$
K_{0,N} = M_0 + M_1 +\ldots + M_{N-1}.
$$
Notice that
$$
[x_\alpha (\CE^I)^\alpha_\mu,x_\beta (\CE^I)^\beta_\nu] = 
x_\gamma (K_{I,N-I})^\gamma_{\mu\nu}, \,\,\,\,I=0,1,\ldots,N.
$$
{\bf Lemma.} For $L = 0,1,2,\ldots$ and $1\leq \mu,\nu,\gamma\leq n$
\begin{equation}\label{eq:shiftd}
\CE^\rho_\nu C^*_{\mu\rho} (\CE^{L})^\gamma_* 
- \CE^\rho_\mu C^*_{\nu\rho} (\CE^{L})^\gamma_* 
= C^\sigma_{\mu\nu} (\CE^{L+1})^\gamma_\sigma
\end{equation}
\begin{equation}\label{eq:ensuresOdd}
 C^\gamma_{\mu\rho} (\CE^{L})^\rho_\nu 
+ \left\lbrack  \CE^\rho_\nu \delta_\rho (\CE^{L})^\gamma_\mu  \right\rbrack
-(\mu\leftrightarrow \nu) = 2 C^\sigma_{\mu\nu} (\CE^{L})^\gamma_\sigma
\end{equation}

{\it Proof.} Equation~(\ref{eq:shiftd}) follows easily 
by applying the Jacobi identity in form~(\ref{eq:1-2prim}) to 
the expression $\CE^\rho_\nu C^*_{\mu\rho}-\CE^\rho_\mu C^*_{\nu\rho}$
and contracting with $\CE^L$.

Equation~(\ref{eq:ensuresOdd}) will be proved by induction. 
For $L = 0$ it boils down to~(\ref{eq:1-2prim}). 
We need to verify directly also $L = 1$ because this will
also be used in the proof for the step of induction. 

Suppose that~(\ref{eq:ensuresOdd}) holds for $L$. 
Then for $L+1$, by Leibniz rule, the LHS is
$$\begin{array}{lcl}
C^\gamma_{\mu\rho} (\CE^{L+1})^\rho_\nu 
+ \CE^\rho_\nu \CE^*_\mu \delta_\rho (\CE^{L+1})^\gamma_{*} 
+ \CE^\rho_\nu C^*_{\mu\rho} (\CE^{L})^\gamma_* - (\mu\leftrightarrow\nu)
\end{array}$$
According to~(\ref{eq:shiftd}) the term 
$\CE^\rho_\nu C^*_{\mu\rho}(\CE^{L})^\gamma_*
- (\mu\leftrightarrow\nu)$ contributes to exactly one half of required RHS. 
Hence it is sufficient to prove that
\begin{equation}\label{eq:ensuresAlt}
\begin{array}{lcl}
C^\gamma_{\mu\rho} (\CE^{L+1})^\rho_\nu
+ \CE^\rho_\nu \CE^*_\mu \delta_\rho (\CE^{L})^\gamma_{*} 
- (\mu\leftrightarrow\nu) = C^\sigma_{\mu\nu} (\CE^{L+1})^\gamma_\sigma
\end{array}\end{equation}
This equation is then proved by induction. Suppose it holds for $L$,
for all $\gamma$. Then multiply by $\CE^\tau_\gamma$ and sum over $\gamma$
to obtain $C^\sigma_{\mu\nu} (\CE^{L+2})^\tau_\sigma$ at RHS and
$$
\CE^\tau_\gamma C^\gamma_{\mu\rho} (\CE^{L+1})^\rho_\nu
+ \CE^\rho_\nu \CE^*_\mu \delta_\rho (\CE^{L+1})^\tau_* 
- \CE^\rho_\nu (\CE^{L+1})^*_\mu C^\tau_{*\rho}
-(\mu\leftrightarrow\nu)
$$
at LHS (we used the Leibniz rule again). Using the antisymmetrization
in $(\mu\leftrightarrow\nu)$ and renaming some dummy indices we obtain
$$
(C^*_{\mu\rho} \CE^\tau_* 
+ \CE^*_\mu C^\tau_{\rho *}) (\CE^{L+1})^\rho_\nu 
+ \CE^\rho_\nu \CE^*_\mu \delta_\rho (\CE^{L+1})^\tau_*
- (\mu\leftrightarrow\nu).
$$
Using~(\ref{eq:1-2prim}) we can sum inside the brackets to obtain
$$
C^\tau_{\mu *} (\CE^{L+2})^*_\nu 
+ \CE^\rho_\nu \CE^*_\mu \delta_\rho (\CE^{L+1})^\tau_*
- (\mu\leftrightarrow\nu)
$$
as required. 

{\bf Corollary.} {\it
If $A_0 = 1$ and $A_1 = 1/2$, $A_{2K+1}=0$ for $K = 1,\ldots,K_0-1$ 
and the relation~(\ref{eq:diffphi}) holds for $N$ odd where
$N = 2K+1$ with $K<K_0$,
then the relation~(\ref{eq:diffphi}) also holds for $N = 2K_0+1$.
}

{\bf Corollary.} {\it In particular, 
relation~(\ref{eq:diffphi}) holds for N odd
if $A_K = (-1)^K B_K/K!$.
}

The even case will require much longer calculation.
Before that, observe that $K_{I,N-I} = K_{N-I,I}$, and regarding that 
for $N>2$, the term $A_1 A_{N-1} = 0$, hence all odd-label terms
are zero, hence the equation~(\ref{eq:withK}) for even $N\geq 4$ reads
\begin{equation}\label{eq:evenKeq} 
-A_N (M_0 + M_1 +\ldots + M_{N-1}) = 
\sum_{k=1}^{\frac{N}{2}} A_{2k} A_{N-2k} K_{2k,N-2k}.
\end{equation}

\section{Hierarchy of formulas and the basis of identities}
\label{sec:hierarchy}

Throughout this section the order $N \geq 2$.  

\nxpoint ({\bf $Z$-tensors, $\hat{Z}$-tensors})
For $1\geq \gamma, \mu,\nu \leq n$, $\mu\neq\nu$ 
define the (components of) $Z$-tensor 
\begin{equation}
(Z^{l,m,k})^\gamma_{\mu\nu} := 
(\CE^l)^*_\mu (\CE^m)^*_\nu C^*_{**} (\CE^k)^\gamma_*
- (\mu\leftrightarrow\nu),
\end{equation}
where $l,m,k\geq 0$, $l+m+k+1 = N$. Recall that $C^r_{su}$ are either
the structure constants of a  Lie algebra of rank $n$ as a $\genfd$-module,
or, in the universal case, the generators of $\genfd[{\mathcal C}_n]$.
However, the RHS makes also sense when the 
$C^r_{su}$, for $r,s,u = 1,\ldots,n$, are simply the $n^3$ 
generators of a free commutative $\genfd$-algebra 
(no Jacobi, no antisymmetry). In that
case, the LHS will be denoted $(\hat{Z}^{l,m,k})^\gamma_{\mu\nu}$.

From now on, whenever $\mu,\nu,\gamma$ are of no special 
importance, we write simply $Z^{l,m,k}$, and $\hat{Z}^{l,m,k}$, 
i.e. $\gamma,\mu,\nu$ will be skipped from the notation 
whenever they are clear from the context. 
As before, $*$-s are dummy indices, and upper $*$-s are 
contracted to lower stars pairwise in left-to-right order. 
$(Z^{l,m,k})^\gamma_{\mu\nu}$ are components of a rank 3 tensor,
antisymmetric in lower indices. Each of the two summands
is a contraction of $l+m+k$ $\CE$-s and one $C$. 
Such expressions appear in our analysis when $N = l+m+k+1$. 
By obvious combinatorial arguments, this tensor $T$ is 
contraction-connected
(the copies of the generators 
($C$-s and $\partial$-s) involved can not be
separated into two disjoint subsets 
without a contraction involving
elements in different subsets). In a universal case, the components
of this tensor lie in $\genfd[{\mathcal C}_n][\partial^1,\ldots,\partial^n]$.
The $Z$-tensors may be called also ``star-tensor'', what points 
to a useful graphical notation in which 
the three branches $\CE^l,\CE^m,\CE^k$  are drawn respectively
left, down and right from the central ``node'' $C$, 
attached by lines denoting contractions to the three indices of $C$.

\nxpoint ({\bf Special cases of $Z$-tensors: $b_i$, $M_j$})
Any of $l,m,k$ may take value zero: $(\CE^0)^i_j = \delta^i_j$
is then the Kronecker tensor. This mean that a ``branch'' was cut and
we have monomials in the tensor which are linear ``$M$-'' and ``$b$-'' chains
denoted $M_0,\ldots, M_{N-1}, b_0, b_1,\ldots, b_{N-1}$, where
for all $0\leq k\leq N-1$, 
$$\begin{array}{l}
b_k := Z^{k,N-k-1,0} = (\CE^k)^*_\nu C^\gamma_{**} (\CE^{N-k-1})^*_\mu 
- (\mu\leftrightarrow\nu)\\
M_k := Z^{k,0,N-k-1} = (\CE^{k})^*_\mu C^*_{*\nu} (\CE^{N-k})^\gamma_* 
 - (\mu\leftrightarrow\nu)\\
\end{array}$$
Then
\begin{equation}
\label{eq:KvsZ}\begin{array}{lcl} 
M_0 &=& Z^{0,0,N-1} = 
C^*_{\mu\nu} (\CE^{N-1})^\gamma_* - (\mu\leftrightarrow\nu)\\
M_{N-1} &=& Z^{N-1,0,0} = (\CE^{N-1})^*_\mu C^\gamma_{*\nu} 
- (\mu\leftrightarrow\nu) = b_0 = b_{N-1} \\
K_{I,N-I} &=& 
\left[\delta_\rho (\CE^I)^\gamma_\mu\right](\CE^{N-I})^\rho_\nu -
\left[\delta_\rho (\CE^I)^\gamma_\nu\right](\CE^{N-I})^\rho_\mu \\
&=& \sum_{l=0}^{I-1} (\CE^l)^*_\mu C^*_{*\rho} 
(\CE^{I-l-1})^\gamma_* (\CE^{N-I})^\rho_nu - (\mu\leftrightarrow\nu)\\
&=& \sum_{l=0}^{I-1} Z^{l,N-I,I-l-1}.
\end{array}\end{equation}
\nxpoint ({\bf Spaces of $Z$-tensors})
The $\genfd$-span of all $Z$-s is denoted 
${\mathcal Z}_N = 
{\mathcal Z}_{\genfd, N,\mu,\nu}^\gamma \subset \genfd[{\mathcal C}_n]
[\partial^1,\ldots,\partial^n]$. Similarly, the $\genfd$-span of all 
$\hat{Z}$-s is denoted by $\hat{\mathcal Z}_{\genfd, N,\mu,\nu}^\gamma$.
First of all we need 

{\bf Lemma.} {\it (Before we quotient out by Jacobi identities and
antisymmetry), if $n$ is sufficiently big, all $\hat{Z}^{l,m,k}$
are $\genfd$-linearly independent.
}

This lemma is in the setting of the free polynomial algebra on $C^r_{su}$
(tensored by the symmetric algebra in $\partial^i$-s), 
hence finding the monomial summands in $\hat{Z}^{l,m,k}$, which comprise
a part of the standard basis, and which are not summands 
in any other $\hat{Z}^{l',m',k'}$, is easy to do, if there are sufficiently
many distinct indices to choose from, say, one distinct from each
contraction. On the other extreme, if $n=1$, all contractions involve
the same index, and hence there are degeneracies. We leave more precise
argument (proof of the lemma) as a combinatorial exercise for the reader.

Clearly, $({\mathcal Z}_N)^\gamma_{\mu\nu} = \hat{\mathcal Z}_N/(J_N+J'_N)$, 
where $J_N$, $J'_N$ are the two submodules defined as follows.
$J_N$ consists of all elements of the form  
$\sum_{ijk} r^p_{ijk} (Q^{ijk}_p)^\gamma_{\mu\nu}$
where $r^p_{ijk} = C_{ij}^* C^p_{*k} + C_{jk}^* C^p_{*i} + C_{ki}^* C^p_{*j}$,
and $Q^{ijk}_p$ is a tensor involving $(N-2)$ $C$-s and $(N-1)$ $\partial$-s
(those numbers are fixed, because both the Jacobi identities 
and the antisymmetry are {\it homogeneous} relations), 
and of external $GL_n(\genfd)$-covariance $()^\gamma_{\mu\nu}$. 
Contraction-connected tensors are in generic case linearly independent from
disconnected, hence by counting free and contracted indices, we observe that
in each monomial involved in $Q^{ijk}_p$, at least one of the indices
$i,j,k$, is contracted to a $\partial$-variable.  
Similarly, using $J'_N$, one handles the antisymmetry. 

Our strategy to determine the structure of $({\mathcal Z}_N)^\gamma_{\mu\nu}$
is as follows: we start with $\hat{\mathcal Z}_N$
and then determine the submodule of relations $J_N+J'_N$, trying to eliminate
superfluous generators to the point where algebraic analysis will
become explicit enough. We notice first, that the terms where
$r_{ijk}$ is coupled to $\partial^i$ and $\partial^j$ simultaneously
are superfluous: by antisymmetry (using $J'_N$) this will be the contraction
of a symmetric and antisymmetric tensor in one of the summands, and
the two others are equal in $\hat{\mathcal Z}_N/J_N'$. Similarly, with other
double contractions of $r_{ijk}$ with $\partial$. In other words, 
Jacobi identities always come in the covariant combination obtained by
contracting exactly one of the three indices to $\partial$, as in 
\eqref{eq:1-2prim}.
In degree $N$, exactly the following relations are of that Jacobi type:
\begin{equation}\label{eq:central}
 Z^{l,m+1, k} = Z^{l,m,k+1} - Z^{l+1,m,k}.
\end{equation}
By this identity we can do recursions in $l, m$ or $k$. 
When $l$ or $m$ is zero $Z$-s are $M$-s. Hence $M$-chains span
the whole ${\mathcal Z}_N$ (but they are not independent).
The antisymmetry in lower indices
of $C^i_{jk}$ is nontrivial only if neither of the two indices are contracted
to $\partial$-s. But that means, that this is the central node of the
star-tensor, and the interchange of the two indices may be traded for
the interchange of nodes of $\mu$ and $\nu$ branch, what results in
the identities
\begin{equation}\label{eq:sym}
Z^{l,m,k} = Z^{m,l,k}
\end{equation}
thereafter called {\bf symmetries}. In particular, $b_k = b_{N-k-1}$.

In addition to ${\mathcal Z}^\gamma_{\genfd N \mu\nu}$
we introduce also spaces ${\mathcal F}^\gamma_{\genfd N}$ and
${\mathcal Z\mathcal F}_{\genfd N}$ 
as follows:
${\mathcal F}_{\genfd, N}$ is the free $\genfd$-module
generated by the symbols $F^{l,m,k}$ where $l,m,k \geq 0$
and $l+m+k+1 = N$; furthermore,
${\mathcal Z\mathcal F}_{\genfd N} := 
{\mathcal F}_{\genfd N}/I$, where $I$ is the submodule
generated by symmetries $F^{l,m,k} - F^{m,l,k}$ and
relations $F^{l,m+1,k}+ F^{l+1,m,l} - F^{l,m,k+1}$.
Image of $F^{l,m,k}$ in ${\mathcal Z\mathcal F}_{\genfd N}$ 
will be denoted by $Z_F^{l,m,k}$.

The summary of the above discussion may be phrased as follows:

\nxpoint \label{p:suffbig}
{\bf Proposition.} {\it 
Let $N\geq 1$, $n$ sufficiently big, $1\leq \mu,\nu,\gamma \leq n$ and
$\mu\neq \nu$. 
The correspondence $Z^{l,m,k}_F \mapsto (Z^{l,m,k})^\gamma_{\mu\nu}$
extends to a $\genfd$-module isomorphism 
${\mathcal Z\mathcal F}_{\genfd N}\cong 
{\mathcal Z}^\gamma_{\genfd N\mu\nu}$.
}

Thus from now on, we may work with the presentation for 
${\mathcal Z}^\gamma_{\genfd, N,\mu,\nu}$. 
Using~\eqref{eq:central} we can do recursions to express $Z$-s in 
terms of special cases when one of the labels is zero: namely $M$-s or $b$-s.
We can choose a distinguished way to do recursion, e.g. in each step
lower $m$ in each monomial. One ends with $Z$ expressed in terms of $M$-s
in a distinguished way, and each symmetry will be an identity between 
$M$-s, and there are no other identities in ${\mathcal Z}_N$. If we 
mix various recursions, this is the same as doing the algorithm of
distinguished  recursion, but at some steps intercepted
by applying a symmetry. 

\nxpoint {\bf Lemma.} {\it 
The distinguished way of recursion, always lowering $m$, yields
\begin{equation}\label{eq:canonForm}
Z^{l,m,k} = \sum_{j=0}^m (-1)^j {m \choose j} M_{l+j} 
\end{equation}
The alternative way of recursion, always lowering $k$, yields
\begin{equation}\label{eq:altcanonForm}
Z^{l,m,k} = \sum_{j=0}^k  {k \choose j} b_{l+j} 
\end{equation}
\begin{equation}\label{eq:Mb}
M_{N-k-1} = \sum_{i=0}^k {k \choose i} b_{i} 
\end{equation}
Moreover,
\begin{equation}\label{eq:K0N}
M_0 + M_1 + \ldots + M_{N-1} = 
\delta_{N,\mathrm{odd}} {N\choose \frac{N+1}{2}} b_{(N-1)/2} +
\sum_{i=0}^{\lfloor N/2\rfloor -1} {N+1\choose i+1} b_i 
\end{equation}
and for any $N$ even and $1 \leq I\leq N/2$,
\begin{equation}\label{eq:Keqf1}
K_{I,N-I} = \sum_{i=0}^{I-1} {I\choose i} b_i 
\end{equation}
}
\newline
{\it Proof.} The frist two formulas follow by easy induction.
Too see~\eqref{eq:K0N} we use~\eqref{eq:altcanonForm}: 
$$M_{M-k-1} = Z^{N-k-1,0,k} = \sum_{j=0}^k {k \choose j} b_{N-k-1+j}
= \sum_{i=0}^k {k \choose i} b_{N-i-1}= \sum_{i=0}^k {k \choose i} b_{i}$$
To see~\eqref{eq:K0N}, notice that from~\eqref{eq:Mb} we have
$$\begin{array}{lcl}
\sum_{k=0}^{N-1} M_{N-k-1} &=& b_0 + (b_0 + b_1) + (b_0 + 2b_1 + b_2) 
+\ldots  
\\ &=& Nb_0 + {N\choose 2} b_1 +\ldots + b_{N-1} =
\sum_{i=0}^{N-1} {N\choose i+1} b_i \\ 
&=& \delta_{N,\mathrm{odd}} {N\choose (N-1)/2} b_{(N-1)/2}
+ \sum_{i=0}^{\lfloor N/2\rfloor -1} 
\left({N\choose i+1}+{N\choose N-i}\right) b_i
\end{array}$$
with~\eqref{eq:K0N} immediately.
To prove~\eqref{eq:Keqf1} we proceed as follows
$$\begin{array}{lcl}
K_{I,N-I} &\stackrel{\eqref{eq:KvsZ}}=& \sum_{l = 0}^{I-1} Z^{l,N-I, I-l-1} =
\sum_{j = 0}^{I-1} Z^{I-j-1,N-I,j} 
\stackrel{}=  b_{I-1} + (b_{I-1} + b_{I-2}) \,+ 
\\ && \,\,\,\,\,\,\,\,\,+\, (b_{I-1} + 2b_{I-2} + b_{I-3})
+\ldots +\, (b_{I-1} +\ldots + (I-1)b_1 + b_0)
\\ &=& I b_{I-1} + {I\choose 2} b_{I-2} +\ldots + b_0
= \sum_{j =1}^{I} {I\choose j} b_{I-j} =
\sum_{i=0}^{I-1} {I\choose i} b_{i} . 
\,\,\,\,\,\,\,\mathrm{Q.E.D.}
\end{array}$$

\nxpoint This result enables us that we can effectively do
all computations with $Z$-tensors either 
in terms of $b$-s or in terms of $M$-s. 
However, both sets of variables have internal linear
dependences which we will now study. 

First the case of $b$-s, which is much simpler. To study the relations
among the relations we introduce $\genfd$-module 
${\mathcal F}_{b\genfd N}$ as the free module on $N$ 
symbols $\ovb_0,\ldots, \ovb_{N-1}$. modulo
the relations $\ovb_i = \ovb_{N-i-1}$. 

The relations ~\eqref{eq:altcanonForm} will be now taken as the definitions
of $Z$-variables. More precisely, define 
$Z_b^{l,m,k} = \sum_{j=0}^k \ovb_{l+j} 
\in {\mathcal F}_{b\genfd N}$. The only relations among
$Z_b^{l,m,k}$ are symmetries 
$X_{b,k}^{(s)} := (Z_b)^{l,l+s,N-2l-s-1} - (Z_b)^{l+s,l,N-2l-s-1}$.
Now 
$$X_{b,k}^{(s)} = \sum_{j=0}^{N-2l-s-1} 
\bar{b}_{l+j}-\sum_{i = 0}^{N-2l-s-1}\bar{b}_{l+s+i}
= \sum_{j=0}^{N-2l-s-1} (\bar{b}_{l+j}-\bar{b}_{N-l-j-1}),$$
where we replaced $i$ by $j = N-2l-s-1-i$ in the second sum.
Thus every symmetry is a linear combination of the $\lfloor N/2\rfloor$
relations $\ovb_i - \ovb_{N-i-1}$, which are linearly
independent in ${\mathcal F}_{b\genfd N}$. This result, together 
with the way definitions were set implies
 
\nxpoint {\bf Theorem.} {\it The correspondence
$Z_F^{k,l,m}\mapsto Z_b^{k,l,m}$ extends to a unique isomorphism
$${\mathcal F}_{\genfd N}\stackrel\cong\longrightarrow 
{\mathcal F}_{b\genfd N}/
\langle \ovb_i - \ovb_{N-i-1}, i = 0,\ldots, \lfloor (N-1)/2\rfloor\rangle.$$
}

Consequently, using~\refpt{p:VFonW} we get,

{\bf Corollary.} {\it For sufficiently large $n$, 
${\mathcal Z}{\mathcal F}_{\genfd N}$ is canonically isomorphic
to ${\mathcal F}_{b\genfd N}/
\langle \ovb_i - \ovb_{N-i-1}\rangle$. In particular,
$\dim_{\mathbb Q} {\mathcal Z}_{{\mathbb Q}N\gamma\mu\nu} = 
N - \lfloor N/2\rfloor = \lfloor (N+1)/2\rfloor = \lceil N/2 \rceil$.
}

In plain words, there are no relations among the $b_i$-s
except (consequences of) the $\lfloor N/2 \rfloor$ 
symmetries for $b$-chains: $b_i = b_{N-i-1}$.
The remainder of this section will be dedicated to a much harder
analogue of above calculation concerning $M$-variables (not used in
further sections). 

\nxpoint Let ${\mathcal F}_{M\genfd N}$ be the free
$\genfd$-module on $N-1$ symbols $\bar{M}_i$, $i = 0,\ldots, N-1$ 

Consider the ``special symmetries'' 
\begin{equation}\label{eq:Xdef}\begin{array}{lcl}
X_k &:=& Z^{k, k+1, N-2k-2}_M - Z^{k+1, k, N-2k-2}_M,
\,\,\,\,\,\,\,k = 0,1,\ldots, \lfloor (N-1)/2\rfloor,
\end{array}\end{equation}
where $Z^{k,l,m}_M$ is defined to mimic $Z^{k,l,m}$ 
expressed in terms of $M$-s, 
using the distinguished recursion lowering $m$, 
i.e. by the following version of~\eqref{eq:canonForm}:
$$
Z^{k,l,m}_M := \sum_{j=0}^m (-1)^j {m \choose j} \ovM_{l+j}.
$$ 
$$
X_k = \ovM_k + \sum_{j= 0}^k (-1)^{j+1}\left[ 
{k\choose j} + {k+1\choose j+1}
\right] \ovM_{k+j+1},\,\,\,\,\,\,\,k\leq N/2-1
$$
$$\begin{array}{lcl}
X_0 = \ovM_0 - 2\ovM_1  &=& \ovM_0\uparrow (1,-2)\\ 
X_1 = \ovM_1 - 3\ovM_2 + 2\ovM_3 &=& \ovM_1\uparrow (1,-3,2)\\
X_2 = \ovM_2 - 4\ovM_3 + 5\ovM_4 - 2\ovM_5 &=& \ovM_2 \uparrow(1,-4,5,-2) \\
X_3 = \ldots & = &\ovM_3\uparrow (1,-5,9,-7,2)\\
X_4 = \ldots & = &\ovM_4\uparrow (1,-6,14,-16,9,-2)\\
\end{array}$$
and in general
\begin{equation}
X_k = \ovM_{k} \uparrow \left(1,-k,\ldots, 
(-1)^{k-1}\left({k \choose 2}+{k+1\choose 2}\right),
(-1)^k(2k+1),(-1)^{k+1} 2\right)
\end{equation}
Our aim now is to show that all other symmetries are expressed
in terms of $X_k$. For $s\geq 0$, these formulas generalize to
\begin{equation}\label{eq:symmVersusX}
X^{(s)}_j := Z^{j,j+s,*}_M - Z^{j+s,j,*}_M = 
\sum_{k = 0}^{\lfloor\frac{s-1}{2}\rfloor} 
(-1)^k {s-k-1\choose k} X_{j+k}
\end{equation}
where $2j+s+ * = N$ and $\lfloor \rfloor$ denotes the greatest integer part.
In particular,
$$\begin{array}{l}
X^{(0)}_j = 0,\,\,\,\,X^{(1)}_j = X^{(2)}_j = X_j,\\
X^{(3)}_j = X_j - X_{j+1},\\ X^{(4)}_j = X_j - 2X_{j+1},\\
X^{(5)}_j = X_j - 3X_{j+1} + X_{j+2},\\
X^{(6)}_j = X_j - 4 X_{j+1}+ 3X_{j+2},\\
X^{(7)}_j = X_j - 5X_{j+1} + 6X_{j+2}-X_{j+3},\\
X^{(8)}_j = X_j - 6X_{j+1} + 10X_{j+2} - 4X_{j+3},\\
X^{(9)}_j = X_j - 7X_{j+1} + 15X_{j+2} - 10X_{j+3} + X_{j+4}
\end{array}$$ 
as it is easy to verify directly. 
Formula~\eqref{eq:symmVersusX} in general follows from a cumbersome induction.
Alternatively, by~\eqref{eq:canonForm} and~\eqref{eq:Xdef}
the equation~\eqref{eq:symmVersusX} is equivalent to the following
statement in the free $\genfd$-module on $\ovM$-symbols:
$$\begin{array}{l}
\sum_{i=0}^{j+s} (-1)^i {j+s\choose i} \ovM_{j+i} -
\sum_{i=0}^j (-1)^i {j\choose i} \ovM_{j+s+i} 
=\sum_{k=0}^{[\frac{s-1}{2}]} (-1)^k {s-k-1\choose k} \,\times \\ 
\,\,\,\,\,\,\,\,\,\,\,\,\,\,\,\,\,\,\,\,\,\,\,\,\,
\times \,\left(\sum_{i=0}^{j+k+1} (-1)^i {j+k+1\choose i} \ovM_{j+k+i}
- \sum_{i=0}^{j+k} (-1)^i {j+k\choose i} \ovM_{j+k+i+1}\right).
\end{array}$$
Equate the coefficients in front of $(-1)^{i} \ovM_{j+i}$ 
to obtain that~\eqref{eq:symmVersusX} 
is equivalent to the assertion that, for all $i$,
\begin{equation}\label{eq:jsi}
{j+s\choose i} - (-1)^s{j\choose i-s} =
\sum_{k=0}^{\lfloor\frac{s-1}{2}\rfloor}
{s-k-1\choose k}\left({j+k+1\choose i-k}+{j+k\choose i-k-1}\right),
\end{equation}
with the convention ${a\choose b} = 0$ whenever $b$ falls out of the range
$0\leq b\leq a$. 
E.g. if $s=5$, $j=4$, the equality reads $126=126$ for $i=4$
and $88=88$ for $i=6$.
\nxpoint {\bf Proposition.} {\it Equation~\eqref{eq:jsi} holds for all 
integers $i>0, j>0, s>0$. } 

{\it Proof.} Label the equation~\eqref{eq:jsi} with the triple $(j,s,i)$.
Roughly speaking, we do the induction on $s$, however being a bit careful
with the choice of $j$ and $i$ in the induction step (it seems
that the direct implication $(j,s,i)\Rightarrow (j,s+1,i)$ is far too
complex to be exhibited). 


If $i = 1$, then~\eqref{eq:jsi} is the tautology $j+s = j+s$.
If $s = 1, i>0,j>0$, the equation~\eqref{eq:jsi} 
becomes the tautology 
${j+1\choose i} + {j\choose i-1} = {j+1\choose i} + {j\choose i-1}$,
and for $s = 2$ the simple identity
${j+2\choose i} - {j\choose i-2} = {j+1\choose i} + {j\choose i-1}$.
If $i\geq 2$, then adding LHS for $(j,s,i)$ 
and for $(j+1,s-1,i-1)$ together, we get 
${j+s\choose i} + {j+s\choose i-1} - (-1)^{s-1} 
\left({j+1\choose i-s} - {j\choose i-s}\right) = 
{j+s+1\choose i} - (-1)^{s+1} {j\choose i - s-1} 
$, what is the LHS of $(j,s+1,i)$. When adding the RHS 
of $(j,s,i)$ and of $(j+1,s-1,i-1)$ add pairwise the summands
for label $k$ in $(j,s,i)$ and those for label $k-1$ in $(j+1,s-1,i-1)$:
for such combination of $k$-s the expression in the brackets on the RHS is
identical and the prefactors involving $s$ add nice.  


The identity~\eqref{eq:jsi} therefore holds
for all $j>0, s>0, i>0$.~~~Q.E.D.

We have thus proved

\nxpoint {\bf Theorem.} {\it 
Equation ~\eqref{eq:symmVersusX} holds for all $s>0$. 
In particular, $X_j^{(s)}$ belong to the $\genfd$-linear span 
of special symmetries $X_0,X_1,\ldots,X_{\lfloor (N-1)/2\rfloor}$.
The formulas for $X_j = X_j^{(1)}$ show that they are $\genfd$-linearly 
independent in ${\mathcal F}_{M\genfd N}$.
The correspondence $Z^{k,l,m}_F\mapsto Z_M^{k,l,m}$ 
extends to a well-defined isomorphism of $\genfd$-modules
${\mathcal F}_{\genfd N}\to {\mathcal F}_{M\genfd N}/
\langle X_k, k=1,\ldots, \lfloor (N-1)/2\rfloor \rangle$.  
}

This gives another way to see the dimension of the space of $Z$-tensors 
over a field: $\dim_{\mathbb Q} {\mathcal Z}_{{\mathbb Q}N\gamma\mu\nu}
= \#\ovM_i - \# X_k = N - \lfloor N/2 \rfloor = \lfloor (N+1)/2\rfloor$. 
With respect to the obvious natural  $GL_n(\genfd)$-actions, 
the isomorphism in the theorem are also $GL_n(\genfd)$-equivariant. 

\section{Formula involving derivatives of $\coth(x/2)$}

We return to proving~\eqref{eq:evenKeq} for Ansatz $A_K = B_K/K!$
for $N\geq 4$ even. This reads
\begin{equation}\label{eq:Nthordercondition}
\beta_N K_{0,N} +
\sum_{k = 1}^{\frac{N}{2}} 
\beta_{2k} \beta_{N-2k} K_{2k,N-2k} = 0,
\,\,\,\,\,\mathrm{where}\,\,\,\beta_{2k} := \frac{B_{2k}}{(2k)!}.
\end{equation}
Let $I \leq N/2-1$. Using ~\eqref{eq:K0N} and~\eqref{eq:Keqf1}
and replacing $b_i = b_{N-i-1}$ for $i\geq N/2$ we get
$$
\sum_{i=0}^{N/2-1}{N+1 \choose i+1} \beta_N b_i
+ \sum_{k = 0}^{N/2} \sum_{i = 0}^{2k-1} {2k\choose i}
\beta_{2k}\beta_{N-2k} b_i = 0
$$
Now we recall that not all $b_i$-s are independent. 
In the case of $\gg = su(2)$ one has $b_i = 0$ for 
$1\leq i \leq N/2-1$ and one has only the terms with $b_0$.
In the generic/universal case,
the $M_{N-1}=b_0, b_1,\ldots, b_{N/2-1}$-s are independent, hence
we need $\alpha_i = 0$ for all $i$. 
But it is easy to replace $b_{i}$ by $b_{N-i-1}$ whenever $i\geq N/2$ 
and then we get an equation of the form
$\sum_{i = 0}^{N/2-1} \alpha_i b_i = 0$,
where $\alpha_0,\ldots, \alpha_{N/2-1}$ are rational numbers
depending on $N$ and $i$. Thus in universal case, 
~\eqref{eq:Nthordercondition} boils down to $\alpha_i = 0$
for all $0 \leq i \leq N/2-1$. 
If $i < N/2$, then the coefficient in front of $b_i$ in the second 
(double) sum is
$$
\sum_{k = \left\lfloor \frac{i}{2} \right\rfloor +1}^{N/2 -1}
\beta_{2k}\beta_{N-2k}
{2k\choose i}
+ \sum_{l = \left\lfloor \frac{N-i-1}{2}\right\rfloor+1}^{N/2-1}
\beta_{2l}\beta_{N-2l}{2l\choose N-i-1}
$$
In these two sums, it will be convenient to count the summation index from
$N/2$ downwards ($k\mapsto N/2 -k$ etc.). 
Thus, ~\eqref{eq:Nthordercondition} is equivalent to 
$$
{N+1 \choose i+1} \beta_N +
\sum_{k = 1}^{N/2 -1-\left\lfloor \frac{i}{2}\right\rfloor}
\beta_{2k}\beta_{N-2k}
{N-2k\choose i} 
+ \sum_{l = 1}^{\left\lfloor \frac{i}{2}\right\rfloor}
\beta_{2l}\beta_{N-2l}{N-2l\choose i-2l+1} = 0
$$
for all $0 \leq i \leq N/2-1$. 
Now we subtract ${N \choose i} \beta_N = {N \choose i} \beta_N \beta_0$ 
from the first summand and absorb it into the first sum, 
as the new $k = 0$ summand. Then, the remainder of the first summand is
${N\choose i+1}  \beta_N$ what can be absorbed into the second sum,
as the new $l = 0$ summand. We obtain
\begin{equation}\label{eq:alphai} 
\sum_{k = 0}^{N/2 -1-\left\lfloor \frac{i}{2}\right\rfloor}
\beta_{2k}\beta_{N-2k}
{N-2k\choose i} 
+ \sum_{l = 0}^{\left\lfloor \frac{i}{2}\right\rfloor}
\beta_{2l}\beta_{N-2l}{N-2l\choose i-2l+1} = 0
\end{equation}
The generating function for the even Bernoulli numbers is
\begin{equation}\label{eq:f}
f(x) := (x/2)\coth(x/2) = \sum_{J = 0}^\infty \frac{B_{2J}}{(2J)!} x^{2J}
= \sum_{J = 0}^\infty \beta_{2J} x^{2J}.
\end{equation}
It can be easily checked that the equation~\eqref{eq:alphai} 
follows from the following functional equation for $f$:
\begin{equation}\label{eq:ffunct}
\frac{1}{i!} ff^{(i)} 
+ \sum_{k = 0}^{ \left\lfloor \frac{i}{2}\right\rfloor}
\beta_{2k} \frac{x f^{(i-2k+1)}}{(i-2k+1)!}
- \delta_{i,{\rm even}} \beta_i f = 0 
\end{equation}
where $\delta_{i,{\rm even}} = 1$ if $i$ is even and vanishes otherwise.
To see this we write the whole LHS as a power
series expansion in $x$ and then for fixed $N$ even number with $N>i$ 
we equate the coefficient in front of  $x^{N-i}$ with zero. So obtained
equation multiply with $(N-i)!$ and we obtain~\eqref{eq:alphai}.
The summand with $\delta_{i,{\rm even}} = 1$ 
was added to offset the extension of
the first sum in~\eqref{eq:alphai} to $k = N/2-i/2$,
for $i$ even, which is needed to make the correspondence of
that sum with $\frac{1}{i!} ff^{(i)}$.

The rest of the section is devoted to the proof of the functional
equation~\eqref{eq:ffunct}.

Denote by $g := \coth(x/2)$ and $g^{(j)} := \frac{d^j}{dx}\coth(x/2)$.
Then $f = xg/2$ and the equation~\eqref{eq:ffunct} becomes
$$\begin{array}{l}
\left\lbrack \frac{1}{2} \frac{1}{i!} gg^{(i)} + 
\sum_{k=0}^{\left\lfloor \frac{i}{2}\right\rfloor} \beta_{2k}
\frac{g^{(i-2k+1)}}{(i-2k+1)!}\right\rbrack x^2
+  \\ \,\,\,\,\,\,\,\,\,\,\,\,\,
+ \left\lbrack \frac{1}{2}\frac{1}{(i-1)!} gg^{(i)}
+ \sum_{k=0}^{\left\lfloor \frac{i}{2}\right\rfloor}\beta_{2k} 
\frac{g^{(i-2k)}}{(i-2k)!} - \delta_{i,{\rm even}}\beta_i g \right\rbrack x = 0
\end{array}$$
If we denote the expression in brackets in front of $x$ by $I_i$
then the expression in brackets in front of $x^2$ equals $I_{i+1}$
(the $\delta_{i,{\rm even}}$-term effectively cancels the $k = i/2$ summand). 
Note that $I_i$ and $I_{i+1}$ are hyperbolic functions, and 
$x$ and $x^2$ are linearly independent under the action by multiplication
of the algebra of all hyperbolic functions on the space of all functions
of one variable. Thus we have the equality {\it only if} $I_i = I_{i+1} = 0$. 
We will show $I_i = 0$ for all $i$ by induction. 
We rephrase this as the following statement:

{\bf
Suppose $i\geq 2$. Then the following identity holds:}
$$
\fbox{$\frac{i}{2} \coth(x/2) \frac{d^{i-1}}{dx}\coth(x/2) + 
\sum_{k = 0}^{\left\lfloor\frac{i-1}{2}\right\rfloor}
{i \choose 2k} B_{2k}\frac{d^{i-2k}}{dx}\coth(x/2) = 0.$}
$$

To simplify notation denote $g^{(j)} := \frac{d^j}{dx}\coth(x/2)$.
Then we can rewrite the statement as
\begin{equation}\label{eq:Bercothderiv}
\frac{i}{2} g g^{(i-1)} + 
\sum_{k = 0}^{\left\lfloor\frac{i-1}{2}\right\rfloor}
{i \choose 2k} B_{2k} g^{(i-2k)} = 0
\end{equation}

In fact, this series of equations (as well as the proof below) holds
also for $th(x/2)$, and, more generally, for 
$g = \frac{Ae^{x/2}+ e^{-x/2}}{Ae^{x/2}- e^{-x/2}}$
where $A$ is any constant.

{\it Proof.} {\it Basis of induction.}
Substituting $\coth(x/2)$ for $g$ observe that $2g' + g^2 - 1 = 0$, hence
$g''  + gg' = 0$, what is our identity~\eqref{eq:Bercothderiv} for $i = 2$.

{\it The step of induction.} Suppose the identity holds for $i$
and act with $d/{dx}$ to the whole identity to obtain
$$
(i/2) g' g^{(i-1)} + (i/2) g g^{(i)} + 
\sum_{k = 0}^{\left\lfloor\frac{i-1}{2}\right\rfloor}
{i \choose 2k} B_{2k} g^{(i-2k+1)} = 0.
$$
Now replace $g'$ by $(1/2) (1-g^2)$. 
By the induction hypothesis, we may also replace
$(-i/4) g^2 g^{(i-1)}$ by 
$gg^{(i)}/2 + (g/2)\sum_{s = 1}^{\left\lfloor\frac{i-1}{2}\right\rfloor}
{i\choose 2s}B_{2s} g^{(i-2s)}$, thus obtaining
$$
\frac{i+1}{2} g g^{(i)} + (i/4)g^{(i-1)} + \frac{1}{2}
\sum_{s = 1}^{\left\lfloor\frac{i-1}{2}\right\rfloor}{i\choose 2s}
B_{2s} g g^{(i-2s)} + \sum_{k = 0}^{\left\lfloor\frac{i-1}{2}\right\rfloor}
{i \choose 2k} B_{2k} g^{(i-2k+1)} = 0.
$$
Further replace $g g^{(i-2s)}/2$ 
by the sum expression given by the induction hypothesis to obtain
$$\begin{array}{ll}
\frac{i+1}{2} g g^{(i)} + (i/4)g^{(i-1)} &
-\sum_{s = 1}^{\left\lfloor\frac{i-1}{2}\right\rfloor}
\frac{1}{i-2s+1}{i\choose 2s} B_{2s}
\sum_{r = 0}^{\left\lfloor\frac{i-2s}{2}\right\rfloor}
{i - 2s+1 \choose 2r } B_{2r} g^{(i-2s-2r+1)} \\  
& + g^{(i+1)} + \sum_{k = 1}^{\left\lfloor\frac{i-1}{2}\right\rfloor}
{i \choose 2k} B_{2k} g^{(i-2k+1)} = 0.
\end{array}$$
Since the condition  
$r \leq \left\lfloor\frac{i-2s}{2}\right\rfloor$,
is equivalent to $r+s \leq  \left\lfloor\frac{i}{2}\right\rfloor$, we have
$$\begin{array}{ll}
\frac{i+1}{2} g g^{(i)} + g^{(i+1)} + 
(i/4)g^{(i-1)} & - \sum_{l = 1}^{\left\lfloor\frac{i}{2}\right\rfloor}
\sum_{s = 1}^{l} \frac{i!}{(i-2l+1)!} 
\frac{B_{2s}}{(2s)!}\frac{B_{2l - 2s}}{(2l-2s)!} g^{(i-2l+1)} 
\\ & +  \sum_{k = 1}^{\left\lfloor\frac{i-1}{2}\right\rfloor}
{i \choose 2k} B_{2k} g^{(i-2k+1)} = 0.
 \end{array}$$
Using the identity $\sum_{s = 1}^{l} \frac{B_{2s}}{(2s)!}
\frac{B_{2l - 2s}}{(2l-2s)!} = \frac{-B_{2l}}{(2l-1)!} + 
\frac{1}{4}\delta_{l,1}$, valid for $l>0$ 
(proof: use the obvious identity $xf' = f - f^2 + (x/2)^2$ 
for the generating function $f$ in~\eqref{eq:f})) 
and $B_2 = 1/12$, we then obtain
$$\begin{array}{ll}
\frac{i+1}{2} g g^{(i)} + g^{(i+1)} 
+ \sum_{l = 1}^{\left\lfloor\frac{i}{2}\right\rfloor}
{i \choose 2l-1} B_{2l} g^{(i-2l+1)} 
+  \sum_{k = 1}^{\left\lfloor\frac{i-1}{2}\right\rfloor}
{i \choose 2k} B_{2k} g^{(i-2k+1)} = 0.
\end{array}$$
Using ${i \choose 2k} + {i\choose 2k-1} = {i+1\choose 2k}$ we 
recover the equation~(\ref{eq:Bercothderiv})
with $i$ replaced by $i+1$. This finishes the proof. 

\section{ Prerequisites on formal group schemes in functorial approach }
\label{s:Formal}
In the next three sections we present an invariant derivation of
our main formula by the first author ({\sc N. Durov}), 
valid over arbitrary rings containing $\bbQ$. 
The idea of this proof is the following. Given 
an $n$-dimensional complex Lie algebra~$\gg$, we can always find a complex 
Lie group $G$ with Lie algebra equal to~$\gg$. 
Elements $X\in\gg$ correspond to right-invariant 
vector fields $X_G$ on $G$, and clearly $[X,Y]_G=[X_G,Y_G]$, 
so we get an embedding $\gg\to\Vect(G)$ of $\gg$ into the Lie algebra 
$\Vect(G)$ of vector fields on~$G$. 
If $e\in U\subset G$, $U\to\bbC^n$ is a coordinate neighborhood 
of the identity of~$G$, we get an embedding
$\gg\to\Vect(U)$, $X\mapsto X_G|_U$, and these vector fields can be 
expressed in terms of these coordinates, i.e.~they can be written as some 
differential operators in $n$ variables with analytic coefficients. 
Basically we obtain an embedding of~$\gg$ 
into some completion of the Weyl algebra in $2n$ generators. 
Of course, this embedding depends on the 
coordinate chart; a natural choice would be to take the chart given by 
the exponential map 
$\exp:\gg\to G$. The embedding thus obtained, when written in coordinates, 
turns out to be exactly the one defined by the Main Formula of this paper. 
However, we would like to deduce such a formula over 
any ring containing $\bbQ$, 
where such complex-analytic arguments cannot be 
used. We proceed by replacing in this argument all Lie groups and complex 
manifolds by their analogues in formal geometry 
-- namely, formal groups and formal schemes. 
This determines the layout of the next three sections.

Section~\ref{s:Formal} is dedicated to some generalities 
on formal schemes. Most notions and notations are variants 
of those developed in SGA~3 for the case of group schemes. 
We develop a similar functorial formalism for formal schemes, 
suited for our purpose, without paying too much attention 
to representability questions. 
Our exposition differs from most currently used approaches since we never 
require our rings to be noetherian. Besides, this section contains the 
construction of a formal group with given Lie algebra as well as some 
computations of tangent spaces and vector fields.

Section~\ref{s:WeylAlg} contains some generalities on Weyl algebras and 
completed Weyl algebras; the aim here is to develop some invariant 
descriptions 
of well-known mathematical objects, valid over any commutative base ring 
and for any finitely generated projective module. Besides, we establish 
isomorphisms between Lie algebras of derivations of (completed) symmetric 
algebras and some Lie subalgebras of the (completed) Weyl algebras.

Finally, in Section~\ref{s:EndOfProof} we deal with some questions of 
pro-representability, and use the results of the previous two sections 
to obtain and prove an invariant version of the Main Formula 
(cf.~\refpt{p:MainFormula}). Of course, explicit computation in terms of 
a chosen base of the Lie algebra $\gg$ and its structural constants 
gives us again the formula already proved in the first part of 
this paper by other methods.

So we proceed with an exposition of our functorial approach to formal schemes.

\nxpoint\label{p:pushouts}
Fix a base ring $\genrg \supset \bbQ$.
We will use the category $\calP$ defined as follows:
$$\begin{array}{l}
\Ob \calP := \{ (R,I) \,|\,\mbox {$I\subset R$ -- a nilpotent ideal in  
a commutative ring $R$}\}\\
\Hom_\calP\bigl((R,I), (R',I')\bigr) := 
\{\mbox{ring homomorphisms }\phi : R\to R'\,|\,\phi(I) \subset I'\}.
\end{array}$$
Category $\calP$ has pushouts (amalgamated sums) constructed as follows:
$$\xymatrix{
(R,I)\ar[r]\ar[d] & (R',I')\ar[d] \\
(R'',I'')\ar[r] & \bigl(R'\otimes_R R'', I'(R'\otimes_R R'')
+ I''(R'\otimes_R R'')\bigr)
}$$
There is a distinguished object $\genrg=(\genrg,0)$ in $\calP$. We will
consider the category $\relPk$ of morphisms in $\calP$ with source $\genrg$.
Sometimes we denote $(R,I)$ simply by $R$ and $I$ is assumed; 
then we denote $I$ by $I_R$. When we write tensor products in $\relPk$, they
are usually understood as amalgamated sums just described, i.e.\ 
$I_{R'\otimes R''}=I_{R'}\cdot(R'\otimes R'')+I_{R''}\cdot(R'\otimes R'')$.

\nxpoint
Our basic category of interest 
will be $\calE := \functcat(\relPk,\setscat)$ whose objects 
will be called presheaves (of sets) and
should be viewed as `formal varieties'; 
${\rm Grp}(\calE)=\calE^{\rm Grp} := \functcat(\relPk,{\rm Grp})$
is the category of group objects in $\calE$.
Its objects should be viewed as
presheaves of groups or `formal group schemes'. 

For any $(R',I')$ denote by $\Spf(R',I')$ the corresponding
representable functor $(R,I)\mapsto \Hom_\calP((R',I'),(R,I))$;
if $I' = 0$, $\Spf(R',I')$ is also denoted by $\Spec R'$.

Consider the following examples of rings and groups in~$\calE$
that might clarify the relationship 
with other approaches to formal groups:
$$\begin{array}{ll}
\uuO : (R,I) \mapsto R & \mbox{ (a ring in $\calE$)}\\
\widehat{\bbG}_a : 
(R,I)\mapsto I &\mbox{ (a group in $\calE$)}\\
\widehat{\bbG}_m : (R,I)\mapsto (1 + I)^\times 
&\mbox{ (a group in $\calE$)}
\end{array}$$
Only the first of these examples will be used in the sequel.

For any $\genrg$-module $M$ consider the $\uuO$-modules
$$
\uW(M): (R,I)\mapsto R\otimes_{\genrg} M,\quad
\uW(M)\supset \uWom(M) : (R,I)\mapsto I\cdot(R\otimes_{\genrg} M)
$$
If $M$ is free or projective, $\uWom(M):(R,I)\mapsto I\otimes_{\genrg} M$;
in this case one should think of $\uW(M)$ as 
``vector space $M$ considered as a manifold'', and
of $\uWom(M)\subset \uW(M)$ as ``formal neighborhood of zero in
$\uW(M)$''.

\nxpoint Given a morphism $\genrg\stackrel\pi\to \genrg'$ in $\calP$, 
it induces a functor 
$\pi_! : \relPkp\to \relPk$,
$\phi\mapsto \phi\circ\pi$, hence the restriction or base change functors 
$\pi^* : \calE = \functcat(\relPk, \setscat)
\to  \functcat(\relPkp, \setscat)=:\calE_{R_1}$,
$F\mapsto F\circ \pi_!$. Functor $\pi^*(F)$ is usually denoted by 
$F|_{\genrg'}$ or
$F_{(\genrg')}$. Functor $\pi^*:\calE\to\calE_{\genrg'}$ 
is exact and has a right adjoint 
$\pi_* : \calE_{\genrg'}\to\calE$ computed as follows:
for any $F:\relPkp\to \setscat$ we define 
$\pi_*F : \relPk\to\setscat$ by 
$\pi_* F : R\mapsto F(R\otimes_{\genrg}\genrg')$; tensor products are 
understood as coproducts in $\calP$ (as explained above in \refpt{p:pushouts});
$\_\mapsto\_\otimes_{\genrg}\genrg'$ is a functor
$\relPk\to \relPkp$.

\textbf{Notation.}  $\pi_* F$ is usually denoted by $R_{\genrg'/\genrg}(F)$
({\em ``Weil scalar restriction''}) or $\prod_{\genrg'/\genrg} F$.

\nxpoint Projective limits (e.g.\ direct products, fibered products,
kernels) are computed in $\calE$, $\calE^{\rm Grp}$ etc.\ 
componentwise, e.g.\ $F\times G : R\mapsto F(R)\times G(R)$, 
for $F,G \in \Ob\calE$. 
Category $\calE$ is in fact a closed Cartesian category. 
In particular, it has inner homs:
given $F,G \in \Ob\calE$, the inner hom is the presheaf 
$$\iHom(F,G): R \mapsto \Hom_{\calE_{R}}(F|_R, G|_R)$$
with the characteristic property
$\Hom_\calE(F,\iHom(G,H)) \cong \Hom_\calE(F\times G,H)$.
There are also canonical maps $\iHom(F,G)\times F\to G$,
$\iHom(F,\iHom(G,H))
\cong\iHom(F\times G,H)$ and so on. 
If $F$ and $G$ are presheaves of groups, one defines similarly
$\iHom_{\rm Grp}(F,G)\subset \iHom(F,G)$. 
There is an obvious subfunctor
$\iIsom(F,G)\subset \iHom(F,G)$ and
the special cases $\iEnd(F):=\iHom(F,F)$
and $\iAut(F):=\iIsom(F,F)$.

Functor $\iAut(F)$ is a group in $\calE$; any group homomorphism
$G\stackrel\rho\to\iAut(F)$ gives a group
action $\nu : G\times F\to F$ and conversely. 

We have the global sections functor $\Gamma:\calE\to \setscat$,
$F\mapsto F(\genrg)$ which is exact and satisfies $\Gamma(\iHom(F,G))
\cong \Hom_\calE (F,G)$.

\nxpoint (Tangent spaces)

Consider the dual numbers algebra $\genrg[\epsilon] = \genrg[T]/(T^2)$. 
In our setup, we have in fact two different versions of this: 
$\genrg[\epsilon]:=(\genrg[\epsilon],(\epsilon))$
and $\genrg[\epsilon]^\omega := (\genrg[\epsilon],0)$.
Notice also the following canonical morphisms in $\calP$:
$\genrg[\epsilon]^\omega \stackrel\nu\rightarrow \genrg[\epsilon]$,
inclusions $\genrg\stackrel{i}\hookrightarrow \genrg[\epsilon]^\omega$,
$\genrg\hookrightarrow \genrg[\epsilon]$ and the projections
$p : \genrg[\epsilon]^\omega$ and 
$\genrg[\epsilon]\to \genrg$.

For any functor $F\in \Ob\calE$ consider two new functors 
$TF := \prod_{\genrg[\epsilon]/\genrg} (F|_{\genrg[\epsilon]})$ and
$T^\omega F:=\prod_{\genrg[\epsilon]^\omega/\genrg}(F|_{\genrg[\epsilon]^\omega})$.
Then $TF : R\mapsto F(\genrg[\epsilon]\otimes_{\genrg} R)$,
$T^\omega F := R\mapsto F(\genrg[\epsilon]^\omega \otimes_{\genrg} R)$, or, more
precisely, $TF : (R,I)\mapsto F((R[\epsilon],I+R\epsilon))$,
$T^\omega F : (R,I) \mapsto F((R[\epsilon], I+I\epsilon))$.

{\bf Notation:} $R[\epsilon]:= \genrg[\epsilon]\otimes_{\genrg} R$,
$R[\epsilon]^\omega := \genrg[\epsilon]^\omega \otimes_{\genrg}R$ where 
tensor products are understood as in \refpt{p:pushouts}.

There is a canonical map $T^\omega F\stackrel{\nu_*}\to TF$,
as well as maps $\pi:= p_* : TF \to F$ and
$s := i_* : F\to T^\omega F\to TF$.
(One should think of $F$ as a ``manifold'', $TF$ -- its tangent bundle,
$\pi : TF\to F$ its structural map, $s : F\to TF$ -- zero section,
$T^\omega F\hookrightarrow TF$ -- formal neighborhood
of zero section in $TF$). 

If $G$ is a group, then $TG$ and $T^\omega G$ are also groups 
(in $\calE$) and $\pi$, $\nu$, $s$ -- group homomorphisms;
we define $\uLie(G)$ (resp.\ $\uLie^\omega(G)$)
to be the kernel of $TG\stackrel\pi\to G$ 
(resp.\ $T^\omega G\stackrel\pi\to G$):
$$\xymatrix{
0 \ar[r] &\uLie^\omega(G)\ar[r]\ar[d] & 
T^\omega G \ar[r]^\pi \ar[d] & G\ar[r]\ar@{=}[d]\ar@/^/[l]^s& 0
\\
0 \ar[r] &\uLie(G)\ar[r] & T G \ar[r]^\pi & G\ar[r]
\ar@/^/[l]^s& 0
}$$
$\uLie^\omega(G), \uLie(G)$ are clearly 
$\uuO$-modules; if $G$ is a ``good'' presheaf of groups (cf.\ SGA~3 I),
e.g.\ pro-representable, $\uLie^\omega(G)$ and $\uLie(G)$
have natural $\uuO$-Lie algebra structures in $\mathcal E$.

Put $\gg := \Gamma\uLie(G)$; this is a Lie algebra over $\genrg$, hence
by adjointness of functors $\uW\dashv \Gamma$ 
we have a canonical map $\uW(\gg)\to \uLie(G)$.
In most interesting situations this map is an isomorphism, and it 
maps $\uWom(\gg)\subset\uW(\gg)$ into $\uLie^\omega(G)\subset\uLie(G)$; 
then we identify $\uLie(G)$ with $\uW(\gg)$ and 
 $\uLie^\omega(G)$ with $\uWom(\gg)$.

\nxpoint For any $F\in \Ob\calE$ we define
$\uVect(F) := \iHom_F(F,TF)
: R\mapsto \{ \phi \in \Hom_\calE(F|_R,TF|_R)
\colon \pi|_R\circ\phi = \id_{F|_R}\}$
(sections of $TF$ over $F$); and similarly 
$\uVect^\omega(F):= \iHom_F(F,T^\omega F) \to \uVect(F)$; elements of 
$\uVect(F)(R)$ are ``vector fields'' on $F$ defined over $R$.

\nxpoint\label{p:vectfgrp} Short exact sequence of groups
$$\xymatrix{
0 \ar[r] &\uLie(G)\ar[r]& 
TG \ar[r]^\pi & G\ar[r]\ar@/^/[l]^s& 0
}$$
splits; so any $\bar{x} \in TG(R)$ can be written in form
$s(g)\cdot X$ where 
$X \in \uLie(G)(R)$ and
$g \in G(R)$; this decomposition is unique since
necessarily $g = \pi(\bar{x})$, $X = s(g)^{-1}\cdot\bar{x}$.
This splitting gives us an isomorphism 
$TG \stackrel\sim\to G\times \uLie(G)$,
and similarly 
$T^\omega G \stackrel\sim\to G\times \uLie^\omega(G)$;
in interesting situations 
$\uLie(G)\cong \uW(\gg)$ and $\uLie^\omega(G)\cong \uWom(\gg)$,
hence $TG \cong G\times \uW(\gg)$
and $T^\omega G\cong G\times\uWom(\gg)$.

Recall that $TG$ is a group, hence $G$ acts on $TG$ (say, from the left) 
by means of $G\stackrel {s}\to TG$; hence $G$ acts also on 
$\uVect(G) = \iHom_G(G,TG) \cong \iHom_G (G, G\times \uW(\gg))
\cong \iHom(G,\uW(\gg))$.
Constant maps of $\Hom(G, \uW(\gg))$
correspond to left-invariant vector fields under this identification. 

Right-invariant vector fields give us another isomorphism
$TG\stackrel\sim\to\uLie(G)\times G$,
corresponding to decomposition $\bar{x} = X\cdot s(g)$.

\nxpoint (Tangent spaces of $\uWom(M)$ and $\uW(M)$) 
Let $M$ be a $\genrg$-module. Then 
for any $R =(R,I)\in \Ob(\relPk)$ we have
$$\begin{array}{rcl}
\uWom(M)(R)
&=& I\cdot(R\otimes_{\genrg} M)\subset \uW(M)(R)
= R\otimes_{\genrg} M\\
\uW(M)(R[\epsilon]^\omega)
&=& \genrg[\epsilon]\otimes_{\genrg} M = \uW(M)(R)\oplus
\epsilon\cdot\uW(M)(R)\\
\uWom(M)(R[\epsilon])
&=& I\cdot(R\otimes_{\genrg}M)
+\epsilon R\otimes_{\genrg} M 
= \uWom(M)(R)\oplus\epsilon\cdot
\uW(M)(R)\\
\uWom(M)(R[\epsilon]^\omega)&=&
\uWom(M)(R)\oplus\epsilon\cdot\uWom(M)(R)
\end{array}$$
From this we get the following four {\em split\/} 
exact sequences of abelian groups in $\calE$:
$$\begin{array}{l}
0\to \uW(M)\stackrel{\cdot\epsilon}\to T\uW(M) \to
\uW(M)\to 0\\
0\to\uW(M)\stackrel{\cdot\epsilon}\to 
T^\omega\uW(M)\to
\uW(M)\to 0\\
0\to \uW(M)\stackrel{\cdot\epsilon}\to T\uWom(M) \to
\uWom(M)\to 0\\
0\to \uWom(M)\stackrel{\cdot\epsilon}\to T^\omega\uWom(M) \to
\uWom(M)\to 0
\end{array}$$
We deduce from these sequences canonical isomorphisms\\
$T\uWom(M)\cong\uWom(M)\times \uW(M)$,
$T^\omega \uW(M)\cong \uW(M)\times \uW(M)$ and so on. 

\nxpoint (Formal groups with given Lie algebra; exponential map)

{\bf From now on we assume} $\genrg\supset \bbQ$. Let $\gg$ be a Lie
algebra over $\genrg$, projective and finitely generated as
an $\genrg$-module (free of finite rank
suffices for most applications). Denote by $\Uenv(\gg)$ 
its universal enveloping $\genrg$-algebra, and by $\Uenv_i(\gg)$
its increasing filtration. The PBW theorem implies that
$\Uenv_i(\gg)/\Uenv_{i-1}(\gg)
\cong S^i(\gg)$. 
The diagonal map $\gg\to \gg\otimes\gg$
induces the comultiplication $\Delta : \Uenv(\gg)
\to \Uenv(\gg\oplus\gg)
\cong \Uenv(\gg)\otimes_{\genrg}\Uenv(\gg)$,
and the map $\gg\to 0$ induces the counit
$\counit : \Uenv(\gg)\to\Uenv(0) = \genrg$.
Thus $\Uenv(\gg)$ is a bialgebra and even a Hopf algebra.
For any $(R,I)/\genrg$ we have a canonical isomorphism
$\Uenv(\gg)_{(R)} = R\otimes_{\genrg} 
\Uenv(\gg) \stackrel\sim\leftarrow 
\Uenv\bigl(\gg_{(R)}\bigr)$, and $I\cdot \Uenv(\gg)_{(R)}$
is a nilpotent two-sided ideal in this ring. Now consider two formal groups
$$\begin{array}{rl}
&\Expp(\gg): (R,I)\mapsto
\{\alpha \in I\cdot \Uenv(\gg)_{(R)}\,|\,
\Delta(\alpha) = 1\otimes\alpha +\alpha \otimes 1, \,\counit(\alpha)=0\}\\
G := &\Expt(\gg)
: (R,I) \mapsto \{\alpha \in 1+I\cdot \Uenv(\gg)_{(R)}\,|\,
\Delta(\alpha) = \alpha\otimes\alpha, \,\counit(\alpha)=1 \}
\end{array}$$
Group operation on $\Expt(\gg)$ is induced
by the multiplication of $\Uenv(\gg)_{(R)}$,
and that of $\Expp(\gg)$ 
is determined by the requirement
that $\exp :\Expp(\gg)
\to \Expt(\gg)$ be a {\em group\/} isomorphism
where $\exp : \alpha \mapsto \sum_{n\geq 0}\frac{\alpha^n}{n!}$ is defined
by the usual exponential series; it makes sense for 
$\alpha \in I\cdot \Uenv(\gg)_{(R)}$ since
$I\cdot \Uenv(\gg)_{(R)}$ 
is a nilpotent ideal and $\bbQ\subset R$;
it is well known ([Bourbaki], Ch.~II) that 
$\exp_{(R,I)} : \Expp(\gg)(R,I)
\stackrel\sim\to \Expt(\gg)(R,I)$
is an isomorphism (i.e.\ it is bijective -- inverse is
given by the ${\rm log}$, and
$\Delta(\alpha) = 1\otimes \alpha + \alpha \otimes 1$, 
$\counit(\alpha) = 0$
iff $\Delta(\exp(\alpha))= \exp(\alpha)\otimes\exp(\alpha)$,
$\counit(\exp(\alpha)) = 1$).
Recall ([Bourbaki], Ch.~II) that the {\bf Campbell--Hausdorff series}
is a formal Lie power series $H(X,Y)$ in two variables $X$, $Y$ 
with rational coefficients, defined by the formal equality
$\exp(H(X,Y)) = \exp(X)\exp(Y)$
in the Magnus algebra 
$\Ahat_\bbQ(X,Y) = \hat{U}_\bbQ(L(X,Y))$. 
This implies that the group law on 
$\Expp(\gg)(R,I)$
is given by $H(X,Y)$, i.e.\ $\alpha\star \beta
= H(\alpha,\beta)$. [Since $\alpha$ and $\beta$ lie in a nilideal,
$H(\alpha,\beta)$ is a finite sum.]

There is also a canonical map $\nu : \uWom(\gg)
\to \Expp(\gg)$ defined as follows:
$\nu_R: \uWom(\gg)(R) = I_R\cdot(R\otimes_{\genrg} \gg)
= I_R\cdot \gg_{(R)} \to \Expp(\gg)(R)
= \{\alpha \in I\cdot \Uenv(\gg)_{(R)}\,|\,
\Delta(\alpha) = 1\otimes \alpha + \alpha \otimes 1,\, \counit(\alpha) = 0\}$
is induced by the canonical embedding 
$\gg_{(R)}\hookrightarrow \Uenv(\gg_{(R)})
\cong \Uenv(\gg)_{(R)}$. 
Actually, $\nu$ is an {\em isomorphism}:

In [Bourbaki], II, it is proved that for $R\supset\bbQ$ all
primitive elements of $\Uenv(\gg)_{(R)}$ come 
from $\gg_{(R)}$. By the PBW theorem, 
$\Uenv(\gg_{(R)})/\gg_{(R)}$ is flat, hence 
$I\cdot \Uenv(\gg_{(R)})\cap \gg_{(R)} = I\cdot\gg_{(R)}$, 
hence all $\nu_R$ are isomorphisms, i.e.~$\nu$ is an isomorphism:
$$\xymatrix{
\uWom(\gg) \ar[rr]^\nu_\sim \ar[rd]^\sim_{\exp'}
&& \Expp(\gg) \ar[ld]^\exp_\sim \\
& G = \Expt(\gg) &
}$$
We put $\exp':=\exp\circ\nu:\uW^\omega(\gg)\stackrel\sim\to G=\Expt(\gg)$ and 
call $\exp'$ {\em the exponential map\/} for $\gg$. 
Note that in the diagram above only $\exp$ is a {\em group\/} isomorphism.

\nxpoint (Lie algebra and tangent space computations for the
exponential map)
Now we are going to compute the maps $T(\nu)$, $T(\exp)$ and $T(\exp')$
from the previous diagram as well as 
$T(\uWom(\gg))$, $T(\Expp(\gg))$ and 
$T(\Expt(\gg))$.
In particular, we shall prove 
$\uLie(G)\cong \uW(\gg)$ and $\uLie^\omega(G)\cong \uWom(\gg)$.

\nxsubpoint
We know $T(\uWom(\gg))\cong\uWom(\gg)\times\uW(\gg)$ and
$T^\omega(\uWom(\gg))\cong\uWom(\gg)\times \uWom(\gg)$ 
with the first projection as the structural map 
$T(\uWom(\gg))\to\uWom(\gg)$ (resp\dots).

\nxsubpoint We compute $\uLie(\Expp(\gg))(R)$ by
definition (here $R = (R,I)$):
$$\begin{array}{lcl}\uLie(\Expp(\gg))(R)& =& 
\Ker\bigl(\Expp(\gg)(R[\epsilon])\to\Expp(\gg)(R)\bigr)
=\\ &=& \{ \alpha\epsilon\in \Uenv(\gg)_{(R[\epsilon])}
\,|\,\Delta(\alpha\epsilon) 
= \alpha\epsilon\otimes 1 + 1\otimes \alpha\epsilon,
\counit(\alpha\epsilon) = 0\}\\ &\cong& \uW(\gg)(R)\quad,\end{array}$$
hence $\uLie(\Expp(\gg))=\uW(\gg)$,
and similarly $\uLie^\omega(\Expp(\gg)) = \uWom(\gg)$.

For any formal group, hence for $\Expp(\gg)$, we have (cf.\ \refpt{p:vectfgrp})
$T(\Expp(\gg))\cong 
\Expp(\gg)\times\uLie(\Expp(\gg))\cong 
\Expp(\gg)\times\uW(\gg)$, and similarly for $T^\omega(\Expp(\gg))$.

\nxsubpoint Since $\exp:\Expp(\gg)\to G=\Expt(\gg)$ is
an isomorphism, we can expect similar descriptions for
$\uLie(G)$ and $\uLie^\omega(G)$. One can also check directly 
$\uLie(G)(R)=\{1+\alpha\epsilon\in\Uenv(\gg)_{(R[\epsilon])}
\,|\,\Delta(1+\alpha\epsilon)=(1+\alpha\epsilon)\otimes (1+\alpha\epsilon),
\quad\counit(1+\alpha\epsilon) = 1 \} \cong \uW(\gg)(R)$ and 
similarly for $\uLie^\omega(G)$.
Note that 
$\uLie(\exp) :\uLie(\Expp(\gg))\to \uLie(\Expt(\gg))$,
$\alpha \epsilon \mapsto 1 +\alpha\epsilon$ is an isomorphism
even if $\genrg$ does not contain $\bbQ$. If we identify 
$\uLie(\Expp(\gg))$ and $\uLie(G)$ with
$\uW(\gg)$, $\uLie(\exp)$ is
identified with the identity map, and 
similarly for $\uLie^\omega$; hence we get
\begin{equation}\label{eq:lieexpsix}  
\xymatrix{
T(\Expp(\gg))\ar[rr]^\sim \ar[d]_{T(\exp)}
&&
\Expp(\gg)\times \uW(\gg)  \ar[r]   
\ar[d]_{\exp\times\id}
&
\Expp(\gg) \ar[d]^\exp
\\
TG \ar[rr]^\sim
&&
G\times \uW(\gg)\ar[r]
& G
}
\end{equation}

\nxsubpoint\label{p:LieBr} One checks that our identifications
$\uLie(\Expp(\gg))\cong \uLie(G)\cong \uW(\gg)$
and their analogues for $\uLie^\omega$ are compatible
with the original Lie algebra structure on $\gg$. Indeed,
consider two elements $1+X\epsilon$, $1+Y\eta$ in $G\bigl(R[\epsilon,\eta]/
(\epsilon^2,\eta^2)\bigr)$; then by definition
$$\begin{array}{lcl}
1 + [X,Y]_{\uLie(G)}\cdot \epsilon\eta & = &
(1+X\epsilon)(1+Y\eta)(1+X\epsilon)^{-1}(1+Y\eta)^{-1}
\\ &=&
  (1 + X\epsilon +Y\eta +XY\epsilon\eta)
(1+X\epsilon + Y\eta +YX\epsilon\eta)^{-1}
\\ & = &  1 + (XY-YX)\epsilon\eta
= 1 + [X,Y]_\gg\cdot\epsilon\eta.
\end{array}$$

\nxsubpoint\label{p:diffofnu} Consider now $\nu :\uWom(\gg)
\to \Expp(\gg)$ (this is an isomorphism, but not an isomorphism of
groups):
\begin{equation}\label{eq:lieexpsix2}\xymatrix{
T(\uWom(\gg)) \ar[rr]^\sim\ar[d]_{T(\nu)}
&&
\uWom(\gg)\times\uW(\gg)
\ar[rr]\ar@{-->}[d]^{\nu'}
&&
\uWom(\gg)\ar[d]^\nu
\\
T(\Expp(\gg)) \ar[rr]^\sim
&&
\Expp(\gg)\times\uW(\gg)
\ar[rr]^<>(0.5)\pi
&&
\Expp(\gg)\ar@/_1pc/@{.>}[ll]_s
}\end{equation}
In this diagram we identify $\uLie(\Expp(\gg))$ with
$\uW(\gg)$ as explained above.
We want to compute the map $\nu'$.
Let $R = (R,I)$, 
$X \in \uWom(\gg) = I\otimes_{\genrg}\gg$ and
$Y \in \uW(\gg)(R) = R\otimes_{\genrg}\gg$; 
the corresponding element of 
$T(\uWom(\gg))(R) =\uWom(\gg)(R[\epsilon])$ 
is given by $X+Y\epsilon$, and
$T(\nu)_R(X+Y\epsilon) = \nu_{R[\epsilon]}(X+Y\epsilon) = X + Y\epsilon$
considered as a primitive element of $\Uenv(\gg)_{(R[\epsilon])}$.
Clearly $\pi_R(X+Y\epsilon)=X$, $s_R\pi_R(X+Y\epsilon)=X$, 
and we want to find 
$Z \in \uLie(\Expp(\gg))(R) = \uW(\gg)(R)$,
such that $(X+Y\epsilon) = X \star Z\epsilon$ inside
$T(\Expp(\gg))(R)$. Since $(-X)\star X = 0$, we have
$Z\epsilon = (-X)\star (X+Y\epsilon)$; classical formula for
$H(X+Y,-X) \bmod \deg_Y 2$ (or for $H(-X, X+Y)$) 
([Bourbaki], II 6.5.5) gives us 
$Z = \sum_{n\geq 0} \frac{(-\ad X)^n}{(n+1)!} (Y)$. 
Hence $\nu'$  is given by
$\nu'_R : (X,Y)\mapsto (X,Z)$ with 
$Z$ defined by the above formula.

{\bf Remark.} If we consider the other 
canonical splitting of $T(\Expp(\gg))$ by
right-invariant vector fields, we obtain almost the same formula
for $Z$, but without the $(-1)^n$ factors:
$Z = \sum_{n\geq 0} \frac{(\ad X)^n}{(n+1)!} (Y)$. 

\nxsubpoint\label{p:compTexp} We now consider the exponential map 
$\exp' = \exp\circ\nu : \uWom(\gg)\to G$.
By composing
~\eqref{eq:lieexpsix} and~\eqref{eq:lieexpsix2} we get
the commutative diagram
\begin{equation}\xymatrix{
T(\uWom(\gg))\ar[rr]^\sim
\ar[d]^{T(\exp')}&&
\uWom(\gg)\times\uW(\gg)
\ar[rr] 
\ar@{-->}[d]^\tau&&
\uWom(\gg)\ar[d]^{\exp'}
\\
TG \ar[rr]^\sim &&
G\times \uW(\gg)
\ar[rr]^<>(.5)\pi
&&
G\ar@/_1pc/@{.>}[ll]_s}
\end{equation}
The middle vertical map is given by
$\tau:(X,Y)\mapsto(\exp(X),\sum_{n\geq 0}
 \frac{(-\ad X)^n}{(n+1)!}(Y))$.
One can check this directly, without referring to \refpt{p:diffofnu}
and the properties of Campbell--Hausdorff series. Indeed, if
$\tau(X,Y) = (\exp(X),Z)$ for some 
$X \in \uWom(\gg)(R)$,
$Y,Z\in \uW(\gg)(R)$,
we must have 
$\exp(X+Y\epsilon) = \exp(X)(1+Z\epsilon)$ inside $G(R[\epsilon])
\subset \Uenv(\gg)_{(R[\epsilon])}$.
Since $\epsilon^2 = 0$, we have
$(X+Y\epsilon)^n = X^n + \sum_{p+q = n-1}
X^p Y X^q \cdot \epsilon$, hence
$\exp(X)\cdot Z = \sum_{p,q\geq 0}
\frac{X^p Y X^q}{(p+q+1)!}$. One then checks directly that
$Z = \sum_{n\geq 0} \frac{(-\ad X)^n}{(n+1)!}(Y)$
satisfies this equality. 
The argument like the one in the proof of 
[Bourbaki], II, 6.5.5 is better: one checks almost
immediately that
$(\ad X)(Z) = (1- e^{-\ad X})(Y)$, and then divides formally
by $\ad X$, considering both sides as elements
of the completed free Lie algebra in $X$ and $Y$ over $\bbQ$.

\nxpoint (Formal completions) Suppose we are given 
morphism $H\stackrel\phi\to F$ in $\calE=\calE_{\genrg}$.
In most cases of interest $\phi$ will be a monomorphism, i.e. $H$ can be
identified with a subfunctor of $F$. 
Denote the natural map $(R,I)\to(R/I,0)$ by $\pi_{(R,I)}$.
We say that $H$ {\bf is complete over}
(or {\bf in}) $F$ if for any $(R,I)$ in 
$\relPk$ the following diagram is cartesian
\begin{equation}\label{eq:complcart}\xymatrix{
H(R,I) \ar[rr]^{H(\pi_{(R,I)})} \ar[d]^{\phi_{(R,I)}}&& H(R/I,0) 
\ar[d]^{\phi_{R/I}}\\
F(R,I) \ar[rr]^{H(\pi_{(R,I)})} && F(R/I,0)
}\end{equation}
Arbitrary morphism $H\stackrel\phi\to F$ may be factored
as $H\stackrel{u}\to\hat{F}_H\stackrel\kappa\to F$ where
$\hat{F}_H\stackrel\kappa\to F$ is complete, and this {\bf completion} 
$\hat{F}_H$ is universal in the sense that if
$F'\stackrel{\kappa'}\to F$ is complete and $H\stackrel{u'}\to F'$ 
is such that $\kappa'\circ u'=\phi$, then there is a unique map
$F'\stackrel{\chi}\to \hat{F}_H$ 
such that $\chi\circ u'=u$ and $\kappa\circ\chi=\kappa'$.
To this aim, define $\hat{F}_H(R,I)$ to be the fibered product
of $F(R,I)$ and $H(R/I,0)$ over $F(R/I,0)$:
\begin{equation}\label{eq:complconstr}\xymatrix{
\hat{F}_H(R,I) \ar[rr] \ar[d]^{\kappa_{(R,I)}}&& H(R/I,0) 
\ar[d]^{\phi_{R/I}}\\
F(R,I) \ar[rr] && F(R/I,0)
}\end{equation}
If $\phi :H\to F$ is monic, then $\kappa$ and $u$ are also.
$\hat{F}_H$ is called the {\bf formal completion of $F$ along $H$}.
If $F$ and $H$ are groups (or $\uuO$-modules)
and $H\to F$ a morphism of groups (resp.~$\uuO$-modules)
then $\hat{F}_H$ is also, 
and $H\stackrel{u}\to \hat{F}_H\stackrel\kappa\to F$
will be morphisms of such. 

{\bf Examples.} a) For any $\genrg$-module $M$, $\uWom(M)$ 
is the formal completion of $\uW(M)$ 
along $0\subset \uW(M)$. Indeed, the universality of%
~\eqref{eq:complconstr} with $F = \uW(M)$,
$\hat{F}_H = \uWom(M)$, $H=0$ is immediate.

b) Suppose $F$ is left-exact functor (e.g. pro-representable).
Then $T^\omega F$ is the completion of $TF$ 
along the zero section $s:F\to TF$.

c) If $G$ is a group in $\calE$, left exact as a functor,
then $\uLie^\omega(G)$ is
the completion of $\uLie(G)$ along $0$.

\section{Weyl algebras}\label{s:WeylAlg}
\nxpoint Let $\genfd$ be a commutative ring, $M$ a $\genfd$-module,
$\Phi : M\times M\to\genfd$ a bilinear form 
which is symplectic: $\forall x \in M$, $\Phi(x,x) = 0$ (this implies
$\forall x,y\in M$, $\Phi(x,y) = -\Phi(y,x)$). 
{\em We do not require the nondegeneration.}
Consider the category 
${\mathcal C}^\Phi_M$, objects of which are pairs $(A,\lambda_A)$ where
$A$ is an associative $\genfd$-algebra and $\lambda_A:M\to A$ is a 
$\genfd$-linear map, such that $\forall x,y\in M$,
$[\lambda(x),\lambda(y)]_A = \lambda(x)\lambda(y)-\lambda(y)\lambda(x)
= \Phi(x,y)\cdot 1$; morphisms $(A,\lambda_A)\to (B,\lambda_B)$ are just
$\genfd$-algebra homomorphisms $f: A\to B$ compatible with \hbox{$\lambda$-s}:
$\lambda_B = f\circ \lambda_A$. 

{\bf Definition.} The universal (i.e.~initial) object of the category 
${\mathcal C}^\Phi_M$ will be denoted by $(SW(M,\Phi),i_M)$ and
it will be called the {\bf symplectic Weyl algebra of} $(M,\Phi)$.

In particular, $SW(M)=SW(M,\Phi)$ is an associative $\genfd$-algebra
and $i_M : M\to SW(M)$ is a $\genfd$-linear map such that
$\forall x,y\in M$, $[i_M(x),i_M(y)]=\Phi(x,y)$.

Of course, $SW(M)$ exists, it may be constructed as
a quotient of tensor algebra, namely $T(M)/I_\Phi$, where 
$T(M) = \oplus_{n\geq 0} M^{\otimes n}$ is the tensor algebra of $M$, and
$I_\Phi\subset T(M)$ is the two-sided ideal of $T(M)$ generated
by the elements of the form $x\otimes y - y \otimes x -\Phi(x,y)\cdot 1$
for all $x,y \in M$. 

Note that if $\{x_\alpha\}$ is a system of generators of $M$, 
by bilinearity
it is sufficient to require $[i_M(x_\alpha), i_M(x_\beta)] = 
\Phi(x_\alpha,x_\beta)$ in the definition of ${\mathcal C}^\Phi_M$ and
$SW(M)$, and $I_\Phi$ is generated by 
$x_\alpha\otimes x_\beta - x_\beta\otimes x_\alpha - 
\Phi(x_\alpha,x_\beta)$.

{\bf Remark.} Symplectic Weyl algebras are very similar to Clifford 
algebras, and, in some other respects, 
to the universal enveloping algebras
of Lie algebras. 

\nxpoint (Basic properties of symplectic Weyl algebras)

a) $SW(0) = \genfd$, $SW(M,0)$ is the symmetric algebra $S(M)$.

b) (functoriality) Given a 
$\genfd$-linear map $f : (M,\Phi) \to (N,\Psi)$ 
with $\Psi\circ (f\times f) = f\circ \Phi$ there is a unique
map $SW(f) : SW(M,\Phi)\to SW(N,\Psi)$ satisfying the
obvious conditions.

c) If $M = \genfd e$ is free of rank 1, 
then necessarily $\Phi = 0$ and $SW(M)\cong S(M)\cong \genfd[e]$.

d) Functor $SW:(M,\Phi)\mapsto SW(M,\Phi)$ 
commutes with filtered inductive limits in $(M,\Phi)$.

e) If $M$ is an orthogonal direct sum of $M_1$ and $M_2$:
$M = M_1 \oplus M_2$ and $\Phi = \Phi_1\oplus \Phi_2$, then 
the canonical map $SW(M_1)\otimes_\genfd SW(M_2)\to SW(M_1\oplus M_2)$
(induced by canonical maps $SW(M_i)\to SW(M_1\oplus M_2)$
and multiplication in $SW(M_1\oplus M_2)$)
is an {\em isomorphism of ${\genfd}$-algebras} (this is immediate
by checking that  $SW(M_1)\otimes_\genfd SW(M_2)$ satifies the
universal property). 

f) $SW(M,\Phi)$ commutes with base change:\\ 
For any $\genfd'/\genfd$
we have $SW(M_{(\genfd')},\Phi_{(\genfd')}) = SW(M,\Phi)_{(\genfd')}$.

\nxpoint (Filtration and $\bbZ/2\bbZ$-grading)
Consider the natural grading of the tensor algebra 
$T(M) = \oplus_{n\geq 0} T^n(M) :=  \oplus_{n\geq 0} M^{\otimes n}$
and the corresponding increasing filtration  $T_n(M)=T_{\leq n}(M) 
:=\oplus_{m\leq n} T^m(M)$; clearly $T_0(M) = \genfd$,
$T_1(M) = \genfd\oplus M$. $T(M)$ has canonical supergrading:
$T^+(M) := \oplus_n T^{2n}(M)$,
$T^-(M) := \oplus_n T^{2n+1}(M)$. Since $I_\Phi \subset T(M)$ is
generated by even elements $x\otimes y- y\otimes x -\Phi(x,y)\in T^+(M)$,
$SW(M):= T(M)/I_\Phi$ inherits a canonical supergrading:
$SW^{\pm}(M) = \pi(T^\pm (M))$, $SW(M) = SW^+(M)\oplus SW^-(M)$
where $\pi :T(M)\to SW(M)$ is the projection.
There is also the canonical filtration (sometimes called ``Bernstein''
or ``arithmetic'') on $SW(M)$ 
-- the image under $\pi$ of the filtration of
$T(M)$, namely $SW_n(M) := \pi(T_n(M))$.
We have a canonical surjective map $\gr(\pi): T(M)=
\gr_F(T(M))\to \gr_F(SW(M))$.
Clearly, $\gr_F(SW(M))$ is generated by the image of $M$ in
$\gr ^1_F(SW(M))$, and for any two elements
$x,y\in M\subset T_1(M)$, their images $\bar{x}$, $\bar{y}$ in $\gr^1_F(SW(M))$
commute since $xy-yx \in SW_1(M)$. 
This means that $\gr(\pi)$ 
factorizes through $T(M)\to S(M)$, 
hence we obtain surjective maps
$S(M)\stackrel\phi\to \gr_F(SW(M))$
and $S^n(M)\stackrel{\phi^n}\to \gr_F^n(SW(M))= 
SW_n(M)/SW_{n-1}(M)$.

\nxpoint\label{p:SW-PBW} {\bf Theorem.} {\em 
If $M$ is a flat $\genfd$-module (e.g. free or projective), then
$\phi$ and
$\phi^n : S^n(M)\to SW_n(M)/SW_{n-1}(M)$
are isomorphisms for all $n\geq0$.
}

{\bf Proof.} a) First assume $M$ free with a base 
$(e_\alpha)_{\alpha\in \Lambda}$. By Zermelo lemma 
we may assume that $\Lambda$ is well-ordered set. 
Denote by $x_\alpha$ and $z_\alpha$ the images of $e_\alpha$ in $SW(M)$
and $S(M)$ respectively. Consider the set of nondecreasing sequences
$\Seq = \{I = (i_1,\ldots,i_n)\,|\, i_k \in \Lambda,
i_1\leq i_2\leq \ldots \leq i_n\}$ with notation $|I| = n$
whenever $I = (i_1,\ldots, i_n)$. If $J = (i_2,\ldots, i_n)$,
we write $I = (i_1,J)$. Notation $\lambda \leq I$ means
$\lambda \leq i_k$, for all $k$, and similarly $\lambda < I$.
Now, for any $I = (i_1,\ldots, i_n)\in \Seq$ define
$x_I := x_{i_1}\ldots x_{i_n} \in SW(M)$,
$z_I := z_{i_1}\cdots z_{i_n} \in S(M)$.
Clearly, $\{z_I\}$ form a basis of $S(M)$,
$\{z_I\}_{|I|=n}$ form a basis of $S^n(M)$, and
$\{x_I\}$ generate a $\genfd$-submodule in $SW(M)$
 containing the image
of $M$ and closed under multiplication (to this aim
check by induction on $|J|$ that $x_\lambda x_{J}$ lies in this submodule,
and then by induction on $|I|$ that $x_I x_J$ lies there as well).
Consequently, it is sufficient to show that the $\{x_I\}$ are linearly 
independent. To prove this we proceed as in the proof of PBW theorem in
[Bourbaki], ch.~I. We construct by induction a family of compatible 
bilinear maps $\rho_n : M\times S_n(M)\to S_{n+1}(M)$ enjoying the
following properties:

$(A_n)$ $\rho_n(x_\lambda, z_I) = z_\lambda z_I$ if $\lambda \leq I$,
$|I|\leq n$ (or $I=\emptyset$); $\rho_n(x_\lambda,z_I)\equiv z_\lambda 
z_I\pmod{S_n(M)}$ if $|I|=n$ for any $\lambda\in \Lambda$;

$(B_n)$ $\rho_n(x_\lambda,z_I) = \rho_{n-1}(x_\lambda,z_I)$ 
if $|I|\leq n-1$;

$(C_n)$ $\rho_n(x_\lambda,\rho_{n-1}(x_\mu,z_I)) =\rho_n(x_\mu,
\rho_{n-1}(x_\lambda,z_I)) +\Phi(e_\lambda,e_\mu)\cdot z_I$ if
$\lambda,\mu\in\Lambda$, $|I|\leq n-1$.

We see that such a family is uniquely defined: if $\lambda\leq I$,
$\rho_n(x_\lambda,z_I)$ is defined by $(A_n)$;
if $|I|\leq n-1$ by $(B_n)$; if $I=(\mu,J)$, $|J|=n-1$ and $\lambda>\mu$
by $(C_n)$; then one checks directly that the maps $\rho_n$ so defined
do satisfy $(A_n)$, $(B_n)$ and $(C_n)$.

This map $\rho : M\times S(M)\to S(M)$ induces a map
$\tilde\rho : M\to \End_\genfd(S(M))$, such that
$[\tilde\rho(e_\lambda),\tilde\rho(e_\mu)] = \Phi(e_\lambda,e_\mu)$;
hence by universal property of Weyl algebras 
$\tilde\rho$ induces a map 
$\tilde\rho' : SW(M)\to \End_\genfd(S(M))$.
One sees that the map $SW(M)\to S(M)$ defined by
$x\mapsto \tilde\rho'(x)(1)$ maps $x_I$ into $z_I$, and $\{z_I\}$ are
linearly independent ; hence $\{x_I\}$ are also independent. 

b) If $M$ is flat, it can be written, 
by a classical result of {\sc Lazard}, 
as a filtered inductive limit of free modules: 
$M = \injlim M_\alpha$, where all $M_\alpha$ are free;
consider on each $M_\alpha$ the pullback $\Phi_\alpha$
of the form $\Phi$ with respect to the natural 
map $M_\alpha\to M$; then $SW(M) = \injlim SW(M_\alpha)$;
everything commutes with $\injlim$ and
$S^n(M_\alpha)\to \gr^n(SW(M_\alpha))$ are isomorphisms,
hence the same holds for $M$. 

\nxpoint (Consequences)

a) If $M$ is flat then all $SW_n(M)/SW_{n-1}(M)\cong S^n(M)$ are flat,
hence also $SW_n(M)$ (by induction on $n$) and
$SW(M) = \injlim SW_n(M)$ are flat.

b) If $M$ is projective, all $SW_n(M)/SW_{n-1}(M)\cong S^n(M)$ 
are projective,
hence $SW_{n-1}(M)$ has a complement $K_n$ in 
$SW_n(M)$; $K_n \cong SW_n(M)/SW_{n-1}(M)$ is projective, hence
$SW_n(M) = \oplus_{m\leq n} K_m$ and $SW(M) = \oplus_{m\geq 0} K_m$
are projective as well.

b') If $M$ is projective and finitely generated, the same can be said
about all $SW_n(M)/SW_{n-1}(M)$ and
all $SW_n(M)$ [ but {\em not\/} about $SW(M)$! ].

c) If $M$ is free then the same is true for
all $SW_n(M)/SW_{n-1}(M)$, all $SW_n(M)$ and for $SW(M)$.

c') If $M$ is free and of finite rank
then the same can be said about 
all $SW_n(M)/SW_{n-1}(M)$ and all $SW_n(M)$

d) If $M$ is flat, then maps $\genfd\to SW(M)$ and 
$M\to SW(M)$ are injective,
so $\genfd$ and $M$ can be identified with their
images in $SW(M)$; moreover, $\genfd\oplus M\to SW(M)$ is injective, 
and its image is $SW_1(M)$.

\nxpoint (Convolution) Fix a symplectic $\genfd$-module $(M,\Phi)$,
and let $i = i_M : M\to SW(M)$ be the canonical map.  
Given $u \in M^* = \Hom_\genfd (M,\genfd)$, the 
{\bf convolution map} $D_u : SW(M)\to SW(M)$ is a $\genfd$-derivation
restricting to $u$ on the image of $M$. In other words, we require
0) $D_u(\genfd\cdot 1) = 0$; 1) $\forall \alpha,\beta \in SW(M)$, 
$D_u(\alpha\cdot\beta) = D_u(\alpha)\cdot\beta +\alpha\cdot D_u(\beta)$;
2) $\forall x\in M$, $D_u (i(x)) = u(x)\cdot 1$. These conditions
imply 
$D_u(x_1 x_2\cdots x_n) = \sum_{i=1}^n u(x_i)\cdot x_1\cdots\widehat{x_i}
\cdots x_n$ for all $n \geq 1$ and for all $x_1,\ldots, x_n \in i(M)$.
This implies the uniqueness of $D_u$. To show the existence, consider
the dual numbers $SW(M)[\epsilon]:=\genfd[\epsilon]\otimes_\genfd SW(M)
= SW(M)\oplus SW(M)\epsilon$, and the $\genfd$-linear map
$i_u : M\to SW(M)[\epsilon]$ given by $i_u : x\mapsto i(x)+u(x)\epsilon$. 
One checks directly $[i_u(x),i_u(y)]=[i(x),i(y)]=\Phi(x,y)$, hence $i_u$ 
induces a $\genfd$-algebra homomorphism $\itilde_u: SW(M)\to SW(M)[\epsilon]$.
Now for any $x\in SW(M)$ we have $\itilde_u(x)=x+D_u(x)\epsilon$ for some $D_u(x)$, and this determines the convolution~$D_u$.

\nxpoint\label{p:SWdeform} (Deformation) 
Suppose $M = M_1 \oplus M_2$, but $M_1$
and $M_2$ are not necessarily orthogonal. Inclusions $M_i\to M$ induce
maps $SW(M_i)\to SW(M)$, hence (together with multiplication
in $SW(M)$) a map $\rho : SW(M_1)\otimes SW(M_2)\to SW(M)$.
In general, this is not a homomorphism of $\genfd$-algebras, 
but still a $\genfd$-linear map. 

{\bf Claim.} {\em If one of $M_1$ or $M_2$ is flat, then $\rho$ 
is an isomorphism of $\genfd$-module.}

{\bf Proof.} a) Suppose $M_2$ is flat; write it in form  
$M_2 = \injlim N_\alpha$ where all $N_\alpha$ are 
free of finite rank and consider
on $M\oplus N_\alpha$ the symplectic form induced from $M$ by 
$M_1\oplus N_\alpha\to M_1\oplus M_2 = M$. Since everything commutes
with inductive limits, it is sufficient to prove that 
$SW(M_1)\otimes SW(N_\alpha)\to SW(M_1\oplus N_\alpha)$
are isomorphisms, i.e.\ we can assume $M_2$ 
to be free of finite rank. An easy induction argument
shows that we can assume $M_2$ free of rank one. 

b) Assume $M_2 = \genfd\theta$ is free of rank 1, 
$M = M_1\oplus \genfd\theta$. Put $u(x) := \Phi(\theta, x)$
for any $x\in M_1$, so we get $u \in M_1^*$. 
Consider the action $\star$ of $M$ on 
$SW(M_1)\otimes_\genfd S(\genfd\theta)
= SW(M_1)\otimes_\genfd \genfd[\theta] = SW(M_1)[\theta]$:
$M_1$ acts on $SW(M_1)[\theta]$ by left multiplication
$x \star \sum_{j\geq 0} \alpha_j \theta^j = 
\sum_{j\geq 0}(x\alpha_j)\theta^j$
for any $x\in M_1$, $\alpha_j\in SW(M_1)$; and 
$\theta \star \sum_{j\geq 0} \alpha_j \theta^j = 
\sum_{j\geq 0} \alpha_j \theta^{j+1} + \sum_{j\geq 0} D_u(\alpha_j)\theta^j
= \sum_{j\geq 0} (\alpha_{j-1} + D_u(\alpha_j)) \theta^j$. So we get
a $\genfd$-linear map 
$\gamma : M = M_1\oplus \genfd\theta\to \End_\genfd(SW(M_1))$. It is
easy to see that $[\gamma(x),\gamma(y)] = \Phi(x,y)$;
if $x,y\in M_1$ or if $x = y = \theta$ this is evident; and in the mixed
case $x = \theta$, $y\in M_1$ we calculate
$$\begin{array}{l}
\gamma(\theta)\gamma(y)(\sum_{j\geq 0}\alpha_j \theta^j) = 
\theta\star(y\star\sum_{j\geq 0} \alpha_j \theta^j) = 
\theta\star \sum_{j\geq 0} (y\alpha_j)\theta^j \\
\quad\quad = \sum_{j\geq 0} y \alpha_j
\theta^{j+1} + \sum_{j\geq 0} D_u(y\alpha_j)\theta^j\quad,
\\
\gamma(y)\gamma(\theta)(\sum_{j\geq 0} \alpha_j\theta^j) = 
\sum_{j\geq 0} y\alpha_j \theta^{j+1} + 
\sum_{j\geq 0} y D_u(\alpha_j)\theta^j\quad,
\\
\left[\gamma(\theta),\gamma(y)\right](\sum_{j\geq 0} \alpha_j\theta^j) =
\sum_{j\geq 0} D_u(y)\alpha_j\theta^j = u(y)\cdot\sum_{j\geq 0}\alpha_j\theta^j
\quad,
\\
\left[\gamma(\theta),\gamma(y)\right] = u(y) = \Phi(\theta,y)\quad.
\end{array}$$
This means that $\gamma$ induces a map 
$\tilde\gamma: SW(M)\to \End_\genfd (SW(M)[\theta])$; 
now it is immediate that $\alpha\mapsto \tilde\gamma(\alpha)(1)$ gives a map
$SW(M)\to SW(M_1)\otimes S(\genfd\theta)$ inverse
to $\rho$, hence $\rho$ is an isomorphism. 

\nxpoint Note that if $M = \genfd e_1 \oplus \ldots \oplus \genfd e_n$, 
\refpt{p:SWdeform} implies that
$S(M)=S(\genfd e_1)\otimes\cdots\otimes S(\genfd e_n)\to SW(M)$ is 
an isomorphism; 
this proves~\refpt{p:SW-PBW} in this case; the general case
of~\refpt{p:SW-PBW} can be deduced from this by taking
inductive limits.

\nxpoint\label{p:SWbilin} 
Now suppose $M=Q\oplus P$ with flat and $\Phi$-isotropic $P$ and $Q$
(i.e.\ $\Phi|_{P\times P} = 0$ and $\Phi|_{Q\times Q} = 0$). 
Then $\Phi$ is uniquely determined by the bilinear form 
$\phi := \Phi|_{Q\times P}:Q\times P\to K$ since
$\Phi(q+p,q'+p') = \phi(q,p')-\phi(q',p)$, and any bilinear
form $\phi$ defines such a $\Phi$. Algebra
$\calD_{Q,P,\phi} = \calD_\phi := SW(Q\oplus P,\Phi)$
is called the {\bf Weyl algebra defined by $\phi$}.
Since $Q$ and $P$ are isotropic, $SW(Q)\cong S(Q)$ and 
$SW(P)\cong S(P)$; since they are flat, the map
$S(Q)\otimes_\genfd S(P) = SW(Q)\otimes_\genfd SW(P)\to SW(Q\oplus P,\phi)
=\calD_\phi$ is an isomorphism of $\genfd$-modules 
by~\refpt{p:SWdeform}.

Any $p\in P$ defines a form $d_\phi(p)\in Q^*$ by the rule
$d_\phi(p): q\mapsto \phi(q,p)$, hence a derivation
(=convolution) $D_p:=D_{-d_\phi(p)}$ on $S(Q) = SW(Q)$ 
(note the minus sign!). On the other hand,
any $q\in Q$ acts on $S(Q)$ by multiplication $L_q :\alpha\mapsto q\alpha$.
Since $[D_p,D_{p'}] = [L_q,L_{q'}]=0$
and $[D_p,L_q] = -\phi(q,p)$, we see that the map
$q+p\mapsto L_q+D_p$ defines a $\genfd$-linear map $\calD_\phi
=SW(Q\oplus P)\to \End_\genfd(S(Q))$,
i.e. a $\calD_\phi$-{\bf module structure on $S(Q)$}. Note
that $S(Q)\subset \calD_\phi$ acts on $S(Q)$ with respect
to this structure in the natural way, and $S(P)\subset \calD_\phi$
acts by convolutions, and, in particular, 
$P\subset S(P)\subset \calD_\phi$ by derivations,
hence $S(Q)\otimes P\subset \calD_\phi$ 
also acts on $S(Q)$ by $\genfd$-derivations. 

Notice that $S^r(Q) \cdot S^n(Q) \subset S^{r+n}(Q)$,
$S^r(P)\cdot S^n(Q)\subset S^{n-r}(Q)$. The construction is preserved
up to signs when we interchange $Q$ and $P$, hence
$\calD_\phi$ acts on $S(P)$ as well. 

\nxpoint\label{p:ComplWeyl} 
(Completed Weyl algebra and its action on the completed symmetric algebra) 
Consider the completed symmetric algebra 
$\Shat(Q) =\prod_{n\geq 0} S^n(Q) = 
\projlim S(Q)/S^{\geq n}(Q) = \projlim S(Q)/(S^+(Q))^n$. Clearly, 
$Q\subset S(Q)\subset \Shat(Q)$;
if $Q$ is a free $\genfd$-module of rank $n$,
then $S(Q)$ is the algebra of polynomials in $n$ variables,
and $\Shat(Q)$ the algebra of formal power series.

Recall that $S(Q)\otimes_\genfd S(P)\stackrel\rho\to \calD_\phi$ is
an isomorphism of $\genfd$-modules; we want to construct an
algebra $\hat\calD_\phi\supset\calD_\phi$ and an isomorphism 
$\hat\rho: \Shat(Q)\otimes_\genfd S(P)\to\hat\calD_\phi$,
compatible with $\rho$ on $S(Q)\otimes_\genfd S(P)\subset 
\Shat(Q)\otimes_\genfd S(P)$ (recall that $S(P)$ is flat!).
Take $\hat\calD_\phi := \Shat(Q)\otimes_\genfd S(P)$,
$\hat\rho := \id$. For any $p\in P$, the convolution
$D_p : S(Q)\to S(Q)$ maps
$S^n(Q)$ into $S^{n-1}(Q)$,
hence it is continuous and lifts to $\hat{D}_p : \Shat(Q)\to \Shat(Q)$.
Since maps $\hat{L}_p :\hat\calD_\phi\to \hat\calD_\phi$ defined by
$\hat{L}_p :\alpha\otimes\delta \mapsto 
\alpha \otimes p\delta +\hat{D}_p(\alpha)\otimes\delta$, considered for 
different $p\in P$, 
mutually commute, they define some $\hat{L}_\delta$ for all $\delta \in S(P)$.

We {\em define\/} the multiplication $\star$ on $\hat\calD_\phi$
by $(\alpha\otimes \delta)\star(\alpha'\otimes\delta') := ((L_\alpha \otimes 1)
\circ \hat{L}_\delta)(\alpha'\otimes\delta')$, where
$L_\alpha : \Shat(Q)\to \Shat(Q)$, $\beta\mapsto \alpha\beta$ is the usual
multiplication map. It is straightforward to check that 
$\hat\calD_\phi$ is an associative $\genfd$-algebra
and that $\calD_\phi\to \hat\calD_\phi$ is
compatible with multiplication. It is enough to check associativity 
for $u,v,w\in \Shat(Q)\otimes S_{\leq n}(P)\subset\hat\calD_\phi$, 
for arbitrary $n\geq 0$. 
For any $k>0$, one can find some 
$u',v',w'\in S_{< k+3n}(Q) \otimes S_{\leq n}(P)$ such that 
$u\equiv u' \pmod{\Shat_{\geq k + 3n}(Q)\otimes S(P)}$ and so on.
Notice that
$(\Shat_{\geq k}(Q)\otimes S_{\leq n}(P))\star (\Shat_{\geq l}(Q)\otimes
S_{\leq n}(P))\subset \Shat_{\geq k + l - n}(Q)\otimes S_{\leq 2n}(P)$
for any $k,l\geq 0$, hence 
$(u\star v)\star w 
\equiv u'v'w' \equiv u\star(v\star w) \pmod{\Shat_{\geq k}(Q)\otimes
S(P)}$ for all $k>0$. In a similar way, we construct an action of 
$\hat\calD_\phi\cong \Shat(Q)\otimes S(P)$ on $\Shat(Q)$ 
(elements of $P$ act on $\Shat(Q)$ by means of the derivation maps
$\hat{D}_p$ constructed above), and an action of
$\hat\calD_\phi$ on $S(P)$ as well. 

\nxpoint\label{p:classicWeyl} (Classical Weyl algebras) Suppose $P$ is projective of finite 
type, and put $Q :=P^*$ and $\phi : Q\times P\to K$
be {\em minus\/} the canonical pairing. 
We put $\calD_P:=\calD_\phi$,
$\hat\calD_P := \hat\calD_\phi$.
Then $\calD_P = \calD_\phi = SW(P^*\oplus P)$ is a filtered
associative algebra, and 
$\gr (\calD_P) \cong S(Q\oplus P) \cong S(Q)\otimes S(P)$; 
$\calD_P$ acts on $S(Q)$ and $S(P)$, and
$\hat\calD_P$ acts on $\Shat(Q)$ and $S(P)$.

\nxpoint If in addition $P$ is free with base 
$\{e_j\}_{j=1}^n$ and $\{e^k\}_{k=1}^n$
is the dual base of $Q=P^*$, then $\hat\calD_P = SW(Q\oplus P)$ 
is a free associative $\genfd$-algebra in 
$x_k := i(e^k)$ and $\partial^j := i(e_j)$ subject to the relations
$[\partial^k,\partial^l] = 0 = [x_i,x_j]$,
$\left[\partial^k,x_j \right] = \delta^k_{j}$. 
In this way, we obtain the classical 
Weyl algebra written in coordinates. $\hat\calD_P$ 
in this situation corresponds to differential operators of the form
$\sum f_{i_1\ldots i_n}(x_1,\ldots, x_n)
(\partial^1)^{i_1}\cdots(\partial^n)^{i_n}$, where $f_{i_1\ldots i_n}$ are
formal power series, all but finitely many equal to zero. 

\nxpoint\label{p:DerSymm} In the situation of \refpt{p:SWbilin}, the 
$\genfd$-submodule $\calL_P := S(Q)\cdot P\subset \calD_P$
is a Lie subalgebra. Indeed, $\forall\alpha,\alpha' \in S(Q)$
$\forall p,p'\in P$, $(\alpha p)\cdot(\alpha'p') = \alpha\alpha' pp'
+\alpha D_p(\alpha')p'$, hence $[\alpha p,\alpha' p']
=\alpha D_p(\alpha')\cdot p' -\alpha' D_{p'}(\alpha)\cdot p\in \calL_P$.
Recall that there is a $\calD_\phi$-module structure on $S(Q)$, 
for which $S(Q)\subset \calD_\phi$ acts by multiplication
and $P\subset \calD_\phi$ by derivations (namely, convolutions),
hence $\calL_P= S(Q)\cdot P\subset \calD_\phi$
acts on $S(Q)$ by derivations. This way we obtain a Lie algebra
homomorphism $\tau : \calL_P\cong S(Q)\otimes P
\to \Der_\genfd(S(Q))$.

{\bf Proposition.} {\em $\tau$ is an isomorphism under assumptions of 
\refpt{p:classicWeyl}. }

{\bf Proof.} Any derivation $D \in \Der_\genfd(S(Q))$ corresponds to
an algebra homomorphism $\sigma := 1_{S(Q)} + D\epsilon :
S(Q)\to S(Q)[\epsilon]$, $\alpha \mapsto \alpha + D(\alpha)\cdot \epsilon$,
such that $\pi \circ \sigma = 1_{S(Q)}$ for $\pi : S(Q)[\epsilon]
\to S(Q)$, $\epsilon\mapsto0$. 
By the universal property of $S(Q)$, the map $\sigma$ is defined
by its restriction $\sigma|_Q : Q\to S(Q)[\epsilon] =
S(Q)\oplus S(Q)\cdot \epsilon$. Clearly, 
$\sigma|_Q (x) = x + \phi(x)\epsilon$ for some map $\phi : Q\to S(Q)$.
Since $P$ and $Q$ are projective of finite rank, 
$\Hom_\genfd(Q,S(Q))\cong S(Q)\otimes P$, so $\phi$ gives
us an element $\tilde\phi\in S(Q)\otimes P$. One checks 
that $\tau(\tilde\phi)= D$ 
(it is enough to check this on $Q\subset S(Q)$ since a derivation 
of $S(Q)$ is completely determined by its restriction on $Q$).
This way we obtain a map $\Der_\genfd (S(Q))\to 
\calL_P\cong S(Q)\otimes P$ inverse to $\tau$.

\nxpoint\label{p:DerComplSymm} 
Similarly, $\hat\calL_P := \Shat(Q)\cdot P \subset\hat
\calD_P$ is closed under Lie bracket, and it acts by derivations
on $\Shat(Q)$. All {\em continuous} derivations of $\Shat(Q)$ arise
in this way. 

\section{Vector fields on formal affine spaces and end of the proof}
\label{s:EndOfProof}
\addtocounter{subsection}{-1}
\nxpoint Fix a projective $\genrg$-module $P$ of finite type, 
put $Q = P^*$.
We are going to compute the $\genrg$-algebras of vector fields on
$\uW(P)$ and $\uWom(P)$. More
precisely, we will identify these vector fields with derivations of
$S(Q)$ (resp.~$\Shat(Q)$), 
hence with elements of $\calL_P = S(Q)\cdot P
\subset \calD_P$ (resp.\ of $\hat\calL_P\subset 
\hat\calD_P$); we will show that this identification respects Lie 
bracket. Then we are going to use this to compute some vector fields
defined in Section~\ref{s:Formal}.

\nxpoint (Representable functors) Suppose $F \in \calE_{\genrg}
= \functcat(\relPk, \setscat)$ is representable by
some $A = (A,J) \in \Ob\relPk$. 
This means that we have an element $X \in F(A)$, such that 
for any $R = (R,I) \in \relPk$ and any 
$\xi \in F(R)$ there is a unique morphism $\phi : A\to R$ in 
$\relPk$, such that $(F(\phi))(X) = \xi$. One can also
write $F(R) \cong \Hom_{\relPk}(A,R)$ or
$F = \Hom(A,-)$. 

Now consider $TF := \prod_{\genrg[\epsilon]/\genrg} F
: (R,I)\mapsto F(R[\epsilon],I\oplus R\epsilon)$ 
together with the projection
$\pi : TF\to F$ induced by $R[\epsilon]\stackrel{p}\to R$, 
$\epsilon \mapsto 0$. By definition, 
$\Vect(F) = \Gamma(TF/F) = \Hom_F(F,TF)$ is the
set of sections of $TF/F$. By Yoneda lemma, any $\sigma \in \Vect(F)$,
i.e.~a section $\sigma : F\to TF$, is
determined by $\sigma_0 := \sigma_A(X) \in TF(A)
= F(A[\epsilon], J\oplus R\epsilon)
\cong \Hom_{\relPk}\bigl((A,J),(A[\epsilon],J\oplus R\epsilon)\bigr)$.
Denote by $\tilde\sigma_0 : (A,J)\to (A[\epsilon],J\oplus R\epsilon)$ the
corresponding morphism in $\relPk$. Since $\sigma$ is a section of
$\pi$ iff $p_A\circ \tilde{\sigma}_0 = \id_A$, then 
$\tilde\sigma_0 = \id_A + \epsilon\cdot D$ for a uniquely
determined $D : A\to A$, and $\tilde\sigma_0$ is 
$\genrg$-algebra homomorphism iff $D$ is a $\genrg$-derivation of $A$:
$\tilde\sigma_0(ab) = (a + \epsilon D(a))(b+\epsilon D(b)) 
= ab + \epsilon(D(a)b + a D(b))$.
We have constructed a bijection 
$\Vect(F)\stackrel\lambda\to \Der_{\genrg}(A)$. 
One sees immediately that $\lambda$ is an isomorphism of $\genrg$-modules, 
where the $\genrg$-structure
on $\Vect(F)$ comes from the $\genrg$-action 
$[c]_R : R[\epsilon]\to R[\epsilon]$,
$x + y\epsilon \mapsto x + cy\epsilon$ for any $c \in \genrg$.

\nxpoint\label{p:VFDer} {\bf Proposition.} 
{\em $\lambda : \Vect(F)\to \Der_{\genrg}(A)$
is an isomorphism of Lie algebras. }

{\bf Proof.} Recall that the Lie bracket on $\Vect(F) = \Gamma(TF/F)$
is defined as follows. Consider three copies of the dual number algebra,
$\genrg[\epsilon]$, $\genrg[\eta]$, $\genrg[\zeta]$, the tensor product
$\genrg[\epsilon,\eta] = \genrg[\epsilon]\otimes_{\genrg} \genrg[\eta]$
and the embeddings of
$\genrg[\epsilon]$, $\genrg[\eta]$, $\genrg[\zeta]$ into $\genrg[\epsilon,\eta]$,
denoted by $\phi_\epsilon, \phi_\eta,\phi_\zeta$, where
the last map is determined by $\phi_\zeta:\zeta\mapsto \epsilon\eta$.
Given a pair of sections $\sigma,\tau : F\to TF$, we consider 
$\sigma$ as a section of $T_\epsilon F$ and
$\tau$ as a section of $T_\eta F$. Here
$T_\epsilon, T_\eta, T_\zeta$ are the corresponding ``tangent bundles'';
of course $T_\epsilon F\cong  T_\eta F \cong  T_\zeta F \cong TF$; 
besides, $T_\epsilon T_\eta F = T_\eta T_\epsilon F = T_{\epsilon,\eta}F
:= \prod_{\genrg[\epsilon,\eta]/\genrg}(F|_{\genrg[\epsilon,\eta]})$. We have two maps
$F\stackrel\sigma\to T_\epsilon F\stackrel{T_\epsilon(\tau)}
\to T_\epsilon T_\eta F = T_{\epsilon,\eta} F$
and 
$F\stackrel\tau\to T_\eta F\stackrel{T_\eta(\sigma)}
\to T_\eta T_\epsilon F = T_{\epsilon,\eta} F$.
The section $[\sigma,\tau]: F\to T_\zeta F$ is defined by 
$(\phi_\zeta)_*([\sigma,\tau]) = T_\eta(\sigma)\circ\tau -
T_\epsilon (\tau)\circ\sigma$. Now $\forall x\in A$,
$\sigma_A(x) = x+\epsilon\lambda(\sigma)x \in A[\epsilon]$,
$T_\epsilon(\tau)_A\circ\sigma_A: x\mapsto x+\epsilon \lambda(\sigma) x +
\eta \lambda(\tau)x + \epsilon\eta\lambda(\tau)\lambda(\sigma) x$,
and, similarly, 
$T_\eta(\sigma)_A\circ\tau_A : x\mapsto x+\epsilon \lambda(\sigma) x 
+ \eta\lambda(\tau) x+ \epsilon\eta\lambda(\sigma)\lambda(\tau)x$.
Therefore $T_\epsilon(\tau)_A\circ\sigma_A-T_\eta(\sigma)_A\circ\tau_A =
[\lambda(\sigma),\lambda(\tau)]\cdot\epsilon\eta$. Q.E.D.

\nxpoint All this can be applied to $T^\omega F$ instead of $TF$: 
in this case we will have $\sigma_0\in \Hom_{\relPk}\bigl((A,J),
(A[\epsilon], J\oplus J\epsilon)\bigr)$, 
so $\sigma_0 = \id_A + \epsilon D$,
where $D$ is a derivation of $A$ such that
$D(J)\subset J$. This gives an isomorphism
$\lambda^\omega : \Vect^\omega F\to 
\{D \in \Der_{\genrg}(A)\,|\, D(J)\subset J\}$.

\nxpoint\label{p:VFonW} Now consider $F:= \uW(P)$ 
for a projective $\genrg$-module 
$P$ of finite type. Functor $F$ is representable by $A:= (S(Q), 0)$, 
where $Q:=P^*$. Indeed,
$\Hom_{\relPk}(A,(R,I))=\Hom_{\genrg-{\rm alg}}
(S(Q),R)=\Hom_{\genrg}(Q,R) \cong Q^*\otimes_{\genrg}R
\cong R\otimes_{\genrg} P = \uW(P)(R)$.
In this way we see that 
$\Vect(\uW(P))\stackrel{\lambda}\cong \Der_{\genrg}(S(Q))
\cong S(Q)\otimes_{\genrg} P\cong S(Q)\cdot P 
=\calL_P\subset \calD_P$
(cf.~\refpt{p:DerSymm}).
Given a $\sigma\in \Gamma(TF/F)$, one obtains the corresponding element
$\bar{D}$ in $S(Q)\otimes P$
as follows: $\sigma_A : F(A)\to TF(A)= F(A[\epsilon])=
\uW(P)(A[\epsilon]) 
= P\otimes_{\genrg} A[\epsilon]$, 
so $\sigma_A(X) = X+\epsilon \cdot \bar{D}'$
for some $\bar{D}'\in P\otimes_{\genrg}A\cong S(Q)\otimes_{\genrg}P$
(recall $X\in F(A)=P\otimes_{\genrg}A= P\otimes_{\genrg}S(Q)$). 
One can check that $\bar{D} = \bar{D}'$ (this also provides an
alternative proof of~\refpt{p:DerSymm}).

\nxpoint (Pro-representable functors) Consider the full
subcategories $\calP_n\subset \calP$ given by
$\Ob\calP_n = \{(R,I)\in \Ob\calP\,|\,
I^n = 0\}$. These subcategories give an exhaustive filtration 
$\calP = \cup_{n\geq 1} \calP_n$; the corresponding
slice categories are $\relPk_n\subset 
\relPk$.
Fix a functor $F\in \Ob\calE$.
It can happen that $F$ is not representable, but its restrictions
$F^{(n)} := F|_{\relPk_n}$ are.
This means that for each $n\geq 0$ we have some 
$A_n = (A_n,J_n)\in \Ob\relPk_n$
(in particular, $J^n_n = 0$) and an element $X_n\in F(A_n)$ such that
for any $R= (R,I)$ in $\relPk_n$ and any
$\xi\in F(R)$ there is a unique morphism
$\phi : A_n\to R$ such that 
$\xi = (F(\phi))(X_n)$. Since $A_n \in \Ob\relPk_n
\subset \Ob\relPk_{n+1}$, 
the universal property of $A_{n+1}$,
$X_{n+1} \in F(A_{n+1})$ gives us a map $\phi_n : A_{n+1}\to A_n$ such that
$(F(\phi_n))(X_{n+1})= X_n$. In this way we obtain a projective 
system $\underline{A} = (\ldots\to A_3\stackrel{\phi_2}\to
A_2\stackrel{\phi_1}\to A_1)$, or even a pro-object
$\uu{A}:=\quotedprojlim A_n$ over $\relPk$. For any $R=(R,I)$ in 
$\relPk$ we get 
$\Hom_{\Pro\relPk}(\uu{A},R) = 
\injlim_n \Hom_{\relPk}(A_n,R) = F(R)$,  
since if $R\in \Ob\relPk_m$, this inductive system
stabilizes for $n\geq m$.

\nxpoint A vector field $\sigma \in \Vect(F) : F\to TF$ is completely 
determined by its values $\sigma_n := \sigma_{A_n}(X_n)$ lying in
$TF(A_n) = F(A_n[\epsilon],J_n+A_n\epsilon)
\cong\Hom_{\relPk}
\bigl((A_{n+1},$ $J_{n+1}), (A_n[\epsilon],J_n+A_n\epsilon)\bigr)$; let
$\tilde\sigma_n : A_{n+1}\to A_n[\epsilon]$ 
be the corresponding morphism in 
$\relPk$ (here we used that $R\in \Ob\calP_n$ 
implies $R[\epsilon]\in \Ob\calP_{n+1}$ for $R = A_n$).
These $\tilde\sigma_n$ satisfy the obvious compatibility relations
$$\xymatrix{ & A_{n+1}\ar[d]^{\phi_n}\ar[dl]_{\tilde\sigma_n} 
& & A_{n+2} 
\ar[r]^{\tilde\sigma _{n+1}} \ar[d]^{\phi_{n+1}}& 
A_{n+1}[\epsilon]\ar[d]^{\phi_{n}[\epsilon]}
\\
A_n[\epsilon]\ar[r]^{p_{A_n}} & A_n && A_{n+1} 
\ar[r]^{\tilde\sigma _{n}}  & A_{n}[\epsilon]\rlap{\quad.}
}$$
We see that $\tilde\sigma_n = \phi_n + \epsilon D_n$ for some 
$D_n : A_{n+1}\to A_n$ satisfying $\phi_n\circ D_{n+1} 
 = D_n\circ \phi_{n+1}$ and 
$D_n(xy) = D_n(x)\cdot \phi_n(y)+\phi_n(x)\cdot D_n(y)$; 
in this way, one can think of $\uu{D} = (D_n)$
as a derivation of the pro-$\genrg$-algebra $\uu{A} =
\quotedprojlim A_n$. Sections $\Vect^\omega(F)$ of
$T^\omega F$ are treated similarly, but we get additional
conditions $D_n(J_{n+1})\subset J_n$. 
Actually $A_n[\epsilon]^\omega\in\Ob\relPk_n$, so in this case 
we get a compatible family of derivations $D'_n:A_n\to A_n$ such that 
$D_n=D'_n\circ\phi_n$ and $D'_n(J_n)\subset J_n$, and the overall
description of $\Vect^\omega(F)$ is even simpler than that of~$\Vect(F)$.

\nxpoint\label{p:TWomP} Let's apply this to $F = \uWom(P)$. 
First of all, $F|_{\relPk_n}$ is representable by $(A_n, J_n) =
\bigl(S(Q)/S^{\geq n}(Q), S^+(Q)/S^{\geq n}(Q)\bigr)$. Namely, for any 
$(R,I) \in \Ob\relPk_n$
we have $\Hom_{\relPk}((A_n,J_n), (R,I)) = 
\bigl\{\phi \in \Hom_{\genrg-{\rm alg}}(S(Q),R)$ such that $\phi(S^+(Q)^n) = 
0$, $\phi(S^+(Q)) \subset I\bigr\}$, what, since $I^n = 0$, 
equals $\{ \phi\,|\,\phi(S^+(Q))\subset I\} \cong
\{ \tilde\phi\in \Hom_{\genrg}(Q,R) \,|\,
\tilde\phi(Q)\subset I\} = \Hom_{\genrg}(Q,I) \cong I\otimes_{\genrg} P$.

We see that the vector fields $\sigma \in \Vect \uWom(P)$
correspond to compatible families 
$\uu{D} = (D_n)$ of ``derivations''
$D_n : S(Q)/S^{\geq n+1}(Q)\to S(Q)/S^{\geq n}(Q)$. One can take the ``true''
projective limit and obtain a continuous derivation 
$D : \Shat(Q)\to\Shat(Q)$, that corresponds by~\refpt{p:DerComplSymm} to
some element of $\hat\calL_P \cong \Shat(Q)\otimes P
\cong \Shat(Q)\cdot P \subset \hat\calD_P$. This is a Lie algebra
isomorphism by the same reasoning as in~\refpt{p:VFDer}. 
Again, given a $\sigma \in \Vect \uWom(P)$, 
we can construct the corresponding element 
$\bar{D} \in \calL_P\cong \Shat(Q)\otimes P$ as follows:
apply $\sigma_{A_n} : \uWom(P)(A_n) = J_n\otimes_{\genrg} P 
\to T\uWom(P)(A_n) = \uWom(P)(A_n[\epsilon])
= J_n\otimes_{\genrg}P \oplus A_n\epsilon\otimes_{\genrg} P$ to
$X_n \in J_n\otimes_{\genrg} P = \uWom(P)(A_n)$ 
and get some element $\sigma_{A_n}(X_n) = X_n + \epsilon\cdot \bar{D}_n$,
$\bar{D}_n \in A_n\otimes_{\genrg} P
= S(Q)/S^{\geq n}(Q) \otimes_{\genrg} P$.
These $\bar{D}_n$ form a compatible family that defines an element of
$\Shat(Q)\otimes_{\genrg} P$; this element is exactly $\bar{D}$ (proof 
is similar to~\refpt{p:VFonW}). 

\nxpoint\label{p:derofAct} 
Suppose $G$ is a group in $\calE$, $F$ is an object
of $\calE$ and $G$ acts on $F$ from the left: we are given some
$\alpha : G\times F\to F$ satisfying usual properties. Since $T$ is left exact,
$T(G\times F) = TG\times TF$, and we get a left action $T\alpha :
TG \times TF\to TF$ compatible with $\alpha$:
$$\xymatrix{
TG\times TF \ar[rr]^{T\alpha}\ar[d]^{\pi_G\times \pi_F} && TF\ar[d]^{\pi_F}\\
G\times F\ar[rr]^\alpha && F
}$$
Hence $\uLie(G)\subset TG$ also acts on $TF$, so we get a map
$\uLie(G)\times TF\stackrel\beta\to TF$. Since $\pi_G$ maps
$\uLie(G)$ into the identity of $G$, the following diagram is commutative:
$$\xymatrix{ \uLie(G)\times TF\ar[d]^{\pr_2}
\ar[rr]^<>(.5){\beta} &&
TF\ar[d]^{\pi_F} \\
TF \ar[rr]^{\pi_F} && F
}$$
After composing $\beta$ with $\id_{\uLie(G)}\times s$ where
$s : F\to TF$ is the zero section of $TF$, we get a map
$\gamma : \uLie(G)\times F\to TF$ over $F$, hence a map
$\gamma^\flat : \uLie(G)\to \iHom_F(F,TF) =\uVect(F)$,
and by taking the global sections (=evaluating at $\genrg$) we obtain
an $\genrg$-linear map
$\Gamma(\gamma^\flat) : \Gamma(\uLie(G))=\Lie(G)
\to \Gamma(\uVect(F))=\Vect(F)$. 
One checks, by means of the description of Lie bracket on~$\uVect(F)$
given in~\refpt{p:VFDer} and a similar description 
of the Lie bracket on~$\uLie(G)$ recalled in~\refpt{p:LieBr},
that $\gamma^\flat$ and $d_e\alpha := \Gamma(\gamma^\flat)$ 
are Lie algebra homomorphisms. 

\nxpoint\label{p:RightInv} Fix a Lie algebra $\gg$ over $\genrg$,
finitely generated projective as an $\genrg$-module, and construct
the formal group $G := \Expt(\gg)$ with
$\uLie(G) = \uW(\gg)$,
hence $\Lie(G) = \Gamma(\uW(\gg)) = \gg$,
as in~\refpt{p:vectfgrp}. Consider first the left action of $G$ on itself
given by the multiplication map $\mu : G\times G\to G$.
According to~\refpt{p:derofAct},
we get a map $d_e\mu : \gg\to \Vect(G) = \Gamma(TG/G)$. It is
clear from the description given in~\refpt{p:derofAct}, that
$Y \in \gg$ maps to the {\em right-invariant vector field\/}
$\sigma_Y : G\to TG$ given by
$(\sigma_Y)_R : g \mapsto (1+Y\epsilon)\cdot s_R(g)$
[here $1+Y\epsilon$ denotes the image of $Y\in \uLie(G)(R)$
in $TG(R)$, $s : T\to TG$ is the zero section of $TG$
and $g$ is an element of $G(R)$, $R\in \Ob\relPk$].
This means that, if we identify $TG$ with
$G\times \uLie(G)\cong G\times\uW(\gg)$ by means of
the map $(g,Y)\mapsto (1+Y\epsilon)\cdot s(g)$, then
$\sigma_Y\in \Hom_G(G,TG)\cong \Hom_G(G,G\times \uLie(G))\cong
\Hom(G,\uLie(G))$ is identified with the constant map
$g\mapsto Y$ in $\Hom(G,\uLie(G))=\Hom(G,\uW(\gg))$.

\nxpoint Recall that the exponential map 
$\exp' : \uWom(\gg)\to G$ gives an isomorphism of
formal schemes, hence we can deduce from $\mu$ a left action 
$\alpha : G\times \uWom (\gg)
\to \uWom(\gg)$:
$$\xymatrix{ G\times \uWom(\gg)
\ar[d]^{\id_G\times \exp'}
\ar[r]^<>(.5){\alpha}
& \uWom(\gg)\ar[d]^{\exp'}_\sim
\\ G\times G\ar[r]^<>(.5){\mu} & G
}$$
This would give us a Lie algebra homomorphism $d_e\alpha : \gg\to
\Vect(\uWom(\gg)) = 
\Hom_{\uWom(\gg)}(\uWom(\gg), T\uWom(\gg))
\cong \Hom(\uWom(\gg), \uW(\gg))$. We want to compute explicitly
the vector fields $\tilde\sigma_Y = (d_e\alpha)(Y)$ in terms
of this isomorphism. 

We have the following diagram (cf.~\refpt{p:compTexp}):
$$\xymatrix{
\uWom(\gg)\ar[r]^{\tilde\sigma_Y}\ar[d]^\sim_{\exp'}&
T(\uWom(\gg))\ar[r]^<>(.5){\sim}
\ar[d]^{T(\exp')}_\sim &
\uWom(\gg)\times\uW(\gg)
\ar[rr] 
\ar@{-->}[d]^\tau&&
\uWom(\gg)\ar[d]^{\exp'}
\\
G\ar[r]^{\sigma_Y} &
TG \ar[r]^<>(.5){\sim} &
G\times \uW(\gg)
\ar[rr]^<>(.5){\pi_G}
&&G\ar@/_1pc/@{.>}[ll]_s
}$$
Here $\tau$ is given by $(X,Z)\mapsto \bigl(\exp(X), 
\sum_{n\geq 0}\frac{(\ad X)^{n+1}}{(n+1)!}(Z)\bigr)$ for any 
$X$ lying in $\uWom(\gg)(R) = I_{R}\cdot \gg_{(R)}$ and any
$Z$ from $\uW(\gg)(R) = \gg_{(R)}$ [ $I_{R}$ is
nilpotent, so is $(\ad X)$, hence the sum is finite; note the
absence of the factor $(-1)^n$ factor in comparison to~\refpt{p:compTexp}; 
this is due to the fact that we have chosen here another splitting 
$TG \cong G\times\uLie(G)$, given by {\em right-invariant\/} vector fields].
On the other hand, by~\refpt{p:RightInv}, $\sigma_Y$ is given by $g\mapsto Y$.
This means that $\tilde\sigma_Y$ maps $X$ into $(X,Z)$ such that
$\sum_{n\geq 0}\frac{(\ad X)^{n+1}}{(n+1)!}(Z) = Y$, i.e.\ 
$P(\ad X)(Z) = Y$, where
$P(T) \in \bbQ[[T]]$ is the series $P(T) = (e^T-1)/T$. 
Therefore, $Z = P(\ad X)^{-1}(Y)$, and classically 
$P(T)^{-1} = T/(e^T-1) = \sum_{n\geq 0} \frac{B_n}{n!} T^n$ 
(this is actually the definition of Bernoulli numbers $B_n$). 

\nxpoint\label{p:DefTheta} (Definition of embedding $\theta$) 
We have just seen that for $Y \in \gg$ the vector field $\tilde\sigma_Y$ 
is given by 
$X\mapsto \bigl(X,\sum_{n\geq 0}\frac{B_n}{n!} (\ad X)^n(Y)\bigr)$. 
On the other hand, by~\refpt{p:TWomP}, 
we know that vector fields on $\uW(\gg)$ 
correspond to continuous derivations of $\Shat(\gg^*)$, or to
the elements of $\Shat(\gg^*)\otimes \gg\cong 
\Shat(\gg^*)\cdot \gg = \hat\calL_\gg\subset 
\hat\calD_\gg$, where $\hat\calD_\gg$ is
the completed Weyl algebra of $\gg$, cf.~\refpt{p:ComplWeyl}, 
\refpt{p:DerComplSymm}.
We want to compute the elements of $D_Y \in \hat\calL_\gg$ that
correspond to $\tilde\sigma_Y\in \Vect(\uWom(\gg))$; 
this would give us a Lie algebra homomorphism 
$\gg\stackrel\theta\to\hat\calL_\gg\subset
\hat\calD_\gg$, $Y\mapsto D_Y$, hence also a homomorphism
$\Uenv(\gg)\stackrel{\tilde{\theta}}
\to\hat\calD_\gg$. We will see in~\refpt{p:injofTheta} that both
$\theta$ and $\tilde{\theta}$ are injective.

\nxpoint Let's apply~\refpt{p:TWomP} for $P = \gg$, 
$Q = \gg^*$, $F = \uWom(\gg)$, 
to compute the element $D_Y \in \Shat(\gg^*)\otimes\gg$ 
corresponding to $\tilde\sigma_Y$ defined by some $Y\in \gg$.
We know that $F$ is pro-representable by 
$\uu{A} = \quotedprojlim (A_n,\phi_n)$, where
$A_n = (A_n,J_n) = (S(\gg^*)/S^{\geq n}(\gg^*), S^+(\gg^*)/
S^{\geq n}(\gg^*))$, $\phi_n : A_{n+1}\to A_n$ is the projection.
We have also the universal elements $X_n \in A_n$; such an 
element is equal to the image of the canonical element 
$c_\gg \in \gg^*\otimes\gg$ in 
$J_n \otimes_{\genrg} \gg =\uWom(\gg)(A_n)
\subset \gg_{(A_n)}$. According to~\refpt{p:TWomP},
we have to apply $(\tilde{\sigma}_Y)_{A_n} : \uWom(\gg)(A_n)
\to T\uWom(\gg)(A_n) = J_n\otimes_{\genrg} \gg\oplus 
\epsilon A_n\otimes_{\genrg}\gg$ to $X_n$ and to take the second
component $D_{Y,n}$. According to~\refpt{p:DefTheta},
$(\tilde{\sigma}_Y)_{A_n}(X_n) = X_n + \epsilon \cdot
\sum_{k = 0}^n \frac{B_k}{k!} (\ad X_n)^k(Y)$; this gives us 
the value of $D_{Y,n}$.
Elements $X_n \in A_n$ define an universal element
$X \in \projlim(J_n\otimes_{\genrg}\gg)\cong \Shat^+(\gg^*)
\otimes_{\genrg}\gg\subset \gg_{(\Ahat)}$,
where $(\Ahat,\hat{J}) = (\Shat(\gg^*),\Shat^+(\gg^*))$
is the ``true'' (topological) projective limit of $\uu{A}$.
Of course, $X$ is still the image of $c_\gg \in \gg^*\otimes
\gg$ in $\hat{J}\otimes\gg\subset \gg_{(\Ahat)}$.
We see that $D_Y = \sum_{k\geq 0} \frac{B_k}{k!}(\ad X)^k(Y)\in 
\gg_{(\Ahat)}\cong \Shat(\gg^*)
\otimes_{\genrg}\gg$. This power series converges
since $X \in \hat{J}\cdot \gg_{(\Ahat)}$ is topologically nilpotent. 
Here $Y\in \gg$ is considered as an element of 
$\gg_{(\Ahat)}\supset \gg$,
and $(\ad X)^k(Y)$ is computed with respect to the 
$\Ahat$-Lie algebra structure on $\gg_{(\Ahat)}$.

\nxpoint\label{p:MainFormula} (Main formula) 

Suppose $\gg$ is a free $\genrg$-module of rank $n$.
Fix a base $(e_i)_{1\leq i\leq n}$ of $\gg$, and consider the
structural constants $C^k_{ij}\in \genrg$ defined by
$[e_i,e_j] = \sum_k C^k_{ij} e_k$. (Here the completed
Weyl algebra is used, hence, unlike in Sections 1--6,
there is no need to use a formal variable $t$,
hence to distinguish $C^i_{jk}$ and $(C^0)^i_{jk}$.)
Denote by $(e^i)$ the dual base
of $\gg^*$, denote by $\partial^i$ the images of
$e_i$ in $S(\gg)\subset \Shat(\gg^*)\otimes_{\genrg} S(\gg)
\cong \hat\calD_\gg$ and by $x_i$ -- the images of
$e^i$ in $\Shat(\gg^*) \subset \hat\calD_\gg$.
(The apparent loss of covariance/contravariance here is
due to the fact we'll need to apply the Weyl algebra automorphism
$x_i\mapsto-\partial^i$, $\partial^i\mapsto x_i$ to recover the
Main Formula in form~\eqref{eq:phi}; there doesn't seem to be
a completely satisfactory way of fixing this.)

Clearly $c_\gg = \sum_i e^i \otimes e_i$, hence
$X = \sum_i x_i e_i \in \gg_{(\Ahat)}=
\gg_{(\Shat(\gg^*))}$ is the universal
element, and for any $Y\in \gg$,
$D_Y$ is given by $\sum_{s\geq 0} \frac {B_s}{s!} (\ad X)^s(Y)$.
In coordinates, $\ad X \in \End_\Ahat(\gg_{(\Ahat)})$
is given by $\ad X : e_j \mapsto \sum_i x_i [e_i, e_j]
=\sum_{i,k} C_{ij}^k x_i e_k$, hence
the matrix $M = (M^i_j)$ of $\ad X$ is given by
$M^i_j = \sum_{k} C_{kj}^i x_k$ (cf.\ with $\CE^i_j$ from
Sections 1--5, which involve $-\partial^k$-s in the place of $x_k$-s).

For $Y = e_j$ we obtain
\begin{equation}\label{eq:mainfDversion}
\begin{split}
D_{e_j} = D_Y = &\textstyle \sum_{s\geq 0} \frac{B_s}{s!} (\ad X)^s (e_j)
= \sum_{s\geq 0} \sum_{i=1}^n \frac{B_s}{s!} (M^s)_j^i e_i \\
= &\textstyle \sum_{i = 1}^n
   \bigl(\sum_{s = 0}^\infty \frac{B_s}{s!} (M^s)_j^i \bigr)\partial^i
\in \calL_\gg\subset\hat\calD_\gg
\end{split}
\end{equation}
Thus we have constructed an explicit embedding
$e_j\mapsto D_{e_j}$ of $\gg$ into the completed
Weyl algebra $\hat\calD_\gg$. Recall that
$\hat\calD_\gg$ is some completion of
the Weyl algebra $\calD_\gg$, and that $\calD_\gg$
in this situation is the free algebra over $\genrg$ generated by
$x_1,\ldots, x_n,\partial^1,\ldots, \partial^n$ subject to the relations
$[x_i,x_j] = [ \partial^i,\partial^j ] = 0$,
$[\partial^k, x_i] = \delta^k_i$,
i.e.\ is the classical Weyl algebra over $\genrg$
with $2n=2\dim\gg$ generators.

\nxpoint\label{p:injofTheta} (Injectivity of $\theta$ and $\tilde{\theta}$)
Now it remains to show that our homomorphisms $\theta:\gg\to\hat\calL_\gg$ and 
$\tilde\theta:\Uenv(\gg)\to\hat\calD_\gg$ are injective. To achieve this we 
consider the ``evaluation at origin map'' 
$\beta: \Shat(\gg^*)\to \Shat(\gg^*)/\Shat^+(\gg^*)=\genrg$ 
and the induced maps 
$\beta\otimes 1_\gg : \hat\calL_\gg\to\gg$ and 
$\beta\otimes 1_{S(\gg)} : \hat\calD_\gg\to S(\gg)$. 
One checks immediately that
$(\beta\otimes 1_\gg)\circ\theta = 1_\gg$, hence $\theta$ is injective; 
for $\tilde\theta$ observe that $(\beta\otimes 1_{S(\gg)})\circ\tilde\theta$ 
maps $\Uenv_n(\gg)$ into $S_n(\gg)$ and induces the identity map between 
the associated graded $\gr(\Uenv(\gg))\cong S(\gg)$ and $\gr(S(\gg))=S(\gg)$, 
hence is injective, hence $\tilde\theta$ is also injective.

\section{Another proof in the language of coderivations}
\label{sec:thirdproof}

\nxpoint In this section, $\genrg\supset \bbQ$. 
If $H$ is any $\genfd$-coalgebra with comultiplication $\Delta_H$, 
then a $\genfd$-linear map $D \in \End_\genfd(H)$ is a {\bf coderivation}
if $(D\otimes 1 + 1 \otimes D)\circ \Delta_H = \Delta_H \circ D$.
We denote by $\Coder_\genfd(H)$ the set 
of all $\genfd$-linear coderivations of~$H$.

\nxpoint Let $\gg$ be a Lie algebra over $\genfd$ finitely generated 
as an $\genfd$-module. Put $H := \Uenv(\gg)$ and
$C := S(\gg)$. Both $C$ and $H$ are Hopf algebras, cocommutative 
as coalgebras. Moreover, there is a unique canonical (functorial in $\gg$)
isomorphism of {\it coalgebras} $\xi :C\stackrel\sim\to H$ 
that is identity on $\gg$ ([Bourbaki], Ch.~II). This 
{\it coexponential map} may be described as follows: it maps
$x_1\ldots x_n \in C^n = S^n(\gg)$ into
$\frac{1}{n!} \sum_{\sigma \in S_n} x_{\sigma(1)} \ldots x_{\sigma(n)}
\in \Uenv(\gg)_n \subset H$.

\nxpoint Consider the left action of $H$ on itself defined 
by its algebra structure: $h \mapsto L_h \in \End_\genfd (H)$,
$L_h : h'\mapsto hh'$. Since $H = \Uenv(\gg)$, this restricts to 
a Lie algebra action of $\gg$ on $H$. 
{\it This is an action by coderivations.} Indeed, for any 
$h \in \gg$, we have $\Delta_H(h) = 1\otimes h + h\otimes 1$,
hence $(L_h\otimes 1 + 1\otimes L_h) (\Delta_H(h')) =
\Delta_H(h)\cdot \Delta_H(h') = \Delta_H(hh') = (\Delta_H\circ L_h) (h').$ 

\nxpoint The coexponential isomorphism $\xi : C\stackrel\sim\to H$ allows us 
to define a left action $h \mapsto D_h := \xi^{-1} \circ L_h \circ \xi$ of 
$H$ on $C$, such that $\gg \subset H$ acts on $C$ by {\it coderivations},
and $D_h (1) = h$ for any $h \in \gg$. Thus we get $\theta : h\mapsto D_h$,
$\theta : \gg \to \Coder_\genfd(C) = \Coder_\genfd(S(\gg))$.
By ``duality'' $S(\gg)^* \cong \hat{S}(\gg^*)$ (we use here that
$\gg$ is projective f.g., $\genfd \supset \bbQ$) since
$S(\gg)^* = (\oplus S^n(\gg))^* = \prod S^n(\gg)^* \cong \prod S^n(\gg^*) 
= \hat{S}(\gg^*)$, and $D_h$ give us derivations
${}^t D_h : \hat{S}(\gg^*) \to \hat{S}(\gg^*)$. Hence we get a Lie 
algebra homomorphism $\theta : \gg\to\Der(\hat{S}(\gg^*))\cong \hat{S}(\gg^*)
\otimes \gg\cong \hat\calL_\gg \subset\hat\calD_\gg$,
where $\hat\calD_\gg$ is the completed Weyl algebra of $\gg$, 
$\theta : h \mapsto {}^t D_h$.

\nxpoint The pairing between $S^n(\gg^*)$ and $S^n(\gg)$ that induces the 
isomorphism $S^n(\gg)^* \cong S^n(\gg^*)$ used above, is given by
$\langle x_1\ldots x_n, u_1\ldots u_n\rangle =
\sum_{\sigma \in S_n} \frac{1}{n!} \langle x_i, u_{\sigma(i)}\rangle$, i.e.
$\langle x^n, u^n\rangle = \langle x, u\rangle^n$. Elements of the form
$x^n$ generate $S^n(\gg)$ and $\xi(x^n) = x^n \in \Uenv(\gg)$ 
$\forall x \in \gg$. We see that  $D_h(x^n) = \xi^{-1}(hx^n)$. Moreover, 
$\xi(x^k h) = \frac{1}{k+1}\sum_{p+q = k} x^p h x^q$ and $(\ad x).\xi(x^k h)
= \frac{1}{k+1} (x^{k+1} h - h x^{k+1}).$ For the following computations,
fix $x \in H$ and denote by $L = L_x$ and $R = R_x$ the left and 
right multiplication by $x$. Since $L$ and $R$ commute 
and $\genfd\supset \bbQ$, the polynomials from $\bbQ[L,R]$ act on $H$.
In particular, 
$$\begin{array}{l}
\xi(x^k h) = \frac{1}{k+1}\sum_{p+q = k} x^p h x^q
= (\frac{1}{k+1}\sum_{p+q = k} L^p R^q) (h),\\
(\ad x).h = xh - hx = (L-R)(h),\\
\xi(x^k (\ad x)^l h) = \frac{1}{k+1}(\sum_{p+q = k} L^p R^q)(L-R)^l (h).
\end{array}$$
We would like to find rational coefficients $a^{(n)}_{k,l} \in \bbQ$,
such that $D_h(x^n) = \sum_{k,l} a_{k,l}^{(n)} x^k (\ad x)^l.h$. This
condition can be rewritten as 
$$
hx^n = \sum_{k,l} a_{k,l}^{(n)} \xi (x^k (\ad x)^l h)\,\,\,\,\, \mbox{ i.e. }
\,\,\,\,\,
R^n(h) = \sum_{k,l} \frac{a^{(n)}_{k,l}}{k+1} \left(\sum_{p+q = k} L^p R^q
\right) (L-R)^l (h).
$$
For this it is enough to require
$$
R^n = \sum_{k,l} \frac{a^{(n)}_{k,l}}{k+1} \left(\sum_{p+q = k} L^p R^q\right) 
(L-R)^l \mbox{ to hold in } \bbQ[L,R].
$$
Since $L-R$ is not a zero divisor in $\bbQ[L,R]$, this identity is 
equivalent to
$$
(L-R) R^n = \sum_{k,l} \frac{{a}^{(n)}_{k,l}}{k+1} (L^{k+1} - R^{k+1}) (L-R)^l
$$
Now the LHS is a homogeneous polynomial of degree $n+1$, hence 
in finding $a^{(n)}_{k,l}$ we may also require that all summands of 
other degrees on the RHS vanish as well, i.e. assume that
$a^{(n)}_{k,l} = 0$ unless $k+l = n$. 
This also allows to simplify the notation,
namely set $a^{(n)}_k := a^{(n)}_{k,n-k}$, $1\leq k\leq n$.
Next, consider the isomorphism $\bbQ[L,R] \cong \bbQ[X,Y]$ given by
$R\mapsto X$, $L\mapsto X+Y$; it allows us to rewrite our identity as
$$
X^n Y = \sum_{k=0}^n \frac{a^{(n)}_k}{k+1} \left(
(X+Y)^{k+1} - X^{k+1}\right) Y^{n-k}
$$
Divide both sides by $Y^{n+1}$ and put $T := X/Y$; so we get the
following identity in $\bbQ[T]$:
\begin{equation}\label{eq:teq}
T^n = \sum_{k=0}^n \frac{a^{(n)}_k}{k+1} \bigl((T+1)^{k+1} - T^{k+1}\bigr)
\end{equation}
For any $P \in \bbQ[T]$ set $\delta P := P(T+1) - P(T)$, $DP := P'(T)$. 
(We don't use the classical notation $\Delta P$ for $P(T+1)-P(T)$ 
since it might be confused with our notation for the comultiplication.)
By Taylor's formula 
$(e^D P)(T) = \sum_{k\geq 0} \frac{P^{(k)}(T)}{k!} = P(T+1)$, hence
$\delta P = (e^D -1)P$, i.e. $\delta = e^D - 1$. 
In terms of $P_n(T) : = \sum_{k=0}^n \frac{a^{(n)}_k}{k+1} T^{k+1}$,
the equation~(\ref{eq:teq}) may be rewritten as
\begin{equation}
\delta P_n(T) = T^n,\,\,\,\,\,\,\,\,P_n(0) = 0.
\end{equation}
This determines $P_n(T)$ uniquely, and the coefficients of $P_n(T)$ can be
expressed in terms of Bernoulli numbers. To obtain this expression,
observe that $\delta = \frac{e^D - 1}{D} D$, 
so $\left(\frac{e^D-1}{D}\right)(DP_n) = T^n$, therefore $DP_n = 
\left(\frac{D}{e^D-1}\right)(T^n)$. By the definition of Bernoulli
numbers $B_k$, we have $\frac{D}{e^D-1} = \sum_{k\geq 0}\frac{B_k}{k!} D^k$,
hence 
$$
DP_n(T) = \sum_{k\geq 0} \frac{B_k}{k!} D^k T^n = 
\sum_{k\geq 0} B_k {n\choose k} T^{n-k} = \sum_{k=0}^n B_{n-k} 
{n\choose k} T^k.
$$
On the other hand, $DP_n(T) = \sum_{k = 0}^n a^{(n)}_k T^k$, hence
$a^{(n)}_k = {n\choose k} B_{n-k}$.
We have proved the following formula:

\nxpoint \label{p:adxpairing}
For any $h \in \gg \subset H$, $x \in \gg \subset C$, we have
\begin{equation}\label{eq:dhxn}
D_h(x^n) = \xi^{-1} (hx^n) = \sum_{k = 0}^n a^{(n)}_k x^k 
\cdot((\ad x)^{n-k}(h)) = \sum_{k =0}^n {n\choose k} B_{n-k} x^k 
\cdot ((\ad x)^{n-k}(h)), 
\end{equation}
or, shortly,
\begin{equation}
D_h(x^n) = \sum_{k = 0}^n {n\choose k} B_k x^{n-k}
\cdot ((\ad x)^k (h)). 
\end{equation}
This implies for any $x \in \gg$, $\alpha \in \gg^*$, $n\geq 0$, 
$$
\langle x^n, {}^t D_h(\alpha)\rangle = \langle D_h(x^n), \alpha\rangle
= \langle B_n \cdot (\ad x)^n (h), \alpha \rangle
= B_n \langle (\ad x)^n (h),\alpha\rangle.
$$
Both sides are polynomial maps in $x$ of degree $n$, by taking 
polarizations we get
$$
\langle x_1 \ldots x_n, {}^t D_h (\alpha) \rangle = 
\frac{B_n}{n!} \sum_{\sigma \in S_n} \left\langle (\ad x_{\sigma(1)})\ldots
(\ad x_{\sigma(n)})(h),\alpha\right\rangle,\,\,\,\,\,
\forall x_1,\ldots, x_n \in \gg.
$$
The above formula is by no means new: it can be found for example
in~\cite{Petracci}, Remark~3.4.

\nxpoint The coderivation $D_h : S(\gg) \to S(\gg)$ corresponds by duality to 
a continuous derivation ${}^t D_h : \hat{S}(\gg^*) \to \hat{S}(\gg^*)$, which
is uniquely determined by its restriction 
${}^t D_h |_{\gg^*} : \gg^* \to \hat{S}(\gg^*)$.
Such a map corresponds to an element 
${}^t \tilde{D}_h \in \hat{S}(\gg^*)\otimes \gg \cong \hat\calL_\gg
\subset \hat{S}(\gg^*) \otimes S(\gg) \cong \hat\calD_\gg.$
Now we want to compute this element $\tilde{D}_h := ({}^t D_h)^\sim
\in \hat\calL_\gg \subset \hat\calD_\gg$ of the Weyl algebra 
since it defines (by means of the usual action of $\hat\calD_\gg$ on
$\hat{S}(\gg^*)$) the derivation ${}^t D_h$ as well as the coderivation $D_h$.

\nxpoint Any element $u \in \gg_{(\hat{S}(\gg^*))} = \hat{S}(\gg^*)\otimes 
\gg = \hat\calL_\gg$ corresponds to a $\genfd$-linear map
$u^\sharp : S(\gg)\to\gg$ and conversely,
since $\gg$ is a projective $\genfd$-module of finite type. 
Note that $\gg_{(\hat{S}(\gg^*))}$ has a $\hat{S}(\gg^*)$-Lie algebra structure
obtained by base change from that of $\gg$. Let's compute $[u,v]^\sharp$
in terms of $u^\sharp$ and $v^\sharp : S(\gg)\to \gg$ for any two elements
$u,v \in \gg_{(\hat{S}(\gg^*))}.$

\nxsubpoint To do this consider the following more general situation:
Let $M,N,P$ be projective $\genfd$-modules of finite rank, 
$B : M\times N\to P$ a $\genfd$-bilinear map,
$\tilde B : M\otimes N\to P$ the map induced by $B$,
$u \in M_{(\hat{S}(\gg^*))}$, $v \in N_{(\hat{S}(\gg^*))}$.
Suppose $u^\sharp : S(\gg)\to M$ and $v^\sharp : S(\gg)\to N$ 
correspond to $u$ and $v$, and we want to compute the map
$w^\sharp : S(\gg)\to P$ corresponding to
$w = B_{(\hat{S}(\gg^*))}(u,v)$ in terms of $u^\sharp$ and $v^\sharp$.

\nxsubpoint The answer here is the following: $w^\sharp = \tilde{B}\circ 
(u^\sharp \otimes v^\sharp)\circ \Delta$ where 
$\Delta : S(\gg)\to S(\gg)\otimes S(\gg)$ is the comultiplication on $S(\gg)$.
Indeed, by linearity, it is sufficient to check this for
$u = x\otimes\phi$, $v = y\otimes\psi$, $x\in M$, $y\in N$,
$\phi,\psi \in \hat{S}(\gg^*) = S(\gg)^*$. Then 
$u^\sharp : \lambda \mapsto \phi(\lambda)x$,
$v^\sharp : \lambda \mapsto \psi(\lambda)y$,
$w= z\otimes \phi\psi$ where $z = B(x,y)\in P$ and
$\phi\psi = (\phi\otimes\psi)\circ\Delta$ by duality
between algebra $\hat{S}(\gg^*)$ and coalgebra $S(\gg)$. Hence 
$w^\sharp : \lambda \mapsto (\phi\psi)(\lambda) \cdot z = (\phi\otimes\psi)
(\Delta(\lambda))\cdot B(x,y) = (\tilde{B}\circ (u^\sharp\otimes v^\sharp)
\circ\Delta)(\lambda)$. Q.E.D.

\nxsubpoint \label{p:unpgequal}
Now we apply this for $M = N = P = \gg$, $B$ -- the multiplication
map. We see that for any two $u^\sharp, v^\sharp : S(\gg)\to \gg$ we have 
$[u,v]^\sharp = \mu_\gg \circ (u^\sharp\otimes v^\sharp)\circ \Delta$
where $\Delta = \Delta_{S(\gg)} : S(\gg)\to S(\gg)\otimes S(\gg)$ is
the comultiplication of $S(\gg)$ and $\mu_\gg : \gg\otimes\gg \to \gg$,
$x\otimes y\mapsto [x,y]$ is the bracket multiplication of Lie algebra
$\gg$. We introduce a Lie bracket $[,]_\gg$ on 
$\Hom(S(\gg),\gg)$ by this rule; 
then $[u,v]^\sharp = [u^\sharp, v^\sharp]_\gg$,
so $\gg_{(\hat{S}(\gg^*))}\to \Hom(S(\gg),\gg)$, $u\mapsto u^\sharp$
is an isomorphism of Lie algebras. 
Observe that there is another Lie bracket $[,]_\calD$ on 
$\Hom(S(\gg),\gg)\cong \Dercont(\hat{S}(\gg))\cong
\Coder(S(\gg))$ given by the usual commutator of coderivations:
$[D_1,D_2]_{\calD} = D_1D_2 - D_2 D_1$. (One must be careful when passing from
derivations to coderivations since $[{}^t D_1,{}^t D_2] = -{}^t[D_1,D_2]$).

\nxpoint We have some specific elements in $\gg_{(\hat{S}(\gg^*))}$, hence in
$\Hom(S(\gg),\gg)$. First of all, any $a \in \gg$ lies in 
$\gg_{(\hat{S}(\gg^*))}$, inducing therefore a map $a^\sharp : S(\gg)\to\gg$.
Clearly, this is the map sending $1 \in S(\gg)$ to $a$, 
and sending all of $S^+(\gg) = \oplus_{n\geq 1} S^n(\gg)$ to zero. 

\nxsubpoint One checks immediately that 
$[a^\sharp, b^\sharp]_\gg = [a,b]^\sharp$, $\forall a,b \in \gg$,
$[a^\sharp, u^\sharp]_\gg = (\ad a)\circ u^\sharp$, $[u^\sharp, a^\sharp]_\gg
= -(\ad a)\circ u^\sharp$, $\forall a\in \gg$, $u^\sharp : S(\gg)\to\gg$.

\nxsubpoint Besides, the ``canonical element'' $X \in \gg\otimes\gg^*
\subset \gg_{(\hat{S}(\gg^*))}$ 
(the image of $\id_\gg$ under the identification 
$\Hom(\gg,\gg)\cong \gg\otimes\gg^*$) provides a map
$X^\sharp : S(\gg)\to \gg$, which is clearly the projection of
$S(\gg) = \oplus_{n\geq 0} S^n(\gg)$ onto $S^1(\gg) = \gg$.
In particular, for all $a \in \gg$, $X^\sharp (a) = a$ if $n= 1$,
and $X^\sharp (a^n) = 0$ otherwise. 

\nxsubpoint Let's compute $[X^\sharp, u^\sharp]_\gg$ for any 
$u \in \Hom(S(\gg),\gg)$. Since elements of the form $a^n$ generate
$S^n(\gg)$, it is sufficient to determine all
$[X^\sharp, u^\sharp](a^n)$. Now, $\Delta(a^n) = \Delta(a)^n
= (a\otimes 1 + 1 \otimes a)^n = 
\sum_{k=0}^n {n \choose k} a^k\otimes a^{n-k}$, hence 
$(X^\sharp \otimes u^\sharp)(\Delta(a^n)) = na \otimes u^\sharp(a^{n-1})$,
so $[X^\sharp, u^\sharp]_\gg(a^n) = n.[a,u^\sharp(a^{n-1})]$ (for
$n=0$ the RHS is assumed to be zero by convention). In particular, 
if $u^\sharp$ was zero restricted to $S^n(\gg)$ for fixed $n$, then 
$[X^\sharp, u^\sharp]_\gg$ is zero on $S^{n+1}(\gg)$.

\nxsubpoint \label{p:sharpconc}  Thus, we have proved 
$(\ad_\gg X^\sharp) (u^\sharp) : a^n \mapsto n(\ad a)(u^\sharp(a^{n-1}))$.
By induction, one obtains $(\ad_\gg X^\sharp)^k(u^\sharp) : a^n \mapsto
\frac{n!}{(n-k)!}(\ad a)^k u^\sharp(a^{n-k})$ for $n\geq k$
(for $n<k$ the RHS is assumed to be zero). 
Using this for $u^\sharp = h^\sharp$ defined by
some $h\in \gg$, we obtain 
$(\ad_\gg X^\sharp)^n (h^\sharp) : x^n \mapsto n! (\ad x)^n(h)$, 
and it is zero outside $S^n(\gg)$. By taking polarizations, we get
$(\ad_\gg X^\sharp)^n(h^\sharp) : x_1 x_2\ldots x_n \mapsto 
\sum_{\sigma \in S_n} (\ad x_{\sigma(1)})(\ad  x_{\sigma(2)})\ldots
(\ad  x_{\sigma(n)})(h).$

\nxpoint \label{p:mf.inv.coder}
Now we can write~\eqref{eq:dhxn} in another way. 
Recall that $\forall h\in \gg$ we have constructed a coderivation
$D_h : S(\gg)\to S(\gg)$ given by~\eqref{eq:dhxn}
such that $h\mapsto D_h$ is a Lie algebra embedding
$\gg\to \Coder(S(\gg))\cong \Hom(S(\gg),\gg)\cong\hat\calL_\gg$.
Recall that the map $\Coder(S(\gg))\to \Hom(S(\gg),\gg)$ maps
a coderivation $D$ into 
$\tilde{D}^\sharp:= X^\sharp\circ D : S(\gg)\to \gg$. Hence
$\tilde{D}_h^\sharp = X^\sharp\circ D_h$ is given by
$$
\tilde{D_h}^\sharp(x^n) = B_n\cdot (\ad x)^n(h).
$$
Comparing with~\refpt{p:sharpconc} we obtain the following equality 
in $\Hom(S(\gg),\gg)$:
$$
\tilde{D_h}^\sharp = \sum_{n\geq 0} 
\frac{B_n}{n!} (\ad_\gg X^\sharp)^n (h^\sharp)
$$
By dualizing and taking into account~\refpt{p:unpgequal} 
we obtain an equality in
$\hat\calL_\gg \cong \hat{S}(\gg^*)\otimes\gg \cong \gg_{(\hat{S}(\gg^*))}$:
$$
{}^t\tilde{D}_h = \sum_{n\geq 0} \frac{B_n}{n!} (\ad_\gg X)^n (h).
$$
This is exactly the main formula~(\ref{eq:mainfDversion}) 
of Section~\ref{s:EndOfProof} in invariant form. 
So we have proved it again in a shorter but less geometric way. 
The above formula has already appeared in a slightly different form in
\cite{Petracci}, Theorem 5.3 and formulas (20), (13) and (15).
We refer to the introduction to the present work for a more detailed
comparison of our results with those of {\em loc.cit.}

\nxpoint In the proof presented in this section, we started from an invariantly
defined isomorphism of coalgebras $\xi : C= S(\gg) \stackrel\sim\to H = \Uenv(\gg)$ and used it to transport the coderivations 
$L_h : x \mapsto hx$ from $H$ onto $C$.
However, we could replace in this reasoning $\xi$ by any other isomorphism of
coalgebras $\xi' : C\stackrel\sim\to H$ and
obtain another embedding $\gg \to \hat\calL_\gg$ in this way, hence another 
formula. Yet another possibility is to consider on $H$ the coderivations
$R_h: x\mapsto -xh$ instead of the $L_h$-s. 
This gives the same formula but with additional
$(-1)^n$ factors in each summand. (One can see this by considering the
isomorphism of $\gg$ onto its opposite $\gg^\circ$ given by
$x\mapsto -x$).

\nxpoint What are the other possible choices of $\xi'' : C\stackrel\sim\to H$ ?
Actually the only choice functorial in $\gg$ is the coexponential map. 
Thus we need some additional data. 
Suppose for example that $\gg$ is free rank $n$ as a $\genfd$-module,
and $e_1,\ldots, e_n$ are a base of $\gg$. Denote by $x_i$ the images of
$e_i$ in $C = S(\gg)$, and by $z_i$ the images
in $H = \Uenv(\gg)$. Then, one can take 
$\xi' : x_1^{\alpha_1} \ldots x_n^{\alpha_n}\mapsto  
z_1^{\alpha_1} \ldots z_n^{\alpha_n}$, for any $\alpha_i \in \bbN_0$.
This is easily seen to be a coalgebra isomorphism, since 
$\Delta_C(x_i) = x_i\otimes 1 + 1\otimes x_i$,
$\Delta_H(z_i) = z_i\otimes 1 + 1 \otimes z_i$,
and these monomials in $z_i$ form a base of $H$ by PBW theorem.
In this way, one obtains another embedding $\gg\hookrightarrow\hat\calL_\gg$,
the lower degree terms of which have to be given by the triangular matrices. 
In the geometric language of previous sections, this corresponds to the
map $\uW(\gg)\stackrel{\exp''}\rightarrow G = \Expt(\gg)$ given by
$x_1 e_1 +\ldots x_n e_n \mapsto
\exp'(x_1 e_1)\cdot \ldots \cdot\exp'(x_n e_n).$

\section{Conclusion and perspectives}

Within the general problematics of finding Weyl-algebra realizations
of finitely generated algebras, a remarkable universal formula
has been derived in three different approaches, 
suggesting further generalizations. 
Our result and proof can apparently also be extended, 
in a straightforward manner, to Lie superalgebras.

More difficult is to classify all
homomorphisms $U(L)\rightarrow A_n[[t]]$, 
which are not universal, but rather defined
for a given Lie algebra $L$. 
We have computed some examples of such representations 
(e.g.~\cite{MS}), but we do not know any classification results.
As usual for the deformation problems, we expect that 
the homological methods may be useful for the treatment of
concrete examples. 

In our representation, $t$ is a formal variable. If $\genfd$ is
a topological ring, then one can ask if our formal series 
actually converges for finite $t$. 
Let $\rho : A_n\to {\mathcal B}(H)$ be a representation of the Weyl
algebra by bounded operators on a Hilbert space $H$, and $t$ fixed. 
Then for any $x \in L$, 
the applying $\rho\circ \Phi_\lambda(x)$ is a power series
in bounded operators. Under the conditions when this is
a convergent series, for  
$\genfd = \mathbb C$, $\lambda = 1$, our formula is known, 
see Appendix 1, formula (1.28) of~\cite{KarMaslov}, 
with a very different proof.

Similarly to the analysis in~\cite{OdesFeigin}, it may be useful to
compute the commutant of the image of $\Phi_\lambda$ in $A_{n}[[t]]$.
It is an open problem if there are similar homomorphisms for the 
quantum enveloping algebras. The approach to formal Lie theory 
taken in~\cite{Holtkamp} may be useful in this regards. 
We expect that our approach may be adapted to the setup of Lie
theory over operads~\cite{Fresse}. In particular, generalizations to
Leibniz algebras (no antisymmetry!) would be very interesting.
There is also an integration theory for Lie algebroids (yielding
Lie groupoids), with very many applications. 
This suggests that the vector field computations
may be adaptable to that case. 

Our main motivation is, however, 
to explore in future similar representations 
in study of possible quantum field theories 
in the backgrounds given by noncommutative spaces, 
where we may benefit on unifying methods and intuition based 
on the exploration of the uniform setup of Weyl algebras. 
Some related physically inspired papers 
are~\cite{AC1,AC2,Luk,LukWor,MS}.

{\bf Acknowledgements.} 
S.M., A.S. and Z.\v{S}. were partly supported by Croatian Ministry Grant
0098003, and N.D. acknowledges partial support from the 
Russian Fund of Fundamental research, grant 04-01-00082a. 
Part of the work was done when Z.\v{S}. was a guest at MPIM Bonn.
He thanks for the excellent working conditions there.
Z.\v{S}. also thanks organizers of the conference {\it Noncommutative
algebras} at Muenster, February 2006, where the work was presented. 
We thank {\sc D. Svrtan} for discussions and references.


\begin{thebibliography}{99}


\bibitem[AmCam1]{AC1}
{\sc G.~Amelino-Camelia, M.~Arzano, L.~Doplicher}, 
Field theories on canonical and Lie-algebra noncommutative spacetimes,
in ``Florence 2001, A relativistic spacetime odyssey'' 
497-512, {\tt hep-th/0205047}

\bibitem[AmCam2]{AC2}
{\sc G.~Amelino-Camelia, M.~Arzano}, {\em Coproduct and star product 
in field theories on Lie-algebra non-commutative space-times}, 
Phys. Rev. D65:084044 (2002) {\tt hep-th/0105120},

\bibitem[Berceanu]{Berceanu}
{\sc S.~Berceanu}, {\em Realization of coherent state Lie algebras by 
differential operators}, {\tt math.DG/0504053}

\bibitem[Bourbaki]{Bourbaki}
{\sc N.~Bourbaki}, {\em Lie groups and algebras} (mainly Ch. I-II)

\bibitem[Dim]{Dim2004}
{\sc M. Dimitrijevi\'{c}, F. Meyer, L. M\"{o}ller, J. Wess},
{\em Gauge theories on the $\kappa$-Minkowski spacetime},
Eur. Phys. J. C Part. Fields  36  (2004),  no. 1, 117--126.

\bibitem[Fresse]{Fresse}
{\sc B. Fresse}, 
{\em Lie theory of formal groups over an operad},
J. Algebra 202 (1998), no. 2, 455--511. \MR{99c:14063}

\bibitem[SGA3]{SGA3} 
{\sc M. Demazure, A. Grothendieck et al.},
{\em Sch\'emas en groupes. I: Propri\'et\'es g\'en\'erales des sch\'emas 
en groupes} (SGA3, vol.1), Springer LNM 151.

\bibitem[Holtkamp]{Holtkamp}
{\sc R. Holtkamp}, 
{\em A pseudo-analyzer approach to formal group laws not of operad type},
J. Algebra 237 (2001), no. 1, 382--405. \MR{2002h:14074}

\bibitem[KarMaslov]{KarMaslov}
{\sc M. Karasev, V. Maslov}, {\it Nonlinear Poisson brackets},
Moskva, Nauka 1991 (in Russian); 
Engl. transl.: AMS Transl. Math. Monog. 119, 1993.

\bibitem[Kathotia]{Kathotia} {\sc V. Kathotia}, 
{\em Kontsevich's universal formula for 
deformation quantization  and the 
Campbell-Baker-Hausdorff formula},  
Internat. J. Math.  11  (2000),  no. 4, 523--551; {\tt math.QA/9811174}
\MR{2002h:53154}

\bibitem[Kontsevich]{KontsDefPoiss} {\sc M.~Kontsevich}, 
{\em  Deformation quantization of Poisson manifolds}, 
Lett. Math. Phys.  66  (2003),  no. 3, 157--216. \MR{2005i:53122}

\bibitem[Lukierski]{Luk} 
{\sc J. Lukierski, H.~Ruegg}, 
{\em Quantum $\kappa$-Poincare in Any Dimensions}, 
Phys. Lett. B329 (1994) 189-194, {\tt hepth/9310117} 

\bibitem[LukWor]{LukWor}
{\sc J.~Lukierski, M.~Woronowicz},
{\em New Lie-algebraic and quadratic deformations of Minkowski space 
from twisted Poincar\'{e} symmetries}, {\tt hep-th/0508083}

\bibitem[MS]{MS}
{\sc S.~Meljanac, M.~Stoji\'{c}}, {\em New realizations of Lie algebra
kappa-deformed Euclidean space}, preprint (2006)

\bibitem[OdesFeigin]{OdesFeigin}
{\sc A.~V. Odesskii, B.~L. Feigin}, {\em Quantized moduli spaces of
the bundles on the elliptic curve and their applications},  
Integrable structures of exactly solvable 
2d models of QFT (Kiev, 2000), 
123--137, NATO Sci. Ser. II Math. Phys. Chem., 35, 
Kluwer 
{\tt math.QA/9812059} \MR{2002j:14040} 

\bibitem[Petracci]{Petracci}
{\sc E. Petracci}, 
{\em Universal representations of Lie algebras by coderivations}
Bull. Sci. Math. 127 (2003), no. 5, 439--465; 
{\tt math.RT/0303020} \MR{2004f:17026}  
\end{thebibliography}
\end{document}